\documentclass[11pt,twoside,reqno]{amsart}
\usepackage{amsmath,amsfonts,amssymb,amsthm}
\usepackage{enumerate}
\usepackage[utf8]{inputenc}
\usepackage{hyperref}
\usepackage[usenames]{color}
\usepackage[english]{babel}
\usepackage{graphicx}
\usepackage[usenames]{color}
\usepackage{tikz-cd}
\usepackage{epstopdf}
\usepackage{subcaption}
\usepackage[margin=1cm]{geometry}
\usepackage{anysize}

\newcommand{\Rtres}{\mathbb{R}^3}
\newcommand{\Rdos}{\mathbb{R}^2}

\newcommand{\zhat}{\frac{\partial}{\partial z}}

\newcommand{\x }{\frac{\partial}{\partial x}}

\newcommand{\cu}{\mathcal{C}_U}

\newcommand{\y }{\frac{\partial}{\partial y}}

\newcommand{\z}{\frac{\partial}{\partial z}}

\newcommand{\rhoo}{\frac{\partial}{\partial \rho}}
\newcommand{\rhoohat}{\frac{\partial}{\partial \rho}}
\newcommand{\thetaa}{\frac{\partial}{\partial \theta}}

\newcommand{\Rcossin}{\R[\cos\theta, \sin \theta]}
\newcommand{\Rcossinz}{\R[\cos\theta, \sin \theta,z]}
\newcommand{\Rcossinzrho}{\R[\cos\theta, \sin \theta,z,\rho]}
\newcommand{\Sing}{\text{Sing}}
\newcommand{\Singad}{\wt{\text{Sing}}}
\newcommand{\fix}{\text{Fix}}
\newcommand{\ord}{\text{ord}}

\newenvironment{fact}{\vspace{1pt}\\ \noindent \textbf{Fact.}}{\\[1pt]}

\newtheorem{thm}{Theorem}[section]
\newtheorem{theorem}[thm]{Theorem}
\newtheorem{defi}[thm]{Definition}
\newtheorem{definition}[thm]{Definition}
\newtheorem{prop}[thm]{Proposition}
\newtheorem{proposition}[thm]{Proposition}
\newtheorem{cor}[thm]{Corollary}
\newtheorem{corollary}[thm]{Corollary}
\newtheorem{lem}[thm]{Lemma}
\newtheorem{lemma}[thm]{Lemma}
\newtheorem{rk}[thm]{Remark}

\theoremstyle{remark}
\newtheorem{remark}[thm]{Remark}
\newtheorem{notation}[thm]{Notation}
\newtheorem{example}[thm]{Example}

\newcommand{\R}{\mathbb{R}}

\newcommand{\SSS}{\mathbb{S}}

\newcommand{\CC}{\mathcal{C}}

\newcommand{\wt}[1]{{{\widetilde{#1}}}}

\newcommand{\g}{\gamma}
\newcommand{\G}{\Gamma}

\renewcommand{\wt}[1]{\widetilde{#1}}

\author[Nuria Corral, María Martín Vega, Fernando Sanz Sánchez]{Nuria Corral \\ \textit{Universidad de Cantabria, Spain} \\
	María Martín Vega\\ \textit{University of Valladolid, Spain}\\
	Fernando Sanz Sánchez\\ \textit{University of Valladolid, Spain}
}

\title[]{Surfaces with central configuration and Dulac's problem for a three dimensional isolated Hopf singularity}
\date{}

\begin{document}
	\maketitle
	\begin{abstract}
	Let $\xi$ be a real analytic vector field with an elementary isolated singularity at $0\in \Rtres$ and eigenvalues $\pm bi,c$ with $b,c\in \R$ and $b\neq 0$. We prove that all cycles of $\xi$ in a sufficiently small neighborhood of $0$, if they exist, are contained in a finite number of subanalytic invariant surfaces entirely composed by a continuum of cycles. In particular, we solve Dulac's problem, i.e. finiteness of limit cycles, for such vector fields.
	\end{abstract}
 	\section{Introduction and statements}
Dulac's problem is a central topic in the study of the dynamics of real analytic vector fields. It consists in proving that there are no infinitely many limit cycles accumulating and collapsing to a singular point. Recall that in general, a {\em cycle} (or a {\em closed orbit}) of a vector field 
in a given manifold $M$ is the image of a non-trivial periodic solution $\g:\R\to M$ (also denoted by $\g$), and a {\em limit cycle} is a cycle possessing a neighborhood free of other cycles. 

In dimension two,  
the problem was answered by Dulac in 1923 (\cite{Duc}), but his proof had an important gap. It was resolved nearly 70 years after by Ilyashenko \cite{Ily} and Écalle \cite{Eca}, with two independent and different proofs, both very intricate. There are recent attempts for alternative proofs, in some particular cases, using {\em o-minimal geometry}, for instance \cite{Kai-R-S,Spe,Gal-K-S,Dol-S}. 
Dulac's result has as a consequence that an analytic vector field in the sphere has finitely many limit cycles which, in turn, proves the finiteness statement of (the second part of) Hilbert's 16th problem. Namely, any polynomial vector field in $\R^2$ has finitely many limit cycles (see Ilyashenko's survey \cite{Ily2} for more information).

Non-accumulation of limit cycles for planar analytic vector fields implies a stronger property: 
either there are none in a neighborhood, or there is a continuum family of nested cycles filling a whole puncture neighborhood. In fact, given a cycle $\g$, we can define the {\em Poincaré first return} map in a transversal segment at some point $p\in\g$. It is an analytic local diffeomorphism at $p$ and its fixed points correspond exactly to cycles in a neighborhood of $\g$, thus there are necessarily finitely many or they form a continuum annulus around $\g$ (a {\em central configuration}). This is also the argument for proving Dulac's result in the easiest case in dimension two (apart, of course, from the trivial hyperbolic or semi-hyperbolic situations, 
when no local cycles exist). It is the case where the linear part of the vector field has purely imaginary non-zero eigenvalues
(the so-called {\em Hopf singularity}). 

\strut

In this paper, we solve Dulac's problem for analytic three-dimensional vector fields with isolated singularity with a pair of conjugated imaginary non-zero eigenvalues (a \emph{three-dimensional Hopf singularity}).
In fact, we determine a finite number of invariant surfaces where local cycles may be placed and theses surfaces present a central configuration. Let us provide precise statements.

Denote by $\mathfrak{X}^{\omega}(\R^3,0)$ the family of germs $\xi$ of analytic singular vector fields at the origin of $\R^3$, that is, $\xi(0)=0$. If $\xi\in\mathfrak{X}^{\omega}(\R^3,0)$ and $U$ is an open neighborhood of $0$ where (a representative of) $\xi$ is defined, we denote by $\mathcal{C}_U=\mathcal{C}_U(\xi)$ the union of all cycles of $\xi|_U$ (that is, entirely contained in $U$). It is called the {\em cycle-locus} of $\xi$ in $U$. Notice that this cycle-locus depends strongly on the neighborhood $U$ and that it does not behave as a germ of a set that we can associate to the germ $\xi$ (i.e., if $U'\subset U$ we can only assert that $\mathcal{C}_{U'}\subset\mathcal{C}_U$, but not $\mathcal{C}_{U'}=U'\cap\mathcal{C}_U$). 

Consider the following family:
$$
\mathcal{H}^3:=\{\xi\in\mathfrak{X}^{\omega}(\R^3,0)\,: \ {\rm Spec}(D\xi(0))=\{\pm bi,c\},\mbox{ where }b,c\in\R\mbox{ and }b\ne 0\}.
$$
Observe that any $\xi\in\mathcal{H}^3$ has a unique formal invariant curve $\widehat{\Omega}=\widehat{\Omega}_\xi$ at $0$, which is non-singular and tangent to the eigenspace corresponding to the  eigenvalue $c$ (see \cite{Bon-D}). It is called the {\em (formal) rotational axis} of $\xi$. When $c\neq 0$, the rotational axis is convergent and provides an analytic invariant curve, since in this case $\widehat{\Omega}$ coincides with the stable or unstable manifold of $\xi$ (see for instance~\cite{Car-S}). On the contrary, if $c=0$ (the \emph{completely hyperbolic case} or \emph{zero-Hopf singularity}), then $\widehat{\Omega}$ may be convergent or not, although there is always an invariant $\mathcal{C}^\infty$-curve whose Taylor expansion at $0$ coincides with that of the formal one. This is a result by Bonckaert and Dumortier in \cite{Bon-D} in the case where $\xi$ has an isolated singularity, and trivially true otherwise since, in this case, $\widehat{\Omega}$ coincides with the singular locus $\Sing (\xi)$, an analytic curve. Notice that in the semi-hyperbolic case $\xi$ has isolated singularity.

\strut 

The main result in this paper can be stated as follows.
\begin{theorem}[Structure of cycle-locus]\label{th:main}
	Let $\xi\in\mathcal{H}^3$ with isolated singularity. Then, there is some neighborhood $U$ of $0\in \Rtres$ where a representative of $\xi$ is defined for which exactly one of the following possibilities holds:
	\begin{enumerate}[(i)]
		\item $\mathcal{C}_U(\xi)=\emptyset$.
		\item There is a finite non-empty family $\mathcal{S}= \{  S_1,...,S_r \}$ of connected regular analytic two-dimensional submanifolds of $U\setminus\{0\}$, mutually disjoint, invariant for $\xi$, subanalytic sets in $U$ and satisfying  $\overline{S_j}\cap U=S_j\cup\{0\}$ for any $j$, such that, for any element $V\subset U$ in some neighborhood basis at $0$, we have
		$$
		\mathcal{C}_{V}(\xi)=(S_1\cup S_2\cup\cdots\cup S_r)\cap V.
		$$
	\end{enumerate} 

\end{theorem}

	\begin{figure}[h]
	\centering

		\centering
	
		\includegraphics[width=0.35\linewidth]{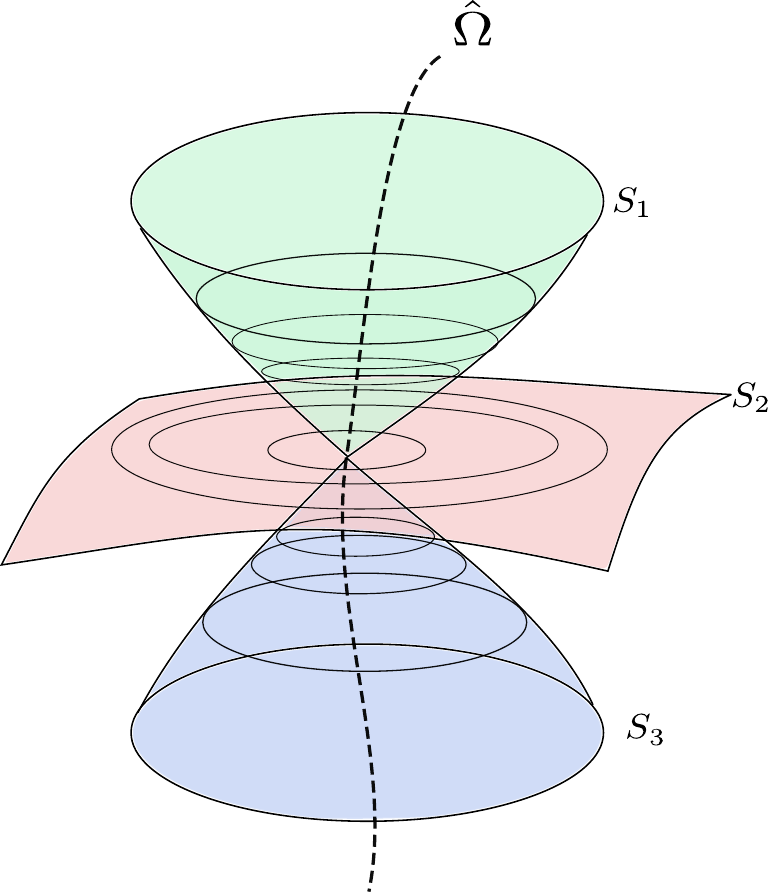}
		\caption{Illustration of case $(ii)$.}
		\label{fig:variassuperficies}

\end{figure}

As a consequence, Dulac's property is true for these vector fields:
\begin{corollary}\label{cor:Dulac}
	Let $\xi\in\mathcal{H}^3$ with an isolated singularity. Then, there are not infinitely many limit cycles of $\xi$ accumulating and collapsing to $0\in\R^3$.
\end{corollary}
In the second possibility \textit{(ii)}, the germs of the surfaces $S_j\in \mathcal S$ are uniquely determined. Each of them, in a sufficiently small neighborhood, is composed by a continuum of nested cycles around the singularity; i.e., each $S_j$ is a surface with a central configuration as in the planar case (although $S_j$ could be singular at the origin). Let us call them {\em limit central surfaces} by analogy with the concept of limit cycle. The following example defines two limit central surfaces, both being singular at the origin.

\begin{example}\label{example:introduction}
	Consider the following vector field in $\mathcal{H}^3$.
	\begin{equation*}
		\xi= (-y-xz^2+x(x^2+y^2))\x+ (x-yz^2+y(x^2+y^2))\y+ (z^3+z(x^2+y^2))\z.
	\end{equation*}
	It has isolated singularity. The two half-cones $S_1=\{ (x,y,z): x^2+y^2-z^2=0, z>0  \}$ and $S_2=\{ (x,y,z): x^2+y^2-z^2=0, z<0  \}$ are invariant. The restriction of $\xi$ to any of the surfaces $S_i$ is $\xi\lvert_{S_i}= -y\x+x\y$, which proves that $\xi$ defines a central configuration in $S_i$, for $i=1,2$.
\end{example}

The result stated in Theorem~\ref{th:main} for the semi-hyperbolic case ($c\ne 0$) has been already proved by Aulbach \cite{Aul}, in a more general situation of $n$-dimensional analytic vector fields with a pair of purely imaginary non-zero eigenvalues and $n-2$ non-zero real eigenvalues. Before Aulbach, the same situation has been considered in the literature by other authors (\cite{Kel,Kir,Schn,Schm}), under the assumption that the vector field has a first integral, as in the classical Lyapunov's result \cite{Lya}. It is worth to comment that, using a basic property from the theory of center manifolds, one obtains that the possibility (ii) can only occur for a unique limit central surface ($r=1$), which coincides with the center manifold $W^c$ of $\xi$ (hence unique, non-singular and analytic in this case).   

Vector fields with a Hopf singularity in the completely non-hyperbolic case ($c=0$) have been studied in the literature. For instance, Dumortier in \cite{Dum} considered such germs  of vector fields of class
$\mathcal{C}^\infty$ at $0\in\R^3$ satisfying two {\em \L ojasiewicz-type} inequalities: one for the vector field itself, which implies that $0$ is an isolated singularity; and a second one for the {\em infinitesimal generator} of the Poincaré first-return map 
associated to the central cycle that appears after the blowing-up of an invariant $\mathcal{C}^\infty$ realization of $\widehat{\Omega}$. He obtains a complete description of the asymptotic behavior of all trajectories in a neighborhood of the origin, as well as a {\em weak topological} classification of the vector field. However, those assumptions prevent the existence of any local cycle (that is, one has only the possibility (i) of Theorem~\ref{th:main}). In our situation where $\xi$ is analytic, \L ojasiewicz's inequality for $\xi$ is equivalent to the property of isolated singularity, but we do not require the second assumption, thus permitting the existence of cycles and hence the possibility (ii).  

We should mention that vector fields in $\mathcal{H}^3$ have also been considered in families for different purposes. We can mention Guckenheimer and Holmes \cite{guckenheimer2013nonlinear}, where there is a complete description of the bifurcation diagrams for small codimension singularities; Baldoma, Ibañez and Martinez-Seara \cite{Bal-Ib-2020} where the appearance of certain chaotic behavior, associated to a \emph{Shilnikov configuration}, is studied; García \cite{IsaacCyclicity}, where it is stablished a bound for the number of limit cycles appearing in certain generic families inside $\mathcal{H}^3$ with a bounded number of turns around the rotational axis.

\strut

Let us summarize the ideas for the proof of Theorem~\ref{th:main} and the plan of the article. 

In section~2, we propose a simple proof in the semi-hyperbolic case, in spite of the existing references already mentioned for this situation. Our aim, apart for the sake of completeness, is to introduce some of the arguments involved in the proof of the general case, absent in Aulbach's proof \cite{Aul} but appearing in Dumortier's work \cite{Dum}. Namely, blowing-up techniques and Poincaré first return map along the cycles emerging from the blowing-up.

The rest of the article is devoted to the proof in the completely non-hyperbolic case. 
We work with formal normal forms $\hat{\xi}$ for $\xi$ (for instance in the sense of Takens \cite{Takens}) and the sequence of analytic vector fields $\{\xi_\ell\}_\ell$ that approximate $\hat \xi$. That is, in terms of jets, they satisfy $j_\ell(\xi_\ell)=j_\ell(\hat \xi)$. Being these vector fields analytically conjugated to the original vector field $\xi$, it is enough to prove Theorem~\ref{th:main} for some $\xi_{\ell}$.

In section~3, we use blowing-ups to study $\hat \xi$ and its jets approximations $\xi_\ell$. The rotational axis $\widehat{\Omega}$ is not necessarily convergent and we cannot blow it up (or any realization of it) if we want to preserve analyticity. Starting from the blowing-up of the origin, we define recursively \textit{sequences of admissible blowing-ups}: a composition of blowing-ups 
centered at either the infinitely close points of $\widehat{\Omega}$ (\textit{characteristic singularities) }
or invariant closed circles of the corresponding strict transforms of $\hat \xi$ (\textit{characteristic cycles}).
The main result of this section is a\textit{ reduction of singularities} of the normal form $\hat{\xi}$ \emph{adapted} to our problem.
This process may be understood as a refinement, for this situation, of Panazzolo's result on reduction of singularties of general three-dimensional analytic vector fields~\cite{panazzolodim3} (notice that a Hopf singularity is already in the final elementary situation in the sense of Panazzolo). Essentially, the formal normal form $\hat \xi$ can be viewed as a vector field of revolution by rotating a planar vector field $\hat{\eta}$; the adapted reduction of singularities corresponds to the reduction of singularities of $\hat \eta$.
 In a second part of section 3, we discuss how to apply sequences of admissible blowing-ups to the jet approximations $\xi_{\ell}.$ In particular, we find a lower bound for $\ell$ so that the characteristic cycles are actual cycles of the (strict transform of) $\xi_\ell$, as well as other dynamical properties that depend on a finite jet are inherited by $\xi_\ell$.

In section~4, we prove that, after any sequence of admissible blowing-ups, the characteristic cycles and the characteristic singularities are the only possible limit sets of families of cycles of the transform of $\xi_{\ell}$ in a neighborhood of the singularity, provided that $\ell$ is large enough. Thus, in order to prove Theorem~\ref{th:main}, we only search for cycles near the characteristic cycles and characteristic singularities.

In section~5, we study the different local situations present after an adapted reduction of singularities $\pi:(M,E)\longrightarrow (\Rtres ,0)$ of $\hat \xi$. We have specific monotonic functions along the trajectories of the transformed vector field $\widetilde{\xi}_\ell= \pi^\ast \xi_\ell$ in neighborhoods of characteristic singularities or corner-characteristic cycles, preventing the existence of cycles of $ \widetilde{\xi}_\ell$ in sufficiently small neighborhoods of them. Around a non-corner characteristic cycle $\gamma$, we work with the associated Poincaré first return map $P_\gamma$ of $\widetilde{\xi}_\ell$. First, we find a formal invariant non-singular surface $S_\gamma$ for $\widetilde{\xi}_\ell$ supported by $\gamma$ and transversal to the divisor using that this is the case for the transform $\pi^\ast \hat \xi$ of $\hat \xi$. This surface $S_\gamma$ provides a formal invariant curve $\Gamma_\gamma$ for $P_\gamma$ and around $\Gamma_\gamma$ we can describe the periodic orbits of $P_\gamma$. Namely, there is a conic neighborhood $\Sigma_\gamma$ around $\Gamma_\gamma$ such that:
 If $\Gamma_\gamma \nsubseteq \text{Fix}(P_\gamma)$, there are not periodic points of $P_\gamma$ inside $\Sigma_\gamma$. If, otherwise, $\Gamma_\gamma \subseteq \text{Fix}(P_\gamma)$, then $\Gamma_\gamma$ is exactly the set of periodic points (thus fixed) inside $\Sigma_\gamma$.

Finally in section~6, we give the proof of Theorem~\ref{th:main} gathering the results of the previous sections.
First, we make a reduction of singularities and fix a vector field $\xi_{\ell}$ to which the sequence of admissible blowing-ups can be applied.
By means of the results in sections~4~and~5, cycles of $\widetilde \xi_{\ell}$ near but not contained in the divisor $E$ (and hence of $\xi_{\ell}$ in a neighborhood of $0\in \Rtres$ under the isomorphism $\pi:M\setminus E \to \Rtres\setminus \{ 0 \}$) are located in neighborhoods of non-corner characteristic cycles $\gamma$. Moreover, the conic neighborhoods $\Sigma_\gamma$ above provide solid conic neighborhoods $\widetilde \Sigma_\gamma$ of $S_\gamma$ in such a way that if a cycle of $\widetilde{\xi}_\gamma$ is inside $\widetilde \Sigma_\gamma$, then the curve $\Gamma_\gamma$ is contained in $\fix (P_\gamma)$ and surpports a continuum of cycles inside an analytic surface around $S_\gamma$. The projection of it under $\pi$ provides a limit central surface. This would finish the proof of Theorem~\ref{th:main} if we could guarantee that there are not other cycles of $\widetilde{\xi}_\ell$ outside $\widetilde \Sigma_\gamma$. We can ``open" these cones to actual neighborhoods by means of further blowing-ups. But then, it is possible that we need a larger jet approximation $\xi_{\ell'}$ with $\ell'\geq \ell$ for which the order of these cones could change, a priori. We overcome this difficulty showing that the order of the cones where the cycles have the desired properties may be uniformly bounded for $\ell'\geq \ell$.

\section{The semi-hyperbolic case}
In this section, we provide a proof of Theorem~\ref{th:main} in the semi-hyperbolic case; i.e., the linear part $D\xi(0)$ has eigenvalues $\{\pm bi,c\}$ with both $b,c$ different from zero.

Assume for instance that $c<0$. Then, the stable manifold $W^s
$ of $\xi$ at $0$ is one-dimensional and, being analytic, it coincides with the rotational axis, which is therefore convergent.  
Fix some center manifold $W^c$ of $\xi$ at $0$ of class $\CC^k$, with $k\ge 2$. In general, it is not analytic, nor unique. But it contains any cycle of $\xi$ that is contained in a sufficiently small neighborhood $U$ of the origin; i.e., $\CC_U(\xi)\subset W^c$ (see \cite{Carr}). 

Take a neighborhood $U_0$ inside which, both the stable manifold $W^s$ and the chosen center manifold $W^c$ are regular embedded submanifolds, and such that $\CC_{U_0}(\xi)\subset W^c$. Let $\pi:M\to U_0$ be the polar blowing-up with center $W^s$. It is a proper analytic map. The divisor $E=\pi^{-1}(W^s)$ is a cylinder and the fiber $\g=\pi^{-1}(0)$ over the origin is a cycle of the transformed vector field $\wt{\xi}:=\pi^*\xi$. The strict transform $\wt{W^c}=\overline{\pi^{-1}(W^c\setminus\{0\})}$ is a surface of class $\CC^{k-1}$, invariant for $\wt{\xi}$ and transversal to $E$. Moreover, $\g=E\cap\wt{W^c}$. 

\begin{figure}[h]
	\centering
	\includegraphics[width=0.97\linewidth]{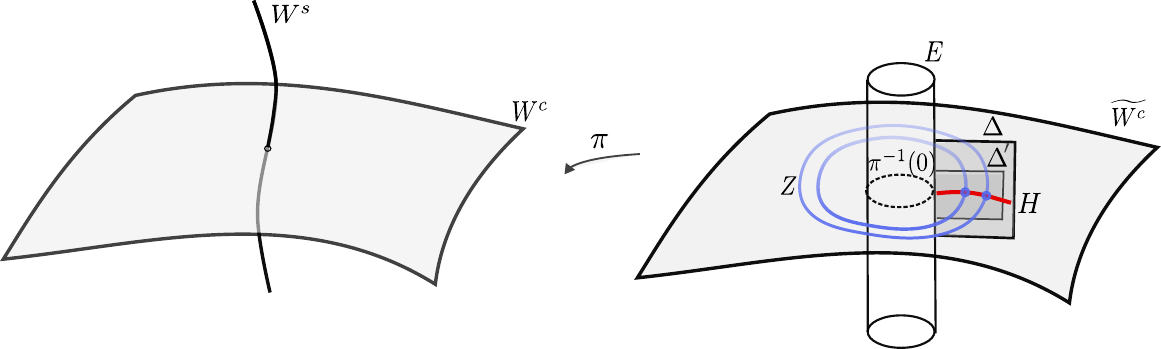}
	\caption{Definition of the Poincaré map $P_\gamma$}
	\label{fig:poincaresemihip}
\end{figure}

Now, consider a point $a\in\g$, and two nested analytic discs $\Delta'\subset\Delta$  transverse to $\wt{\xi}$ at $a$ so that the Poincaré first-return map $P_\g:\Delta'\to\Delta$ of $\wt{\xi}$ associated to $\g$ is well defined and analytic. Notice that if $Z$ is any cycle of $\wt{\xi}$ such that $Z\cap \Delta=Z\cap \Delta'$, 
then the intersection $Z\cap\Delta$ is a periodic orbit of $P_\g$  (see Figure \ref{fig:poincaresemihip}). In particular, if $Z$ is the inverse image by $\pi$ of a cycle inside $\CC_{U_0}(\wt{\xi})$, then, $Z$ is contained in $\wt{W^c}$. Taking into account that $W^c$ is two-dimensional and using classical arguments based on the Jordan Curve Theorem (see for instance \cite{Pal-dM}), we conclude in this case that $Z$ cuts $\Delta'$ in a single point, necessarily a fixed point of $P_\g$. 
Hence, the family of cycles of $\xi$ in a given neighborhood of the origin are in bijection with the set of fixed points of $P_\g$ not in $E$, and this set is contained in the intersection $H=\wt{W^c}\cap\Delta'$. Let us prove now Theorem~\ref{th:main}.

Suppose that item \textit{(i)} does not hold; i.e.,  $\CC_U( \xi)\ne\emptyset$ for any open neighborhood $U$ of $0$. Then we have infinitely many cycles of $\xi$ that accumulate and collapse to $0$. By the above, there are infinitely many fixed points of $P_\g$ in $H$ accumulating to the point $a$. Being $P_\g$ an analytic map, its set of fixed points, $\text{Fix}(P_\g)$, is an analytic set of positive dimension. Since $\text{Fix}(P_\g)\subset H$ and $H$ is a curve of class $\CC^{k-1}$ (transversal intersection of $\wt{W^c}$ and $\Delta'$), we conclude that $H=\text{Fix}(P_\g)$ is an analytic curve.

Let $\wt{U}$ be a neighborhood of $\g$ in $M$ satisfying:
\begin{itemize}
	\item $\widetilde U\cap \Delta=\Delta'$.
	\item $\widetilde U\cap W^c$ is the saturation of $H$ by the flow of $\wt{\xi}$.
	\item $U=\pi (\wt{U})$ is contained in $U_0$.
\end{itemize}
  We get that $U$ is a neighborhood of $0$ and $\CC_{U}(\xi)=W^c\cap U\setminus\{0\}$. Notice also that $\wt{W^c}\cap\wt{U}$ is an analytic set since $H$ is an analytic curve. Since $\pi$ is proper, we conclude that $W^c\cap U$ is a subanalytic set and Theorem~\ref{th:main} is proved. 

\begin{rk}{\em 
		The proof above shows that, in the semi-hyperbolic case, there is at most one limit central surface $S_1$. 
		Moreover, if $S_1$ exists, then $\overline{S_1}=W^c$ is a center manifold which is moreover unique and analytic (using Tamm's Theorem \cite{tamm}, because $W^c$ is  of class $\CC^k$ and subanalytic in this case).
	}
\end{rk}

\section{
Admissible blowing-ups and adapted reduction of singularities
}\label{sec:adaptedres}
	Consider a vector field $\xi$ in the  family $\mathcal{H}^3$ with completely non-hyperbolic linear part, that is, 
	$\text{Spec}(\xi)=\{  \pm bi,0 \}$. 
	In some coordinates the vector field is written as
\begin{equation} \label{eqperturb}
\xi= \left( -by+ A_1(x,y,z)   \right)\x + \left( bx+ A_2(x,y,z)   \right)\y + \left(  A_3(x,y,z)   \right)\z,
\end{equation}
with $A_1,A_2,A_3 \in \mathbb{R}\{x,y,z\}$ of order at least two. 
Without loss of generality for the study of the foliation generated by $\xi$ we will assume $b=1.$

\subsection{Formal normal form and truncated normal forms} \label{sec:formalnormal}

	Using Takens' theorem on normal forms (see \cite{Takens}), 
	there exists a formal automorphism 
	at 0, expressed in terms of the chosen coordinates as 
	$$
	\hat{\varphi}(x,y,z)=(\hat{\varphi}_1(x,y,z),\hat{\varphi}_2(x,y,z), \hat{\varphi}_3(x,y,z)) \in \R[[x,y,z]]^3,$$ with $ D\hat{\varphi}(\underline{0})=Id$, such that the formal vector field $\hat{\xi}=\hat{\varphi}^\ast (\xi)$ is written in the form
	\begin{equation}\label{eq:formalVF}
			\hat{\xi }=
			T(x^2+y^2,z) \cdot \left( -y\x + x\y   \right) +  R(x^2+y^2,z) \cdot \left( x\x +y \y \right)\\
			+ Z(x^2+y^2,z)\z ,
	\end{equation}
	where $R,T,Z\in \mathbb{R}[[u,v]]$ and $T(0,0)=1$. Note that $R(u,v), Z(u,v) \in (u,v)$ and $Z(0,v) \in (v^2)$.
Neither the automorphism $\hat{\varphi}$ need to be convergent, nor the components of $\hat{\xi}$ need to belong to $\mathbb{R}\{ x,y, z  \}$. 
	Any formal vector field $\hat \xi$ as in \eqref{eq:formalVF} obtained as above is called a formal normal form of $\xi$. We remark that $\hat \xi$ is not uniquely determined by $\xi$.
	\begin{remark}
		\label{rmk:nondeg} \label{rmk:line-sing}
		The $z-$axis is sent to the rotational axis $\widehat{\Omega}$ of $ \xi$ by $\hat \varphi$, that is, $\widehat{\Omega}= \hat \varphi (0,0,z)$. On the other hand, since $\hat \varphi$ must preserve the (formal) singular locus, the hypothesis that $\xi$ has isolated singularity implies that $Z(0,v)\neq 0$. 
	\end{remark}

 Along this section, we will work with $k-$jets of the coefficients of vector fields in some coordinates. We remind some well-known properties and definitions just to fix notations. 
 In general terms, if $K$ is a $\R$-algebra and $K[[x_1,\ldots, x_n]]$  is the $\R$-algebra of formal power series with coefficients in $K$, an element $F\in K[[x_1,\ldots, x_n]]$ can be written as
 $$F(x_1,\ldots,x_n)=\displaystyle \sum_{i=0}^{\infty} F_i,$$ where each $F_i\in K[x_1,\ldots, x_n]$ is an homogeneous polynomial of degree $i$. When $F\neq 0$ the first index such that $F_i\neq 0$ is called the order of $F$ denoted by $\nu(F)$. For any $k\geq 0,$ we define the \emph{$k$-jet} of $F$ with respect to $\mathbf{x}=(x_1,\ldots, x_n)$ as $$j_k(F)=j_k^{\mathbf{x}}(F):= \displaystyle \sum_{i=\nu(F)}^{k} F_i  .$$
 
  Suppose that $F$ can be seen as an element $\tilde F\in K[x_1,\ldots,x_{s}][[ x_{s+1}, \ldots, x_n]] \subset K[[x_1,\ldots, x_n]]$, that is, $F$ is the image of some $\tilde F$ by the natural monomorphism $K[x_1,\ldots,x_{s}][[ x_{s+1}, \ldots, x_n]] \hookrightarrow K[[x_1,\ldots, x_n]]$. Denote $\mathbf{y}=(x_1,\ldots,x_{s})$, $\mathbf{z}= (x_{s+1},\ldots,x_{n})$, and $\tilde F(x_1,\ldots,x_n)=F(\mathbf{y}, \mathbf{z})=\displaystyle \sum_{i=0}^{\infty} \sum_{|\alpha|=i}\tilde F_\alpha(\mathbf{y}) \mathbf{z}^\alpha$ where $\mathbf{z}^\alpha=x_{s+1}^{\alpha_{s+1}}\cdots x_{n}^{\alpha_n} $, $|\alpha|=\alpha_{s+1}+ \cdots + \alpha_n$ and $\tilde F_\alpha\in K[\mathbf{y}]$. Then, the 
\emph{$k$-jet of $F$ with respect to $\mathbf{z}$} is $$j_k^{\mathbf{z}}(F):=j_k(\tilde F)= \displaystyle \sum_{i=0}^{k} \sum_{|\alpha|=i}\tilde F_\alpha(\mathbf{y}) \mathbf{z}^\alpha.$$

 In what follows, when this situation occurs, we will use the same notation for $F$ and $\tilde F$, the distinction will be clear from the context. We have the following properties:
\begin{itemize}
	\item $j_k(F\cdot G)= j_k( j_k(F) \cdot  j_k(G)  )$, for $F,G\in K[[x_1,\ldots, x_n]]$. In fact, this property can be refined:
	if $k\geq \max\{ \nu(F), \nu(G)  \},$ then $j_k(F\cdot G)= j_k (j_{k-\nu(G)}(F) \cdot  j_{k-\nu(F)}(G))$.
	\item $j_k(F^{-1})= j_k( (j_k (F))^{-1}  )$ if $F$ is a unit in $K[[x_1,\ldots, x_n]]$.
\end{itemize}
If $F\in K[\mathbf{y}][[\mathbf{z}]]$:
\begin{itemize}
	\item $j_k^{\mathbf{z}}(F(x_1,\ldots, x_{i}+a, \ldots, x_n))= j_k^{\mathbf{z}} (F)(x_1,\ldots, x_{i}+a, \ldots, x_n)$ for $i\leq s$ and $a\in K$.
	\item $j_k(F)=j_k( j_k^{\mathbf{z}}(F))$.
\end{itemize}
	If $\hat \eta$ is a formal vector field written in coordinates $\mathbf{x}$ as:
		$$ \hat \eta= \eta_1 \frac{\partial}{\partial x_1}+ \cdots + \eta_n \frac{\partial}{\partial x_n} ,$$
		where each $\eta_i\in K[[x_1,\ldots, x_n]]$, then the $k-$jet of $\hat \eta$ is the vector field
		$$ j_k(\hat \eta)= j_k(\eta_1) \frac{\partial}{\partial x_1}+ \cdots + j_k(\eta_n) \frac{\partial}{\partial x_n} .$$
	When each coefficient $\eta_j$ belongs to $K[\mathbf{y}][[\mathbf{z}]]$, we can consider the $k-$jet with respect to the coordinate $\mathbf{z}$ defined by 
	$$ j_k^{\mathbf{z}}(\hat \eta)= j_k^{\mathbf{z}}(\eta_1) \frac{\partial}{\partial x_1}+ \cdots + j_k^{\mathbf{z}}(\eta_n) \frac{\partial}{\partial x_n} .$$
	Once we fix a formal normal form $\hat \xi$ of $\xi$ given by $\hat \xi= \hat \varphi^\ast \xi$,
	we can consider truncated normal forms of $\xi$ in the following way. For any $\ell\in \mathbb{N}_{\geq 2}$, let $\varphi_\ell$ be the polynomial tangent to the identity diffeomorphism of $(\mathbb{R}^3,0)$ given by $$\varphi_\ell(x,y,z)=(j_{\ell+1} \hat{\varphi})(x,y,z):=  (j_{\ell+1}(x\circ \hat \varphi), j_{\ell+1}(y\circ \hat \varphi), j_{\ell+1}(z\circ \hat \varphi)).$$ 
	The vector field $\xi_\ell=(\varphi_\ell)^\ast(\xi)$ has the same $\ell-$ jet as the formal one $\hat{\xi}$ in coordinates $(x,y,z)$. That is, $j_\ell(\xi_\ell)=j_\ell(\hat{\xi})$.
	Notice that the vector field $\xi_{\ell}$ for any $\ell$, is analytically conjugated to $\xi$ and formally conjugated to $\hat{\xi}$. More precisely, we have the following formal equation:
	\begin{equation}\label{eq:conjugationpsiell}
	\hat \xi = \psi_\ell^\ast \xi_{\ell}  ,\ \text{where }\psi_\ell:= \varphi_\ell^{-1}\circ \hat{\varphi},
	\end{equation}
	It is sufficient to prove Theorem \ref{th:main} for $\xi_\ell$ for any $\ell$. The strategy is the following: we use $\hat{\xi}$ as a guiding vector field so that, after a sequence of blowing-ups, we get a transform of $\hat{\xi}$ with a specific good expression. Both the choice of the sequence of blowing-ups and the expression of the transform will depend only on a finite jet of $\hat{\xi}$, allowing us to choose $\ell$ sufficiently large so that all the construction is applied to $\xi_\ell$.
	
	Notably, we propose an \emph{adapted reduction of singularities of $\hat{\xi}$} and their corresponding approximations for $\xi_\ell$ with $\ell$ large enough. The blowing-ups will be real ones (thus generating boundary and corners) either with center at a point or at an analytic curve isomorphic to the circle $\SSS^1$. See for instance the work of Martín-Villaverde, Rolin and Sanz \cite{MartinVillaverdeRolinSanz} for intrinsic and general definitions of real blowing-ups.
	
	\subsection{The first blowing-up}\label{sec:firstblowing-up}
	 The first blowing-up will be the real blowing-up $\sigma_0: (M_0,E_0)\longrightarrow (\mathbb{R}^3,0)$ with center at the origin. The blown-up space $M_0$ is a manifold having the divisor $E_0= \sigma_0^{-1}(0)$ as its boundary. This divisor is homeomorphic to a sphere and represents the space of all the half-lines through $0$. The morphism $\sigma_0$ defines an analytic diffeomorphism from $M_0\setminus E_0$ to $\Rtres \setminus\{0  \}$.
	We consider $M_0$ covered by three charts $(C_0 , (\theta,z^{(0)},\rho^{(0)}))$, $(C_\infty ,(x^{(\infty)},y^{(\infty)},z^{(\infty)}))$ and $(C_{-\infty},(x^{(-\infty)},y^{(-\infty)},z^{(-\infty)}))$ where $C_0\simeq \SSS^1 \times \R \times \R_{\geq 0}$ and $C_{\pm\infty}\simeq \Rdos \times \R_{\geq 0}$. In these charts, the expression of $\sigma_0$ is given by:
	\begin{eqnarray} \label{eq:cartaC0}
	&	\text{In }C_0:  & \left\{ \begin{matrix}
			x &=& \rho^{(0)} \cos \theta\\
			y &=& \rho^{(0)} \sin \theta \\
			z &=& \rho^{(0)} z^{(0)}
		\end{matrix} \right.  \qquad  (\cos \theta, \sin \theta)\in \mathbb{S}^1, z^{(0)}\in \mathbb{R}, \rho^{(0)} \geq 0 \\
	\label{eq:cartaCinfty}
		&	\text{In }C_\infty:& \left\{ \begin{matrix}
			x &=& x^{(\infty)} z^{(\infty)}\\
			y &=& y^{(\infty)} z^{(\infty)} \\
			z &=& \space z^{(\infty)}
		\end{matrix} \right.\qquad  x^{(\infty)},y^{(\infty)} \in \mathbb{R}, \space z^{(\infty)}\geq 0
		\\ 	\label{eq:cartamenosCinfty}
	&	\text{In }C_{-\infty}:& \left\{ \begin{matrix}
			x &=& x^{(-\infty)} z^{(-\infty)}\\
			y &=& y^{(-\infty)}z^{(-\infty)}\\
			z &=& -z^{(-\infty)}
		\end{matrix}
		 \qquad
		x^{(-\infty)},y^{(-\infty)}\in \mathbb{R},
		z^{(-\infty)}\geq 0.
		\right.
	\end{eqnarray}
	\begin{remark}\label{rem:theta}
		Strictly speaking, $C_0$ is not the domain of a usual chart of $M_0$, since it is not homeomorphic to an open set of $ \Rdos\times\R_{\geq 0}$. Considering the usual covering $\tilde C_0=\Rdos\times\R_{\geq 0}$ with $\tau:\tilde C_0\longrightarrow C_0$ given by $( \theta, z, \rho)\mapsto(\sin \theta, \cos \theta, z, \rho)$, we can treat $\theta$ as a true coordinate (and we will tacitly do), so that $\sigma_0\circ \tau$ has the expression in \eqref{eq:cartaC0}. This convention justifies our abuse of terminology in the expression of the form ``a chart $(C_0,(\theta, z^{(0)},\rho^{(0)}))$".
	\end{remark}

	The origins of the charts $C_\infty$ and $C_{-\infty}$ will be denoted by $\gamma_\infty$ and $\gamma_{-\infty}$, respectively. 
	They are the points of the divisor $E_0$ corresponding to the half-lines contained in the $z-$axis and they are the only points of $E_0$ not covered by $C_0$.
	More explicitly, $\sigma_0(C_0)=\Rtres \setminus \{  x=y=0 \}$.

	 We define the \emph{(total) transform of $\hat \xi$ by $\sigma_0$ in the chart $C_0$} as the pull-back $$  \hat{\xi}
	^{(0)} := (\sigma_0\lvert_{C_0})^\ast \hat{\xi} . $$
	Using simplified notation, $(z,\rho):= (\rho^{(0)}, z^{(0)})$ and equations \eqref{eq:formalVF} and \eqref{eq:cartaC0}, the vector field $\hat{\xi}
	^{(0)}$ is given by:
	\begin{equation}
		\label{eq:campoprimeraexpformal}
		\hat{\xi}
		^{(0)} 
		=  B_\theta(z,\rho) \thetaa\\
		+ B_z(z,\rho) \z + B_\rho(z,\rho)
		\rhoo ,
\end{equation}
where $B_\theta(z,\rho)= T(\rho^2,\rho z)$, $B_z(z,\rho)= \frac{1}{\rho}Z(\rho^2,\rho z) -zR(\rho^2,\rho z)$ and $B_\rho(z,\rho)=\rho R(\rho^2,\rho z)$. 
Notice that $ B_\theta, B_z,B_\rho\in \mathbb{R}[z][[\rho]]$, $B_\theta(0,0)=1$ and that $( B_z,B_\rho)\neq (0,0)$ since $Z(u,v)\neq 0$ by Remark~\ref{rmk:nondeg}. Notice also that $\rho$ divides $B_z,B_\rho$.

The fact that the coefficient
$B_\theta(z,\rho)$ is a unit in $\mathbb{R}[z][[\rho]]$ allows us to consider $\theta$ as the ``time variable" and, consequently, $\hat{\xi}^{(0)}$ is completely described by the \emph{associated two dimensional} formal vector field $\hat{\eta}_0$ given by the system of formal ODEs:
\begin{equation}
\label{eq:formal2dim}
  \hat\eta_0: \left\{ \begin{array}{lcll}
  		\frac{d z}{d \theta} &=& B_\theta(z,\rho)^{-1}B_z(z,\rho)  &=\rho^{n^{(0)}}
  	A_z(z,\rho)\\
	\frac{d \rho}{d \theta} &= &B_\theta(z,\rho)^{-1}B_\rho(z,\rho)  &=\rho^{n^{(0)}} 
	A_\rho(z,\rho)
.
\end{array} \right. 
\end{equation}
In this expression, $A_i\in \mathbb{R}[z][[\rho]]$ for $i=z,\rho$ and ${n^{(0)}}$ is the maximum exponent $n$ such that $\rho^n$ divides both $B_\rho$ and $B_z$.
The\emph{ associated reduced vector field} is by definition $\hat{\eta}'_0:=\rho^{-n^{(0)}}\hat{\eta}_0$.

There are two possible scenarios as we explain in the following definition:

\begin{definition}\label{def:dicritical}
The blowing-up $\sigma_0$ is called \emph{non-dicritical} if $A_\rho(z,0)\equiv 0$ and \emph{dicritical} if $A_\rho(z,0)\neq 0$. Alternatively, we say that $E_0$ is \email{non-dicritical} or that $E_0$ is \emph{dicritical}, respectively.
\end{definition}

Despite of the formal nature of $\hat \eta_0$, the restriction $\hat \eta_0'|_{F_0}$ to the curve $F_0:=E_0\cap \{\theta=0\}$ is a well defined vector field. This restriction has polynomial coefficients in the coordinate $z$ along $F_0$. Therefore, its singular locus:
$$\text{Sing}(\hat \eta'\lvert_{F_0}):=\{    a\in F_0:\  \hat \eta'_0 \lvert_{F_0}(a)=0\}=\{    (z,0):\  A_\rho(z,0)=A_z(z,0)=0 \}$$
is finite. Singular points are points where we have to focus in order to define successive blowing-ups. But, in the dicritical case, we have to add those non-singular points where the vector field is tangent to the divisor. To be used for later, we recall the definition of such non-transversal points in the general situation of a normal crossing divisor. It can be seen for instance in Cano, Cerveau and Deserti's book \cite{cano-cerv-deserti} in the complex holomorphic context.
\begin{defi}
	Let $F$ be a germ of non-empty normal crossing divisor at $0\in \Rdos$ (that is, $F$ is given locally by $xy^\epsilon=0$) for some coordinates $(x,y)$ where $\epsilon=0$ or $\epsilon=1$). Let $\chi$ be a formal vector field at $0$. We say that 0 is an \emph{adapted singularity of $\chi$ relatively to $F$} if either $\chi(0)=0$ (a singular point) or $\chi(0)\neq 0$ and the unique formal invariant curve $C$ of $\chi$ through $0$ is such that $F\cup C$ has no normal crossings at $0$.
\end{defi}
We will consider this definition in a non-germified situation. Namely, suppose that $F$ is a normal crossing divisor in a real analytic manifold $N$ (possibly with boundary and corners). Let $(U,(x,y))$ be a chart in $N$ and let $\chi$ be a formal vector field in coordinates $(x,y)$ such that for any $a=(a_1,a_2)\in F\cap U$, the transformed vector field $\chi_a:=\chi(\tilde x+ a_1,\tilde y + a_2)$ is well defined as a formal vector field in coordinates $(\tilde x,\tilde y)$. Then, 
 the \emph{adapted singular locus of $\chi$ relatively to $F$} denoted by $\Singad (\chi,{F})$ is the set of points $a\in F\cap U$ for which the origin is an adapted singularity of $\chi_a$ relatively to $F_a$, the germ of $F$ at $a$.

Applied to our reduced vector field $\hat \eta_0'$ and to $F_0$, we have
\begin{itemize}
	\item[a)] If $E_0$ is non-dicritical, then $\Singad(\hat\eta_0', F_0)=\Sing(\hat \eta_0'\lvert_{F_0})$.
	\item[b)] If $E_0$ is dicritical, then $\Singad (\hat \eta_0,F_0)=\Sing(\hat \eta_0'\lvert_{F_0})\cup \{   (z,0): A_\rho(z,0)=0 \}$
\end{itemize}
In both cases, the adapted singular locus $\Singad(\hat\eta_0, F_0)$ is finite.

We define also the transforms $\hat \xi^{(\infty)}:= (\sigma_0\lvert_{C_\infty})^\ast \hat \xi$ and  $\hat \xi^{(-\infty)}:= (\sigma_0\lvert_{C_{-\infty}})^\ast \hat \xi$ of $\hat \xi$ in the charts $C_\infty, C_{-\infty}$, respectively.
The expressions for $\hat{\xi}
^{(\infty)}$, using simplified notation\linebreak 
$(x,y,z):=(x^{(\infty)},y^{(\infty)},z^{(\infty)})$ is the following:
	\begin{equation}\label{eq:campoenCinfty}
		\begin{split}
			\hat{\xi}^{(\infty)}=& R^{(\infty)}({x}^2+{y}^2, z) \left({x}\frac{\partial}{\partial {x}}+ {y}\frac{\partial}{\partial {y}} \right)  + 
			T^{(\infty)}({x}^2+{y}^2, z) \left(-{y}\frac{\partial}{\partial {x} }+ {x}\frac{\partial}{\partial {y}} \right)\\ &+ 
			Z^{(\infty)}({x}^2+{y}^2, z) \frac{\partial}{\partial z},
		\end{split}
	\end{equation}
where $R^{(\infty)},T^{(\infty)},Z^{(\infty)}\in \mathbb{R}[{x}^2+{y}^2][[z]]$ are given by: $$R^{(\infty)}({x}^2+{y}^2, z)=R(({x}^2+{y}^2)z^2, z)-\frac{1}{z}Z(({x}^2+{y}^2)z^2, z),$$ $$T^{(\infty)}({x}^2+{y}^2, z)=T(({x}^2+{y}^2)z^2, z) \text{ and } Z^{(\infty)}({x}^2+{y}^2, z)=Z(({x}^2+{y}^2)z^2, z).$$ In a similar way, we obtain expressions for $\hat \xi^{(-\infty)}$.

\subsection{ Characteristic cycles and successive blowing-ups}\label{sec:sucessive-blow-up}
Recall that the adapted singular locus of $\hat \eta'_0$ is finite. Its elements, belonging to $F_0=\{  \theta=\rho^{(0)}=0 \}$ are determined by the $z^{(0)}-$coordinate in the chart $C_0$. Denote then by
$$\wt{\text{Sing}}(\hat \eta_0', F_0)=\{ (\omega_i^{(0)},0): i=1,\ldots,m_0   \}, \quad \text{ with }\omega_i^{(0)}< \omega_j^{(0)}\text{ if } i<j.$$

\begin{defi}
The \emph{characteristic cycles} of $\hat{\xi} $ in ${M}_0$ are the connected components of the set  $\SSS^1\times \wt{\text{Sing}} (\eta_0'\lvert {F_0}) \subset C_0$, that is, the circles in the divisor $E_0$ given by $\gamma_i:=\{  z^{(0)}= \omega^{(0)}_i, \rho^{(0)}=0 \}$ for $i=1,2,\ldots, m_0$. 
The origins $\gamma_\infty, \gamma_{-\infty}$ of the charts $C_\infty$ and $C_{-\infty}$ (cf. equations \eqref{eq:cartaCinfty} and \eqref{eq:cartamenosCinfty}) are called the \emph{characteristic singularities of $\hat{\xi}$ in ${M}_0$}. We use the term \emph{characteristic elements} to refer to either the characteristic cycles or characteristic singularities. 
\end{defi}

In the rest of this section, we inductively define certain {sequences of blowing-ups} attached to $\hat{\xi}$ starting from the data defined above for the first blowing-up $\sigma_{0}$.

The tuple $\mathcal{M}_0:=(M_0,\sigma_0,\mathcal{A}_0,\mathcal{D}_0)$, where:
\begin{itemize}

	\item  $\mathcal A_0$ is the atlas of $M_0$ composed by the charts $C_{-\infty},C_0,C_\infty$, 
	\item  $\mathcal{D}_0$ is the family of \emph{characteristic elements} of $\hat \xi$ in $M_0$, that is, $\mathcal{D}_0:=\{ \gamma_{-\infty},\gamma_1, \ldots, \gamma_{m_0}  ,\gamma_\infty\}$, 
\end{itemize}
is by definition a \emph{sequence of admissible blowing-ups of length $l=0$  for $\hat{\xi}$}.
Suppose that we have already defined sequences of admissible blowing-ups for $\hat \xi$ of length $l-1$, consisting on tuples $\mathcal{M}=(M,\pi,\mathcal{A}, \mathcal{D})$ satisfying the following hypothesis:
\begin{itemize}
	\item[(H1)] $\pi:(M,E)\longrightarrow(\Rtres,0)$ is a sequence of (real) blowing-ups with smooth analytic centers factorizing through $\sigma_0$ (i.e., $\pi= \sigma_0 \circ \bar \pi$, where $\bar \pi: M\longrightarrow M_0$ is either the identity or a sequence of blowing-ups with smooth analytic centers).
	\item[(H2)] $\mathcal{D}= \{  \gamma_I \}_{I\in \mathcal I}$ is a finite family of disjoint closed subsets of the divisor $E= \pi^{-1}(0)$, such that:
	\begin{itemize}
		\item There are two elements in $\mathcal D$ with indices $I_\infty^{\mathcal M}= (\infty,\stackrel{s}{\ldots}, \infty)$ and $I_{-\infty}^{\mathcal M}= (-\infty,\stackrel{t}{\ldots}, -\infty)$ for some $s,t\in \mathbb{N}_{\geq 1}$,  that are the two points where $E$ intersects the strict transform $\pi^\ast(\{x=y=0\})$ of the $z-$axis. They are called the \emph{characteristic singularities (of $\hat \xi$ in $M$)}.
		\item The rest of the elements $\gamma_I$, with $I\neq I_{-\infty}^{\mathcal M}, I_{\infty}^{\mathcal M} $, are analytic embedded circles called \emph{characteristic cycles (of $\hat \xi$ in $M$)}.
		\item The intersection of any pair of components of $E$ is an element of $\mathcal D.$ Each of them is called a \emph{corner characteristic cycle}. The corner characteristic cycles are those indexed by tuples $I=(i_1,\ldots, i_r)\neq I_{-\infty}^{\mathcal M}, I_{\infty}^{\mathcal M}$ for which $i_r=\pm \infty.$
	\end{itemize}
\item[(H3)] $\mathcal A=\{C_J\}_{J \in \mathcal{J}}$ is an atlas of $M$ with the following properties:
\begin{enumerate}
	\item	There are  charts $(C_J,(x^{(J)} , y^{(J)},z^{(J)}   ))$ centered at 
	the characteristic singularities $\gamma_J$, with $J\in \{ I_\infty^{\mathcal M}, I_{-\infty}^{\mathcal M}     \}$, satisfying $E\cap C_J= \{  z^{(J)} =0 \}$. 
	 Moreover, the expression of $\pi$ in the chart $C_{J}$ with $J=I_{\epsilon \infty}^\mathcal{M}$, for $\epsilon=\pm 1$, is 
	$$\pi(x^{(J)} , y^{(J)},z^{(J)})=((z^{(J)})^rx^{(J)} , (z^{(J)})^ry^{(J)},\epsilon z^{(J)})$$
	 with $r\in \mathbb{N}_{\geq 1}$ ($r$ may depend on $\epsilon$). Furthermore, the coefficients of $\xi^{(J)}:=(\pi\lvert_{C_J})^\ast \hat \xi$ belong to $\R[x ^{(J)},y^{(J)}][[z^{(J)}]]$.

	\item If $J\notin \{ I_\infty^{\mathcal M}, I_{-\infty}^{\mathcal M}     \}$, the chart $(C_J,(\theta,z^{(J)},\rho^{(J)}))$ is defined for  $\theta\in \mathbb R $, $z^{(J)}$ in $\R$ or $\R_{\geq 0}$ and $\rho^{(J)}\in \R_{\geq 0}$ (with the same observation as in Remark~\ref{rem:theta}), and satisfies $E\cap C_J=\{ \rho^{(J)} (z^{(J)})^\epsilon=0 \}$ with $\epsilon=0$ or 1 according to $z^{(J)}$ being defined either in $\R$ or $\R_{\geq 0}$, respectively. 
	In the case $\epsilon=0$, the chart $C_J$ is a \emph{non-corner chart} and the characteristic cycles contained in $E\cap C_J$ are given by equations $\{   z^{(J)}=a_i , \rho^{(J)}=0 \}$, where $\{  a_i \}_i$ is a finite collection of real numbers. In the case $\epsilon=1$, the chart $C_J$ is a \emph{corner chart} and the family of characteristic cycles contained in
	$E\cap C_J$ consists of a unique corner characteristic cycle given by $\{   z^{(J)}=0, \rho^{(J)}=0 \}$ and a collection of non-corner characteristic cycles given either by $\{   z^{(J)}=b_j , \rho^{(J)}=0 \}_j$ for a family $\{b_j\}_j$ of positive numbers or by $\{ z^{(J)}=0 , \rho^{(J)}=c_k \}_k$ for a family $\{c_k\}_k$ of positive numbers.
	\item For any $J\notin \{ I_\infty^{\mathcal M}, I_{-\infty}^{\mathcal M}  \}$ the expression of $\pi$ in $C_J$ is polynomial in $(\cos \theta, \sin \theta, z^{(J)},\rho^{(J)})$ and the transformed vector field $\hat \xi^{(J)}=(\pi\lvert_ {C_J})^\ast \hat \xi$, written with simplified notation $(\theta,z,\rho)=(\theta, z^{(J)}, \rho^{(J)})$ as
	 \begin{equation}
	 	\label{eq:campoformalexplokMhyp}
	 	\hat{\xi}^{(J)}
	 	=
	  B_\theta
	 	^{(J)}(z,\rho) \thetaa\\
	 	+ B_z
	 	^{(J)}(z,\rho) \zhat +	B
	 	^{(J)}_\rho(z,\rho)
	 	\rhoohat,
 \end{equation}
satisfies that $B_i^{(J)}$, for $i=\theta,z,\rho$, belongs to $\R [z][[\rho]]$ if $C_J$ is a non-corner chart,  or to both algebras $\R [z][[\rho]]$ and $\R [\rho][[z]]$, if $C_J$ is a corner characteristic chart. In any case, $B_\theta^{(J)}$ is a unit of that corresponding algebra and $B_\theta(0,0)=1$.
\end{enumerate}
\end{itemize} 
 \strut
 When $J=I_{-\infty}^{\mathcal M}$ or $J=I_{\infty}^{\mathcal M}$, we define $n^{(J)}$ to be the maximum $n$ such that $\xi^{(J)}(z^{(J)})=(z^{(J)})^{n}\cdot \tilde B (x^{(J)} , y^{(J)},z^{(J)})$ with $\tilde B$ an element in $\R [x^{(J)},y^{(J)}][[z^{(J)}]]$.

Observing (H3)-(2), for $J\notin \{ I_\infty^{\mathcal M}, I_{-\infty}^{\mathcal M}     \}$, we define the formal two dimensional vector field $\hat \eta_J$ \emph{associated to the transform $\hat{\xi}^{(J)}$}, 
 as the one given by the following system of ODEs (using \eqref{eq:campoformalexplokMhyp} and simplifying $(z,\rho)=(z^{(J)},\rho^{(J)})$):
\begin{equation}
	\label{eq:formal2dimgeneral}
	\hat\eta_J: \left\{ \begin{array}{lcll}
		
		\frac{d z}{d \theta} & 
		=& B_\theta^{(J)}(z,\rho)^{-1}B_z^{(J)}(z,\rho)  &=\rho^{n_1^{(J)}} z^{n_2^{(J)}} 
		A_z^{(J)}(z,\rho)\\
		\frac{d \rho}{d \theta} &
		= &B_\theta^{(J)}(z,\rho)^{-1}B_\rho^{(J)}(z,\rho)  &=\rho^{n_1^{(J)}} z^{n_2^{(J)}} 
		A_\rho^{(J)}(z,\rho).
	\end{array} \right. 
\end{equation}
Here, $n^{(J)}_1 $ is the maximum $n$ such that $\rho^n$ divides both $B_\rho^{(J)}$ and $B_z^{(J)}$ (and thus, $A_\rho^{(J)}$ and $A_z^{(J)}$ are formal series not both together divisible by $\rho$). If $C_J$ is a the non-corner chart, we take $n_2^{(J)}=0$. If $C_J$ is a corner chart, we take $n_2^{(J)}$ to be the maximum $m$ such that $z^m$ divides both $B_\rho^{(J)}$ and $B_z^{(J)}$. We define $n ^{(J)}:=\max \{    n^{(J)}_1, n^{(J)}_2 \} $ in both cases. 
\strut

\begin{figure}
	\centering
	\includegraphics[width=\textwidth]{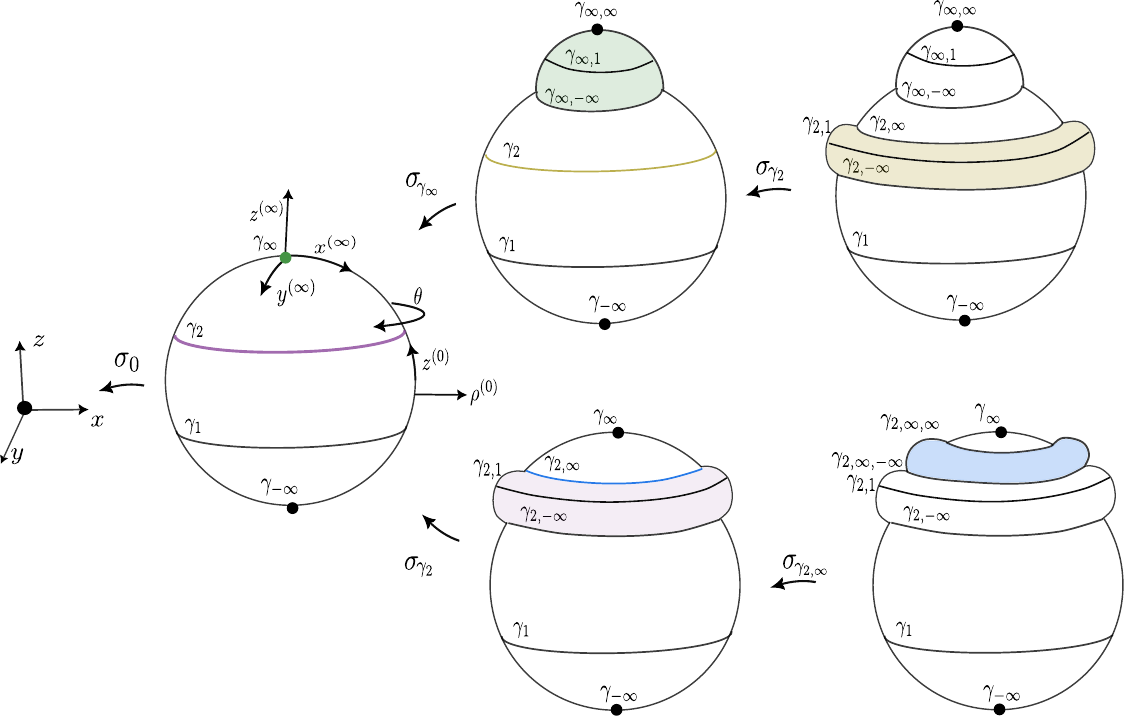}
	\caption{Several sequences of admissible blowing-ups.}
	\label{fig:admblowingups}
\end{figure}

It is clear that $\mathcal M_0$ fulfill (H1-H3). Now, a \emph{sequence of admissible blowing-ups for $\hat \xi$ of length $l$} is a tuple $\mathcal{M'}:=(M',\pi',\mathcal{A}', \mathcal{D}')$ built from a sequence of admissible blowing-ups $\mathcal{M}=(M,\pi,\mathcal{A}, \mathcal{D})$ of length $l-1$ in such a way that $\pi'= \pi \circ \sigma_{\gamma_{I}}$, where
$$ \sigma_{\gamma_I}: M'\longrightarrow M $$
is the blowing-up centered at some $\gamma_I\in\mathcal D$. The expression of $\sigma_{\gamma_{I}}$ in charts and the description of the families $\mathcal D',\mathcal A'$ are exposed in what follows (see Figure~\ref{fig:admblowingups} for an illustration of the different situations).
We consider two cases: the blowing-up $\sigma_{\gamma_{I}}$ is centered at a characteristic singularity or at a characteristic cycle. 

i) \textbf{First case}. The blowing-up
$\sigma_{\gamma_I}$ is centered at the singular point $\gamma_{I}$ with $I=I_\infty^\mathcal{M}$ (or analogously for $I=I_{-\infty}^\mathcal{M}$. Put $J=I$, $J_\infty=(J,\infty)$ and $J_{0}=(J,0)$. The point $\gamma_{I}$ is the origin of a chart $(C_{J},(x^{(J)}, y^{(J)}, z^{(J)}))\in \mathcal A$ where the divisor $E=\pi^{-1}(0)$ is given by $\{  z^{(J)}=0 \}$. 
Then, the exceptional divisor $\sigma_{\gamma_I}^{-1}(\gamma_I)$ is covered by the two charts of $M'$:  $(C_{J_\infty},(x^{(J_\infty)}, y^{(J_\infty)}, z^{(J_\infty)}))$ and $(C_{J_0}, ( \theta
, z^{(J_0)},\rho^{(J_0)}))$, so that $\sigma_{\gamma_I}$ is written as:
\begin{eqnarray}
	\label{eq:cartaC00}
	&\text{In }C_{J_0}:& \left\{ \begin{matrix}
		x^{(J)} &=& \rho^{(J_0)} \cos \theta
		\\
		y^{(J)}  &=& \rho^{(J_0)} \sin \theta
		\\
		z^{(J)}  &=& \rho^{(J_0)}z^{(J_0)}
	\end{matrix} \right.  \qquad  \theta
	\in \R, z^{(J_0)} ,\rho^{(J_0)}  \geq 0\\
	\label{eq:cartaC0infty}
	&	\text{In }C_{J_\infty}: & \left\{ \begin{matrix}
		x ^{(J)}&=& x ^{(J_\infty)} z ^{(J_\infty)}\\
		y^{(J)}  &=& y ^{(J_\infty)} z ^{(J_\infty)} \\
		z ^{(J)} &=& \space z ^{(J_\infty)}
	\end{matrix} \right.\qquad  x ^{(J_\infty)} ,y ^{(J_\infty)}  \in \mathbb{R}, \space z ^{(J_\infty)} \geq 0.
\end{eqnarray}
We set the atlas of $\mathcal M'$ to be $\mathcal A'=(\mathcal A  \setminus \{  C_J \}) \cup \{   C_{J_0}, C_{J_\infty} \}$ (under the identification $\sigma_{\gamma_J}:M'\setminus \sigma_{\gamma_J} ^{-1}{(\gamma_J)}\to M'\setminus {\gamma_J} $). The chart $C_{J_0}$ is a corner chart in this case.

Let us define $\mathcal{D}$. We consider first the vector field $\hat \xi  ^{(J_\infty)}:=(\sigma_ {\gamma_I}\lvert_{C_{J_\infty}})^\ast \hat \xi ^ {(I)}$. 
 Define $n ^{(J_\infty)}$ as the maximum $n\in \mathbb N$ such that $(z ^{(J_\infty)})^{n}$ divides $\hat \xi  ^{(J_\infty)}(z^{(J_\infty)}).$ 
In the chart $C_{J_\infty}$, the expression of the vector field $\hat \xi^{(J_{\infty})}$ is similar to \eqref{eq:campoenCinfty}. The origin of $C_{J_\infty}$ is named $\gamma_{J_\infty}$. 

Consider the formal vector field $\hat{\xi}^{(J_0)}:=(\sigma_{\gamma_I}\lvert_{C_{J_0}})^\ast(\hat{\xi}^{(J)})$. Use the expression in \eqref{eq:campoformalexplokMhyp} for $\mathcal M$, and rename $(z,\rho)=(z^{(J_0)},\rho^{(J_0)})$. Then, $\hat{\xi}^{(J_0)}$ is  given by 
\begin{equation}
	\label{eq:campoformalexplok}
	\hat{\xi}^{(J_0)}
	=
 B_\theta
	^{(J_0)}(z,\rho) \thetaa\\
	+ B_z
	^{(J_0)}(z,\rho) \zhat  +	B
	^{(J_0)}_\rho(z,\rho)
	\rhoohat ,
\end{equation}
where $B_\theta
^{(J_0)}\in \mathbb{R}[z][[\rho]]$ is a unit in the algebra $\R[z][[\rho]]$ (because $B^{(J_0)}_\theta(z,0)=1$) and
\begin{equation}\label{eq:BrhonBtildek}
	B_z
	^{(J_0)}=\rho^{n_1^{(J_0)}} z^{n_2^{(J_0)}}  \tilde{B}_z, \quad
B
^{(J_0)}_\rho=\rho^{n_1^{(J_0)}} z^{n_2^{(J_0)}} \tilde{B}
^{(J_0)}_\rho,
^{(J_0)} \quad \tilde{B}_\rho
^{(J_0)},\tilde{B}_z
^{(J_0)}\in \mathbb{R}[z][[\rho]].
\end{equation}
with $n_1^{(J_0)}, n_2^{(J_0)}$ defined similarly as we have defined $n_1^{(J)}$ and $n_2^{(J)}$.
The {vector field associated to $\hat{\xi}^{(J_0)}$} is the two dimensional vector field $\hat{\eta}_{J_0}$ with coefficients in $\mathbb{R}[z][[\rho]]$,  
 defined in a similar manner as $\hat \eta_0$ in \eqref{eq:formal2dim}. That is,
\begin{equation}\label{eq:twodimVFIinfty}
	\hat{\eta}_{J_0}= \rho^{n_1^{(J_0)}} z^{n_2^{(J_0)}}  \left( A_z^{(J_0)}(z,\rho) \z + A_\rho^{(J_0)}(z,\rho)\rhoo    \right).
\end{equation}
where $A_k^{(J_0)}=B_k^{(J_0)}\cdot(B_\theta^{(J_0)})^{-1}$ for $k=z,\rho.$
The vector field $\hat{\eta}_{J_0}':=\frac{1}{\rho^{n_1^{(J_0)}} z^{n_2^{(J_0)}}}\hat{\eta}_{J_0}$
is called the \emph{associated reduced vector field of $\hat{\xi}^{(J_0)}$}. 
We distinguish two cases:
 \begin{itemize}
 	\item The blowing-up $\sigma_{\gamma_I}$ is \emph{non-dicritical} if $A_\rho^{(J_0)}(z,0)\equiv 0$. In this case, we say that the divisor $\sigma_{\gamma_I}^{-1}(\gamma_I)$ is a \emph{non-dicritical component} of the total divisor $E':=(\pi')^{-1}(0)$.
 	\item The blowing-up $\sigma_{\gamma_I}$ is \emph{dicritical} if $A_\rho^{(J_0)}(z,0)\neq 0$ and $\sigma_{\gamma_I}^{-1}(\gamma_I)$ a \emph{dicritical component} of the total divisor $E'$.
 \end{itemize}

 Put $E_{J_0}:=\sigma_{\gamma_{I}}^{-1}(\gamma_I)$ and $F_{J_0}:=E_{J_0}\cap \{   \theta=0 \}$ and consider $\Singad(\hat \eta_{(J_0)}', F_{J_0})$ the adapted singular locus of $\hat \eta_{(J_0)}'$ relatively to $F_{J_0}.$ Taking into account that the coefficients of $\hat \eta_{(J_0)}'$ belong to $\R[z][[\rho]]$, it is a finite set. Let us denote those contained in the regular part $\dot F_{J_0}:= F_{J_0}\cap \{  z>0 \}$ of $F_{J_0}$ as
$$\wt{\text{Sing}}(\hat \eta_{J_0}', {F_{J_0}})\cap \dot F_{J_0}=\{ (\omega_i^{(J_0)},0): i=1,\ldots,m_{J_0}   \} 
 \text{, with }0<\omega_i^{(J_0)}< \omega_j^{(J_0)}\text{ if } i<j.$$
The circles $\gamma_{I,i}:=\{  z=\omega_i^{(J_0)}, \rho=0 \}\subset E_{J_0}$ for each $i=1,\ldots, m_{J_0}$ are by definition the {\em non-corner characteristic cycles (in  $C_{J_0}$)}. The circle $\gamma_{I,-\infty}:=\{  z=0, \rho=0 \}$ is named \emph{a corner characteristic cycle}. 

We define the family $\mathcal D':= \{ \gamma_I  \}_{I\in \mathcal I'}$ of \emph{characteristic elements of $\mathcal M'$}, where  
$$\mathcal{I}'=\left(\mathcal{I}  \setminus\{  {I}   \}    \right) \cup\left( \displaystyle \bigcup_{i=1}^{m_{J_0}} \{   {(I, i)}   \} \right) \cup  \{  {(I, -\infty),(I, \infty)}   \} .$$
The elements of $\mathcal D'$ are subsets of $E'=(\pi')^{-1}(0)$, once we identify $\gamma_L= \sigma_{\gamma_I}^{-1}(\gamma_L)$ for $L\in \mathcal{I} \setminus \{ I \}$. They are either the two points $\gamma_{\infty, \ldots, \infty}$ and $\g_{-\infty, \ldots, -\infty}$ (whose indices are denoted also by ${I^{\mathcal{ M}'}_{ \infty}}$ and ${I^{\mathcal{ M}'}_{ -\infty}}$, respectively) called the characteristic singularities of $\hat{\xi}$ in $\mathcal{M}'$ or circles (the characteristic cycles of $\hat \xi$ in $\mathcal{ M}'$).

ii) \textbf{Second case}. $\sigma_{\gamma_I}$ is centered at one of the characteristic cycles $\gamma_{I}\in \mathcal D$ with $I=(i_1,\ldots,i_r)$. It can be a corner characteristic cycle (in which case $i_r=\pm \infty$) or not. The charts after the blowing-up $\sigma_{\gamma_I}$ are defined in a different manner in each case. In order to simplify the notation, name $I'=(i_1,\ldots, i_{r-1})$.
\begin{itemize}
\item[a)] When ${\gamma_I}$ is a corner characteristic cycle, it can be seen as $\{\rho^{(J)}=0,z^{(J) }=0\} $ in a chart $C_{J}$ by (H3). 
Put $J_0=(I,0)$ and $J_\infty=(I,\infty)$. The set $\sigma_{\gamma_I} ^{-1}(\gamma_I)$ is covered by two new charts $(C_{J_\infty},(\theta
,z ^{(J_\infty)},\rho ^{(J_\infty)}))$ and $(C_{J_0},(\theta
,z^{(J_0)},\rho^{(J_0)}))$, where the blowing-up $\sigma_{\gamma_I}$ is written as:
\begin{eqnarray}
	  \label{eq:centralchartxx}
	  &\text{In }\quad	C_{J_\infty}: &  \left\{
	  \begin{matrix}
	  	\theta
	  	&=&\theta\\
	  		z^{(J)}&=& z ^{(J_\infty)}\\
	  	\rho^{(J)}&=&\rho ^{(J_\infty)} z ^{(J_\infty)}  	  
	  \end{matrix}
	  \right. , \quad   \theta
	  \in \mathbb R, \ \rho ^{(J_\infty)}, z ^{(J_\infty)}\geq 0 ,\\
	  \label{eq:centralchartxxm}
	 & \text{In }\quad	C_{J_0}: &\left\{
	  \begin{matrix}
	  	\theta
	  	&=&\theta\\
	  		z^{(J)}&=& \rho^{(J_0)} z^{(J_0)}\\
	  	\rho^{(J)}&=&\rho^{(J_0)} 
	  \end{matrix}
	  \right., \quad  \theta
	  \in \mathbb R, z^{(J_0)} ,\rho^{(J_0)},  \geq 0.
\end{eqnarray}
The new atlas is defined by $\mathcal A':=(\mathcal A \setminus \{ C_{J} \})\cup \{  C_{J_0},C_{J_\infty}   \}  $, where we have identified $M'\setminus \sigma_{\gamma_I} ^{-1}(\gamma_I)$ and $M\setminus \gamma_I$ via $\sigma_{\gamma_I}$.

To determine the new family $\mathcal D'$  of characteristic elements in this case, we write the transformed formal vector fields $\hat \xi^{(J_0)}:=(\sigma_{\gamma_{I}}\lvert_{C_{J_0}})^{\ast}\hat \xi^{(J)} $,  $\hat \xi^{(J_\infty)}:=(\sigma_{\gamma_{I}}\lvert_{C_{J_\infty}})^{\ast}\hat \xi^{(J)} $
in the two charts. Both are similar and, in fact, to determine $\mathcal{D}'$ only one of the expressions is sufficient. Considering for instance the chart $C_{J_0}$, and with similar computations and notations as in the precedent paragraphs, we get (simplifying $(z,\rho)= (z^{(J_0)}, \rho^{(J_0)})$)
\begin{equation}
	\label{eq:campoformalexploiia}
	\hat{\xi}^{(J_0)}
	=
	B_\theta
	^{(J_0)}(z,\rho) \thetaa\\
	+ B_z
	^{(J_0)}(z,\rho) \zhat  +	B
	^{(J_0)}_\rho(z,\rho)
	\rhoohat ,
\end{equation}
where  $B_\theta
^{(J_0)},B_z^{(J_0)},B_\rho^{(J_0)}\in \mathbb{R}[z][[\rho]]$ and $B_\theta^{(J_0)}$ is a unit. The \emph{associated two dimensional vector field} is $ \hat \eta_{J_0}= (B_\theta^{J_0})^{-1}B_{z}^{(J_0)} \z + (B_\theta^{J_0})^{-1}B_{\rho}^{(J_0)} \rhoo  $ and we put
\begin{equation*}\label{eq:campoformalexploiia2dim}
	\hat \eta_{J_0}= \rho^{n_1^{(J_0)}} z^{n_2^{(J_0)}}\hat \eta_{J_0}'= \rho^{n_1^{(J_0)}} z^{n_2^{(J_0)}} \left( A_z^{(J_0)}\z+     A_\rho^{(J_0)} \rhoo      \right). 
\end{equation*}
The natural numbers $n_k^{(J_0)}$ for $k=1,2$ are defined as above. The vector field  
$\hat \eta_ {J_0}'$ is the \emph{associated reduced formal vector field} 
We distinguish the cases when $\sigma_{\gamma_{I}}$ (or the component $E_{J_0}:=\sigma_{\gamma_{I}}^{-1}(\g_I)$) is dicritical $ (A_\rho^{(J_0)}(z,0)\neq 0 )$ or non-dicritical $ (A_\rho^{(J_0)}(z,0)\equiv 0 )$. Put $F_{J_0}:= E_{J_0}\cap \{ \theta=0 \}, \ \dot F_{J_0}:=F_{J_0}\cap \{ z>0 \}  $ 
and denote
 $$\wt{\text{Sing}}(\hat \eta_{J_0}', { F_{J_0}})\cap \dot F_{J_0}=\{ (\omega_i^{(J_0)},0): i=1,\ldots,m_{J_0}   \} \text{, with }0<\omega_i^{(J_0)}< \omega_j^{(J_0)}\text{ if } i<j.$$
With these data, we put:
$$  \gamma_{I,i}:= \SSS^1 \times \{  (z^{(J_0)},\rho^{(J_0)})=(\omega_i ^{(J_0)},0)  \} , i =1, \ldots m_{J_0} $$ $$\  \gamma_{I,-\infty}:= \SSS^1 \times \{  (z^{(J_0)},\rho^{(J_0)})=(0,0)  \}  \subset C_{J_0} $$
$$  \gamma_{I,\infty}:= \SSS^1 \times \{  (z^{(J_\infty)},\rho^{(J_\infty)})=(0,0)  \}  \subset C_{J_\infty} $$
and we define the family of charactersitic elements of $\mathcal{ M}'$ as $\mathcal D'= \{   \gamma_I \}_{I\in \mathcal I}$, where
\begin{equation*}\label{eq:indices}
	\mathcal{I}'=\left( \mathcal{I}\setminus \{ I \} \right)\cup  \displaystyle  \{ (I,i)  \}_{i=1}^{m_{J_0}} \cup \{ (I,-\infty)    \} \cup \{ (I,\infty)    \}.
\end{equation*}
 identifying $\gamma_L$ with $\sigma_{\gamma_{I}}^{-1}(\gamma_L)$ for $L\in \mathcal I \setminus \{  I \}$. Notice that, among the new characteristic cycles, $\gamma_{I,\infty}$, $\gamma_{I,-\infty}$ are corner cycles and the other ones are non-corner characteristic cycles.

\item[b)] When $\gamma_I$ is a non-corner characteristic cycle (that is, by (H2), when $I=(i_1, \ldots, i_{r})$ with $i_{r}\neq \pm \infty$), it can be seen as the set $\gamma_I=\{z^{(J )}=\omega^{(J)}_{k}, \rho^{(J)}=0\} $ for some $\omega_k^{(J )}$ in the domain of $z^{(J)}$ of a chart $C_{J}$, by (H3). Set  $J_{-\infty}:=(I,-\infty)$, $J_\infty:=(I,\infty)$ and $J_0:=(I,0)$. The blowing-up $\sigma_{\gamma_I}:(M',E')\longrightarrow (M,\gamma_I)$ of $\gamma_I$ is given in three new charts. 

\begin{eqnarray}
	\label{eq:cartaCIinfty}
  &	\text{In }	C_{J_\infty}:  & \left\{
	\begin{matrix}
		\theta
		&=&\theta\\
		z^{(J)}&=& z^{(J_ \infty)}- \omega^{(J)}_{k}\\
		\rho^{(J)}&=&\rho ^{(J_\infty)} z ^{(J_\infty)}  
	\end{matrix}
	\right. \qquad  \ \theta
	\in \mathbb R, \rho^{(J_ \infty)}, z ^{(J_\infty)}  \geq 0. \\
	\label{eq:cartaCI0}
	&\text{In }	C_{J_ 0}:& \left\{
	\begin{matrix}
		\theta
		&=&\theta\\
		z^{(J)}&=& \rho^{(J_ 0)} (z^{(J_0)}  - \omega^{(J)}_{k} )\\
		\rho^{(J)}&=&\rho^{(J_0)} 
	\end{matrix}
	\right. \theta
	\in \mathbb R,z^{(J_0)} \in \mathbb{R}, \rho^{(J_0)} \geq 0. \\
	\label{eq:cartaCImenosinfty}
	&\text{In }	C_{J_\infty}:     & \left\{
	\begin{matrix}
		\theta
		&=&\theta\\
			z^{(J)}&=& -z^{(J_{ -\infty})}+ \omega^{(J)}_{k}  \\
		\rho^{(J)}&=&\rho^{(J_{-\infty})} z^{(J_{ -\infty})}  
	\end{matrix}
	\right.  \quad \theta
	\in \mathbb R, \rho^{(J_{ -\infty})}, z^{(J_{ -\infty})}  \geq 0.
\end{eqnarray}
The new atlas is $\mathcal A':= \mathcal A\setminus \{  C_J \} \cup \{ C_{J_0}, C_{J_\infty}, C_{J_{-\infty}} \}$. The family $\mathcal D'$ of \emph{characteristic elements} of $\hat \xi$ in $\mathcal M'$ is defined analogously as in case a), studying the corresponding transformed vector field $\hat \xi^{(J_0)}= (\sigma_{\gamma_I}\lvert_{C_{J_0}})^\ast \hat \xi^{(J)}$, its associated two-dimensional vector field $\hat \eta_{J_0}$ and the adapted singular locus of the reduced associated vector field $\hat \eta_ {J_0}'$ relatively to to $F_{J_0}=C_{J_0}\cap \sigma_{\gamma_I}^{-1}(\gamma_I)\cap \{  \theta=0 \}$. We just observe the following:
\begin{itemize}
	\item The charts $C_{J_\infty}$ and $C_{J_{-\infty}}$ are corner charts and the curves $\gamma_{J_\infty}=\{z ^{(J_\infty)}=0,\rho ^{(J_\infty)}=0 \}$ and $\gamma_{J_{-\infty}}=\{z^{(J_{-\infty})}=0, \rho^{(J_{-\infty})}=0\}$  are corner characteristic cycles.
	\item 
	The chart $C_{J_0}$ is a non-corner chart and 
	
contains 
the new non-corner characteristic cycles $\gamma_{I,i}$ where $\gamma_{I,i}=\SSS^1\times \{ (\omega_i^{(J_0),0)}  \}, \ i=1, \ldots, m_{J_0}$
being $(0,\omega_i^{(J_0),0)})\in \Singad (\hat \eta_{J_0}',F_{J_0})$.
\end{itemize}
The set of indices is
\begin{equation*}\label{eq:indices2}
	\mathcal{I}'=\left( \mathcal{I}\setminus \{ I \} \right)\cup \left( \displaystyle \bigcup_{i=1}^{m_{J_0}}  \{ (I,i)  \} \right) \cup \{ (I,-\infty)    \} \cup \{ (I,\infty)    \}.
\end{equation*}
\end{itemize}
From the construction, we can check that the hypothesis (H1-H3) are fulfilled for $\mathcal M'$. Thus, we have defined admissible sequences of blowing-ups of $\hat \xi$ of any length.

\begin{remark}\label{rmk:ciclosendoscartas}
By construction, a non-corner characteristic cycle $\g_I\in \mathcal D$ is defined in some chart $C_J$ by equations $\g_I=\{   z^{( J)}=\omega^{(J)}_k , \rho^{( J)}=0  \}$  for some $k\in \mathbb{N}_{\geq 1}$. It may happen that the same characteristic cycle $\g_I$ is defined in a corner chart $C_{\tilde J}$ by $\{   z^{(\tilde J)}=0 , \rho^{(\tilde J)}=c  \}$, for some $c\in \mathbb R$.
\end{remark}

\begin{remark}\label{rmk:non-zero}
	Notice that for any given sequence of admissible blowing-ups $\mathcal{ M}=(M,\pi, \mathcal A,\mathcal D)$ and for any chart $C_J$ of $\mathcal{A}$, the associated vector field $\hat \eta_{ J}$ is not identically zero. This can be seen by construction of $\mathcal{ M}$ and using Remark~\ref{rmk:nondeg}.
\end{remark}

\subsection{Adapted reduction of singularities.} 
According to Panazzolo's result on reduction of singularities of three-dimensional vector fields~\cite{panazzolodim3}, any vector field $\xi \in \mathcal{H}^3$ or any of its formal normal forms $\hat \xi$ is already in final form, since it has a non-nilpotent linear part. 
Based on theorems of reduction of singularities of two-dimensional vector fields, we state a finer reduction of singularities of $\hat \xi$, adapted to our problem.

Recall that a formal vector field $\chi=A(x,y)\x+B(x,y)\y$ at $(\Rdos,0)$ has a \emph{simple singularity} if $\Sing(\chi)=\{0\}$, at least one of the eigenvalues $\lambda_1,\lambda_2$ of the linear part $D\chi(0)$ is different from zero, for instance $\lambda_2\neq0$, and $\frac{\lambda_1}{\lambda_2}\notin \mathbb{Q}_{>0}$. In this case, $\chi$ has exactly two formal invariant curves, also called \emph{separatrices,} which are tangent to the corresponding eigenspaces, non-singular and mutually transverse (see the book \cite{cano-cerv-deserti}). Moreover, the separatrices are real if and only if $\{ \lambda_1,\lambda_2  \}\subset \R$. We need an extended notion of simple singularity, also taken from that reference, that takes into account the existence of a divisor and the possibility that the singularity is not isolated.
\begin{defi}\label{def:adaptedsimplesing}
	Let $F= \{  xy^\epsilon=0 \}$ where $\epsilon\in \{ 0,1\}$, be a normal crossing divisor at $0\in \Rdos$. A formal vector field $\chi$ at $(\Rdos,0)$ has an \emph{adapted simple singularity relatively to  $F$} $\chi(0)=0$ if one of the two following situations occurs:
	\begin{enumerate}
		\item $\Sing(\chi)=\{0\}$, the singularity is simple and each component of $F$ is invariant for $\chi$ (thus, if $\epsilon=1$, the two components of $F$ are the two separatrices).
		\item If $\epsilon=0$ and $\Sing(\chi)$ is a formal non-singular curve $\Gamma$ transversal to $F=\{x=0\}$, given by an equation $\Gamma=\{  y-\hat{g}(x)=0 \}$ and $\chi=( y-\hat{g}(x))^r\bar{\chi}$ where $r\geq1$ and $\bar{\chi}$ is either non-singular at $0$ and $F$ is the only invariant curve of $\chi$ through 0, or $\bar \chi$ has a simple singularity at $0$ and the set of separatrices of $\bar{\chi}$ at $0$ is $\{ F,\Gamma  \}$.
	\end{enumerate}
\end{defi}
To distinguish the two cases of this definition, in the situation of \textit{(2)}, we say that $\chi$ has a \emph{non-saturated adapted simple singularity}. Usually in this situation, one divides $\chi$ by an equation of $\Sing(\chi)$ to get the situation in \textit{(1)} or a non-singular point. However for us, the vector field $\chi$ will come from some three dimensional vector field, hence it will be important to keep unaltered the singular locus, independently if it is of dimension 1 or 0.

\begin{prop}[Adapted resolution of singularities] \label{prop:adaptedressing}
There exists an admissible sequence of blowing-ups $\mathcal{M}=(M,\pi,\mathcal{A}, \mathcal{D} )$ of $\hat \xi$ with $\mathcal{A}=\{ C_J  \}_{J\in \mathcal J}$, $\mathcal{D}=\{  \gamma_I  \}_{I\in \mathcal I}$ and total divisor $E=\pi^{-1}(0)$ such that
\begin{enumerate}
\item For $J\in \{ I_\infty^{\mathcal M}, I_{-\infty}^{\mathcal M}     \}$, the transformed vector field $\hat{\xi}^{(J)}=(\pi\lvert_{C_J})^\ast \hat \xi$ satisfies $\hat \xi  ^{(J)}(z^{(J)})  =(z^{(J)})^t \cdot G$ where $t\geq 1$ and $G$ is a unit in $\mathbb{R}[[x^{(J)},y^{(J)},z^{(J)}]]$.
\item For any $J\in \mathcal J \setminus \{ I_\infty^{\mathcal M}, I_{-\infty}^{\mathcal M} \}$, the singularities of the reduced two-dimensional vector field $\hat \eta_J'$ are adapted simple relatively to the divisor $E\cap C_J\cap \{ \theta=0 \}$.
\item If $E_J\subset E_J$ is a dicritical component of $E$ 
 then, $E_J$ is isolated as a dicritical component (i.e. any other component that intersects $E_J$ is non dicritical).
 Moreover, one has $\Singad(\hat \eta_J, F_J)=\emptyset$ where $F_J=E_J\cap C_J\cap \{ \theta=0 \}$,
 in particular, $\hat \eta_J$ is everywhere transversal to $F_J$ along $F_J$.
\end{enumerate}
\end{prop}
\begin{proof}

Since $\hat{\xi}$ has isolated singularity, from Remark \ref{rmk:nondeg} there exists a term $c_jz^j$ in the coefficient $\hat \xi (z)$ of $\z$ with $c_j\neq 0$. Assume without loss of generality that $j$ is the minimum exponent with this condition. 
Write $\hat \xi (z)\in\R[[x,y,z]]$ as:
\begin{equation*}
	\hat \xi (z)= Z(x^2+y^2,z)=z^{t_0} G(x,y,z)=z^{t_0}  \displaystyle \sum_{k=\nu (G)}^{\infty} G_k(x,y,z),
\end{equation*}
where $G_k$ is an homogeneous polynomial of degree $k$, $t_0\geq 0$ is defined as the maximum such that $z^{t_0} $ divides $\hat \xi (z)$. Then $G_{j-{t_0} }(x,y,z)$ contains the monomial $c_j z ^{j-t_0}$ (notice that $j\geq t_0$ and the equality holds if and only if $\nu (G)=0$). Consider the first blowing-up $\sigma_0$ and study $\hat \xi^{(\infty)}(z^{(\infty)})
$, where $\xi^{(\infty)}=(\sigma_{0}\lvert_{C_\infty})^\ast \hat{\xi}$. Omitting super-indices for the coordinates $(x^{(\infty)}, y^{(\infty)}, z^{(\infty)})$, we have:
\begin{equation*}
	\hat \xi^{(\infty)} (z)= z^{t_0}\displaystyle \sum_{k=\nu(G)}^{\infty} G_k(x,y,1)z^k=  z^{t_1} \displaystyle \sum_{k=\nu (G)}^{\infty} G_k(x,y,1)z^{k-\nu(G)},
\end{equation*}
	where 
	$t_1=t_0+\nu(G)\geq t_0$. 
	Rewrite the series $G^{(1)}:=\sum_{k=\nu (G)}^{\infty} G_k(x,y,1)z^{k-\nu (G)}$ in homogeneous components:
	\begin{equation*}
		\hat \xi^{(\infty)} (z)=z^{t_1}G^{(1)}(x,y,z) =z^{t_1} \displaystyle \sum_{k=\nu(G^{(1)})}^{\infty} G_k^{(1)}(x,y,z).
	\end{equation*}
	If $j=t_1$, we see that $G^{(1)}_0=c_j$ and thus, $G^{(1)}$ is a unit, which gives statement \textit{(1)} of the proposition for $t=t_1$. Otherwise, if $t_1<j$, we see that $G^{(1)}_{j-t_0}(x,y,z)$ contains a term $c_jz^{j-t_1}$. Notice that in this case we have $t_1\geq t_0$ since, otherwise, if $t_1=t_0$ then $\nu(G)=0$ and $j=t_0=t_1$. Thus, $j-t_0>j-t_1\geq 0$.
By recurrence over $j-t_0$, there exists an admissible sequence of blowing-ups $\widetilde {\mathcal{ M}}= (\wt{M}, \wt{\pi},  \wt{\mathcal{A}}, \wt{\mathcal {D}})$
 with $\widetilde {\pi}$ a composition of $s$ blowing-ups at the corresponding characteristic singularities $\gamma_{I_\infty^{\mathcal{ M}_i}}$ such that, defining $t_0,t_1, \ldots, t_s$ as above, we have $j=t_s$. We 
 conclude \textit{(1)} for $\wt{\pi}^\ast \hat \xi$ at the characteristic singularity $\gamma_{I_\infty^{\widetilde{\mathcal{M}} } }$ with $t=t_l$. Analogously, up to blowing-up also the characteristic singularity $\gamma_{-\infty}^{\tilde{\mathcal M}}$, we may assume that \textit{(1)} holds at $\gamma_{I_{-\infty}^{\widetilde {\mathcal{M}}}}$.

According to the construction of sequences of admissible blowing-ups in the previous section, 
$\wt{\mathcal A}$ is composed by the two charts $C_{I_{-\infty}^{\widetilde{\mathcal M}}}$ and $C_{I_{\infty}^{\widetilde{\mathcal M}}}$ and a finite number of charts named:
\begin{equation}\label{eq:chartsaftertype1}
	C_0,C_{\infty,0},\ldots,C_{\infty,\stackrel{s-1}{\ldots},\infty,0}, C_{-\infty,0},\ldots,C_{-\infty,\stackrel{t-1}{\ldots},-\infty,0},
\end{equation}
with coordinates of the form $(\theta,z,\rho)\in \mathbb{R}\times (\R_{\geq 0})^2$ - except for the first one with $z$ taking values in $\R$. 
Denote by $\widetilde{\mathcal J}_0$ the set in the list of \eqref{eq:chartsaftertype1}. For any $J\in \wt{\mathcal {J}}_0$, consider the transformed vector field $\hat \xi^{(J)}=(\wt{\pi}\lvert_{C_J})^\ast \hat \xi$ and the corresponding reduced associated two-dimensional vector fields $\hat \eta_J'$. 
Denote by $\tilde E= \tilde \pi^{-1} (0)$ the total divisor of $\wt{\mathcal M}$
Notice that the coordinate $\theta$ is well defined in the union $U=\bigcup_{J\in \wt{\mathcal J}_0} C_{J}$ so that $\tilde F:=\tilde E\cap \{  \theta=0 \}\cap U$ has a perfect sense, in fact, $\tilde F =\tilde E \cap \overline{\pi^{-1}(\{y=0, x>0\})}$. Now, given $J\in \wt{\mathcal J}_0$ and $a\in \Singad(\hat \eta'_J, \tilde F)$, there exists a sequence of punctual two dimensional blowing-ups $\tau_a:(N_a, F)\to (C_J\cap \{\theta=0\},a)$ that provides a reduction of singularities of $\hat \eta_J'$ in the following enhanced sense:
\begin{enumerate}
	\item[(\textit{a})] For any point $q\in F=\tau_a^{-1}(a)$, if $\chi_q'$ is the strict transform of $\hat \eta_J'$ by $\tau_a$ at $q$ (that is, $\chi_q'=\frac{1}{f^k} \tau_a^\ast(\hat \eta_J')$, where $f$ is a local reduced equation of $F$ at $q$ and $k$ is maximal so that $\chi_q'$ has no pole), then, $q\in \Singad(\chi_q', F)$ if and only if $q\in \Sing(\chi_q')$.
	\item[(\textit{b})]  If $q\in F$ is a singular point of $\chi_q'$, then $q$ is an adapted singularity relatively to $F$ (c.f. Definition~\ref{def:adaptedsimplesing}).
\end{enumerate}
Recall that $\eta_{ J}'\neq 0$ by Remark~\ref{rmk:non-zero}.
In the case where $\hat \eta_J'$ has as isolated singularity at $a$, the existence of the reduction of singularities $\tau_a$ is a classical result by Seidenberg~\cite{Seidenberg1968ReductionOS}, together with an extended version that can be seen in Cano-Cerveau-Deserti's book \cite{cano-cerv-deserti} that eliminates adapted singularities relatively to the divisor that are not singular points (to get (\textit{a})). In the case where the singular locus $S_a:=\Sing_a(\hat \eta_J') $ at $a$ is not reduced to $\{ a\}$ (thus it is a formal curve), we first consider a reduction of singularities  $\psi_a:( N^{(1)}, F^{(1)})\to (C_J\cap \{\theta=0\},a)$ of $S_a$. Then, we blow up any adapted singularity $q$ of the strict transform $
 \chi_q'$ of $\hat \eta_J'$ by $\psi_a$ that is not a singular point to get (\textit{a}). After that, we may assume that such strict transform $
  \chi_q'$ either has an isolated singularity at $q$ (and hence we apply the procedure in \cite{cano-cerv-deserti}) or $q$ is a point in the strict transform $ S_{a,q}$  of the curve $S_a$ by $\psi_a$.  In this last case, there are coordinates $(x,y)$ at $q$ such that $ F^{(1)}=\{x=0\}$, locally, and $
  \chi_q'$ is written as $
  \chi_q'=(y-\tilde g(x))^r  \bar \chi_q' $, where $\{  y-\hat g(x)=0 \}$ is an equation of $ S_{a,q}$, $r\geq 1$ and $\bar \chi_q'$  has at most an isolated singularity at $q$. If $\bar \chi_q'(q)=0$, after a reduction of singularities of $\bar \chi_q'$, we may assume that $q$ is a simple singularity of $\bar \chi_q'$. By further blowing-ups of $\bar \chi_q'$, we separate $\tilde S_{a,q}$ from the two separatrices of $\bar\chi_q'$ unless one of them coincides with $\tilde S_{q,a}$. We will get in this way adapted simple singularities of (the transform of) $\tilde{\chi}_q'$ either reduced (c.f. Definition~\ref{def:adaptedsimplesing}-\textit{(1)}) or non reduced (c.f. Definition~\ref{def:adaptedsimplesing}-\textit{(2)}). When $\bar \chi_q'(q) \neq 0$, if $\Gamma$ is the formal solution of $\bar{\chi}'_q$ through $q$, a new blowing-up at $q$ produces an adapted simple singularity for the transform of $\bar \chi _q' $ at the point corresponding to the tangent line of $\Gamma$. If $\Gamma$ coincides with $\tilde{S}_{a,q}$, we get an adapted simple singular point for the transform of $\tilde \chi_q$. Otherwise, by further blowing-ups, we separate $\Gamma$ from ${S}_{a,q}$ and we get either adapted simple singularities or points in the situation already treated.

 In the process just described, we start blowing up a point $a\in \Singad (\hat \eta_J, \tilde F)$ which is a singular point of $\hat \eta_J'$, a point in a dicritical component where $\hat \eta_J'$ is not transversal, or a corner point. All these points correspond to characteristic cycles of $\wt{\mathcal{ M}}$ (they belong to $\wt{\mathcal{D}}$). Considering admissible blowing-ups $\sigma_{\gamma_{I}}$ centered at those $\gamma_I\in \wt{\mathcal{D}}$, the restriction $\sigma_{\gamma_{I}}\lvert_{\{\theta=0\}}$ is exactly the blowing-up of the corresponding singularity of the two-dimensional system. 
 Moreover, this property repeats for the subsequent points to be blown up to achieve $\tau_a$ and the corresponding strict transform of $\hat \eta_J'$.
 In other words, having defined the sequence of blowing-ups $\tau_a$ as above, satisfying (a) and (b) for any $a\in \Singad(\hat \eta_J',\tilde F)$ and for any $J\in \wt{J}_0$, the composition of these sequences of two dimensional blowing-ups provides a sequence of admissible blowing-ups $\mathcal{ M}=(M,\pi, \mathcal A, \mathcal D)$ factorizing through $\wt{\pi}$ (i.e. $\pi=\wt{\pi}\circ \pi'$) such that $\mathcal{ M}$ satisfies \textit{(2)} and \textit{(3)} of the statement. Since $\pi'$ does not modify the characteristic singularities $\g_{I_{-\infty}^{\widetilde {\mathcal{M}}}}, \g_{I_{-\infty}^{\widetilde {\mathcal{M}}}}$, we have also \textit{(1)}, and we are done.
\end{proof}
\begin{remark}\label{rmk:non-dicriticalreso}
	Notice that after an adapted reduction of singularities, the non-corner characteristic cycles that we obtain are contained in non-dicritical components of the total divisor.
\end{remark}

\subsection{Analytic approximations of normal forms}\label{sec:approximationandjets}
In this section, we study the effect of sequences of admissible blowing-ups to the analytic approximations $\xi_{\ell}$ of the formal normal form $\hat \xi$ for convenient values of $\ell$. First, we stablish the jet dependence of the transform of $\hat \xi$ by such blowing-ups in the different charts.

\begin{prop}\label{prop:ljets}
	Let $\hat \xi$  be a formal normal form of $\xi\in \mathcal H_3$. Consider an admissible sequence of blowing-ups $\mathcal{M}=(M,\pi,\mathcal A, \mathcal D)$ for $\hat \xi $ of length $l>0$, with $\mathcal A=\{ C_J \}_{J\in \mathcal J}$. For every $J\in \mathcal J$ and for every $k\geq 1$, if $x$ is a coordinate of the chart $C_J$ such that $\{ x=0  \}\subset E=\pi^{-1}(0)$, then we have $$j_{k}^x(\hat \xi^{(J)})=j_k^x((\pi\lvert_{C_J})^{\ast} j_{k+l+1} (\hat \xi)).$$

\end{prop}
\begin{proof}
	The proof of this result is a consequence of the following:\\
\begin{fact}
	 Let $\eta$ be a $K$-derivation with coefficients in $K[[x_1,\ldots, x_n]]$ and let $\tau$ be a quadratic morphism of the form $\tau(x_1,\ldots, x_n)=(x_1x_i,\ldots, x_{i-1}x_i, x_i, x_{i+1}x_i, \ldots, x_nx_i)$. Then, 
	\begin{equation}\label{eq:igualdadjetsexploencarta}
		\begin{split}
			j_k^{x_i}(\tau^\ast \eta)&=j_k ^{x_i} (\tau^\ast j_{k+1}(\eta)) \\
			j_k^{x_j}(\tau^\ast \eta)&= j_k ^{x_j} (\tau^\ast j_{k}^{x_j}(\eta)), \ j\neq i.
		\end{split}
	\end{equation}
\end{fact}
	We proceed by induction on the length $l$ of $\mathcal{ M}$. 	First consider $\mathcal{M}$ has length $l=0$ so that  $\pi= \sigma_{0}$ is the  blowing-up of the origin $0\in \mathbb{R}^3$ described in section~\ref{sec:firstblowing-up}. Then, we have (with simplified notation $\rho:=\rho^{(0)}, z:=z^{(0)}$):
		\begin{equation}\label{eq:igualdaddejetssigma0} \begin{split}
				j^\rho_{k}((\sigma_0\lvert_{C_0})^\ast \hat \xi)&=j^\rho_{k}( (\sigma_0\lvert_{C_0})^\ast   j_{k+1}  (\hat \xi)       ),\\
				j^z_k((\sigma_0\lvert_{C_\infty})^\ast \hat \xi)&= j^z_{k}( (\sigma_0\lvert_{C_\infty})^\ast   j_{k+1}  (\hat \xi)       ),\\
				j^z_k((\sigma_0\lvert_{C_{-\infty}})^\ast \hat \xi)&= j^z_{k}( (\sigma_0\lvert_{C_{-\infty}})^\ast   j_{k+1}  (\hat \xi)       ).
			\end{split}
		\end{equation}
	
	Suppose 
	 $l>0$ and that $\pi=\widetilde  \pi \circ \sigma_{\gamma_I}$ where $\sigma_{\gamma_I}$ is the admissible blowing-up centered at $\gamma_I$ for some characteristic element $\gamma_I$ of a sequence of admissible blowing-ups $\widetilde{ \mathcal{M}}=(\widetilde M,\widetilde \pi,\widetilde{ \mathcal A},
	\widetilde {\mathcal D})$ of length $l-1$. 
	 It is enough to study the transform in the charts $C_J$ when $\sigma_{\gamma_I} ^{-1}(\gamma_I)\cap C_J\neq \emptyset $, since the map $\sigma_{\gamma_I}$ is a diffeomorphism out of $\sigma_{\gamma_I} ^{-1}(\gamma_I)$. According to the construction of $\mathcal{M}$ from $\widetilde{\mathcal{M}}$ and using the same notations as in section~\ref{sec:sucessive-blow-up}, there are several cases:
	\begin{enumerate}
		\item $\g_I$ is a characteristic singularity (for instance $I=I ^{\widetilde{\mathcal{ M}} }_{ \infty}$) and $J=I ^{{\mathcal{ M}} }_{ \infty}=(\infty,\stackrel{s}{\ldots},\infty)$. In this case $x=z^{(J)}$, is the only coordinate of the chart $C_J$ in the conditions of the statement. 
		We prove the result by applying $s$ times the first formula in~\eqref{eq:igualdadjetsexploencarta}. Notice that $s\leq l+1$, so we get, for $x=z^{(J)}$, that $j_{k}^x(\hat \xi^{(J)})=j_k^x((\pi\lvert_{C_J})^{\ast} j_{k+s} (\hat \xi))=j_k^x((\pi\lvert_{C_J})^{\ast} j_{k+l+1} (\hat \xi))$.
		\item The point $\g_I$ is the origin of a chart $(C_{J_I},(x^{(J_I)}, y^{(J_I)}, z^{(J_I)} ))$ of $\widetilde{\mathcal{ A}}$ where $z^{(J_I)}=0  $ is the equation of the divisor $\widetilde{E}\cap C_{J_I}$ and $\sigma_{\gamma_{I}}\lvert_{C_J}: C_J\to C_{J_I}$ has the same expression as \eqref{eq:cartaC0} for $\sigma_0$, considering coordinates $(\theta,z^{(J)},\rho^{(J)})$ for $C_J$ and with the obvious change of notation. Notice that in $C_J$ the two coordinates $x=\rho^{(J)}$ and $x=z^{(J)}$ are in the conditions of the statement. By the induction hypothesis, renaming $z=z^{(J_I)}$ for simplicity, we have, for any $k\geq 1$, that $j^{z}_{k}(\hat \xi ^{(J_I)})= j^{z}_{k} ((\tilde \pi\lvert_{C_{J_I}})^\ast j_{k+l} (\hat \xi ) )$ and by the fact that $j_k(\hat \xi ^{(J_I)})=j_k(j_k^z(\hat \xi ^{(J_I)}))$, we also have $j_{k}(\hat \xi ^{(J_I)})= j_{k} ((\tilde \pi\lvert_{C_{J_I}})^\ast j_{k+l} (\hat \xi ) )$. The result follows from this similarly to the case of the first blowing-up $\sigma_0$ for $x=\rho^{(J)}$. The case $x=z^{(J)}$ is a consequence of the second equation of~\eqref{eq:igualdadjetsexploencarta}.

		\item 
		${\gamma_I}$ is a characteristic cycle of $\widetilde{ \mathcal{M}}$. Taking into account Remark~\ref{rmk:ciclosendoscartas}, we may assume  $\gamma_I\subset \{  \rho^{( J_I)}=0   \}$ inside the domain of a chart $(C_{J_I},(\theta, z^{(J_I)},\rho ^{(J_I)}))\in \tilde {\mathcal A}$.
		Let us put for simplicity $(z,\rho)=(z^{(J_I)},\rho^{(J_I)})$ and denote $\hat \xi^{(J_I)}=(\pi\lvert_{C_{(J_I)}})^\ast \hat \xi$. We distinguish two cases:
		\begin{enumerate}
		\item $\gamma_I$ is a corner characteristic cycle. In this case, $\sigma_{\gamma_I}^{-1}(\gamma_I)$ is covered by two charts $(C_{J}, (\theta,z^{(J)}, \rho^{(J)}))$ with $J=J_0,J_\infty$, for which the expression of $\sigma_{\gamma_I}$ is given by~\eqref{eq:cartaC00} and~\eqref{eq:cartaC0infty}, respectively. By symmetry, both are treated similarly, and we only develop the case $J=J_\infty$. Notice that the coordinates $x=z^{(J)}$ and $x=\rho^{(J)}$ are in the condition of the statement. For $x=z^{(J)}$, we have, for any $k\geq 1$:
		{\small\begin{equation}\label{eq:propjetgeneral2}
				\begin{split}
					j_k^x(\hat \xi^{(J)})&= j_k^x((\sigma_{\gamma_I}\lvert_{C_J})^{\ast} j_{k+1}(\hat \xi^{(J_I)}))= j_k^x((\sigma_{\gamma_I}\lvert_{C_J})^{\ast} j_{k+1}(j_{k+1}^z(\hat \xi^{(J_I)}))) \\
					&= j_k^x( (\sigma_{\gamma_I}\lvert_{C_J})^{\ast}  j_{k+1} ( j_{k+1}^z ( (\widetilde \pi\lvert_{C_{J_I}})^\ast (j_{k+l+1}(\hat \xi))  ) )  )= j_k^x( (\sigma_{\gamma_I}\lvert_{C_J})^{\ast}  j_{k+1} ((\widetilde \pi\lvert_{C_{J_I}})^\ast (j_{k+l+1}(\hat \xi))  )   ),\\ 
					&  =   j_k^x( (\sigma_{\gamma_I}\lvert_{C_J})^{\ast}   (\widetilde \pi\lvert_{C_{J_I}})^\ast (j_{k+l+1}(\hat \xi))     )   = j_k^x((\pi\lvert_{C_J})^\ast (j_{k+l+1}(\hat \xi))),
				\end{split}
		\end{equation}}
		where we have used the first formula of equation~\eqref{eq:igualdadjetsexploencarta} for the quadratic map $\sigma_{\gamma_{I}}$ in the first and fifth equality, general properties of jets presented in section~\ref{sec:formalnormal} in the second equality and the induction hypothesis in the third equality. This proves the proposition for $x=z^{(J)}$.
		For $x=\rho^{(J)}$, we have, for any $k\geq 1$:
		{ \small
			\begin{equation}\label{eq:propjetsgeneral1}
				\begin{split}
					j_k^x(\hat \xi^{(J)})&= j_k^x((\sigma_{\gamma_I}\lvert_{C_J})^{\ast} j_{k}^\rho(\hat \xi^{(J_I)}))= j_k^x((\sigma_{\gamma_I}\lvert_{C_J})^{\ast} j_{k+1}^\rho(\hat \xi^{(J_I)})) \\
					&= j_k^x( (\sigma_{\gamma_I}\lvert_{C_J})^{\ast}  j_{k+1}^\rho (  (\widetilde \pi\lvert_{C_{J_I}})^\ast (j_{k+l+1}(\hat \xi))  )   )= j_k^x( (\sigma_{\gamma_I}\lvert_{C_J})^{\ast}   (\widetilde \pi\lvert_{C_{J_I}})^\ast (j_{k+l+1}(\hat \xi))     ),\\ 
					&= j_k^x((\pi\lvert_{C_J})^\ast (j_{k+l+1}(\hat \xi)))
				\end{split}
		\end{equation}}
		where we have used the second formula of~\eqref{eq:igualdadjetsexploencarta} for the quadratic map $\sigma_{\gamma_{I}}$ in the first and forth equality and the induction hypothesis in the third equality. This proves the proposition for $x=\rho^{(J)}$.
		\item $\g_I$ is not a corner characteristic cycle. In this case $\sigma_{\gamma_{I}}^{-1}(\g_I)$ is covered by three charts $C_{J_{\infty}}, \ C_{J_0}$ and $C_{J_{-\infty}}$, for which the expression of $\sigma_{\gamma_{I}}$ is given by equations~\eqref{eq:cartaCIinfty}, \eqref{eq:cartaCI0} and \eqref{eq:cartaCImenosinfty}, respectively. In the chart $(C_{J_\infty},(\theta, z ^{(J_\infty)}, \rho^{(J_\infty)}))$, the two coordinates $x=z^{(J_\infty)}$ and $x=\rho^{(J_\infty)}$ are in the hypothesis of the statement. The proof of the result is analogous to the one above, namely equations~\eqref{eq:propjetgeneral2} and~\eqref{eq:propjetsgeneral1}. The chart $C_{J_{-\infty}}$ is similar to $C_{J_\infty}$. Finally, in the chart $(C_{J_0},(\theta,z^{(J_0)},\rho^{(J_0)}))$ only $x=\rho^{(J_0)}$ is in the hypothesis of the statement. The proof for this coordinate is just the same sequence of equalities as in~\eqref{eq:propjetgeneral2}, considering $x=\rho^{(J_0)}$ and with the interchange of the role of the coordinates $z$ and $\rho$ in $C_{J_I}$.
		\end{enumerate}
	\end{enumerate}
\end{proof}
Recall that the truncated normal forms $\xi_{\ell}$ of $\xi$ satisfy $j_\ell(\xi_\ell)=j_\ell(\hat \xi)$. Taking this into account, the same proof as above shows: 
\begin{prop}\label{cor:relationjetsabove}
	Let $\mathcal{M}=(M,\pi,\mathcal A, \mathcal D)$ be an admissible sequence of blowing-ups of length $l$ with $\mathcal A=\{ C_J \}_{J\in \mathcal J}$. Then, for $\ell\geq l+1$ and $J\in \mathcal J$ the transform $\xi_{\ell}^{(J)}:=(\pi\lvert_{C_J})^\ast \xi_{\ell}$ is analytic. Moreover, if $k\in \mathbb{N}$, $x$ is a coordinate of $C_J$ such that $\{ x= 0\}\subset E=\pi^{-1}(0)$ and $\ell\geq k+l+1$, then, we have
	\begin{equation}
		j_k^x(\xi_\ell^{(J)})=j_k^x(\hat \xi^{(J)}).
	\end{equation}
\end{prop}
In fact, we need a more precise result. To state it, recall the definition of the associated two dimensional vector fields $\hat \eta_J$ to $\hat \xi^{(J)}$ for $J\in \mathcal A\setminus \{ I_\infty^{\mathcal{M}},I_{-\infty}^{\mathcal{M}}  \}$ and the corresponding reduced vector fields $\hat \eta_J'=(\rho^{n_1^{(J)}} z^{n_2^{(J)}})^{-1} \hat \eta_J$, where $(\theta,z,\rho)$ are the coordinates in $C_J$. Let us write the transform $\xi_\ell^{(J)}$ for $\ell\geq l+1$ as
\begin{equation}
	\xi_\ell^{(J)}= B_{\ell,\rho}^{(J)}(\theta,z, \rho)  \rhoo + B_{\ell,\theta }^{(J)}(\theta,z, \rho) \thetaa+ B_{\ell,z}^{(J)}(\theta,z, \rho) \z
\end{equation}
The \emph{associated (to $\xi_\ell^{(J)}$) system of ODEs} $\eta_{\ell,J}$ is defined as:
\begin{equation}
	\label{eq:analytic2dimgeneral}
	\left\{ \begin{array}{lcll}
		\frac{d z}{d \theta} &=&  B_{\ell,\rho}^{(J)}(\theta,z, \rho)\cdot (B_{\ell, \theta}^{(J)}(\theta,z, \rho)) ^{-1}
		\\
		\frac{d \rho}{d \theta} &=& B_{\ell,z}^{(J)}(\theta,z, \rho)\cdot (B_{\ell, \theta}^{(J)} (\theta,z, \rho) )^{-1} 
	\end{array} \right. 
\end{equation}
\begin{remark}\label{rmk:characteristiccyclesarecycles}
	Since $n_1^{(J)}$ is always greater than 0, we have that every set $\{ \rho^{(J)}=0, z={\omega}   \}$ for any $\omega\in \mathbb{R}$ is a cycle for the vector field $\xi_{\ell}^{(J)}$, in particular, the characteristic cycles of $\hat\xi^{(J)}$ are cycles of $\xi_{\ell}^{(J)}$. For this reason, they are called \emph{characteristic cycles of $\xi_{\ell}^{(J)}$}.
\end{remark}

Recall that if $J\in  \{ I_\infty^{\mathcal{M}},I_{-\infty}^{\mathcal{M}}  \}  $ and we use simplified notation $(x,y,z):=(x^{(J)},y^{(J)},z^{(J)})$ for the coordinates of $C_J$, we have defined that $n^{(J)}$ as the maximum $n\in \mathbb{N}$ such that $\hat \xi ^{(J)} (z)$ is divisible by $z^n$. If $J\in \mathcal{J}\setminus  \{ I_\infty^{\mathcal{M}},I_{-\infty}^{\mathcal{M}}  \}$, we have defined $n^{(J)}:=\max \{  n_1^{(J)}, n_2^{(J)} \}$. With those notations, we have the following Corollary of Proposition~\ref{cor:relationjetsabove}.
\begin{cor}\label{cor:elem}
	Let $\mathcal{M}=(M,\pi,\mathcal A, \mathcal D)$ be an admissible sequence of blowing-ups of length $l>0$ with $\mathcal A=\{ C_J \}_{J\in \mathcal J}$. Define  $\ell_{\mathcal{ M}}:=  \max \{ n^{(J)}: J\in \mathcal{J}    \} +l+1$.
	\begin{enumerate}
		\item Let $k\in \mathbb{N}$ and $J\in \mathcal{J}\setminus \{  I_{-\infty} ^{\mathcal M},  I_{\infty} ^{\mathcal M} \}$. For every $\ell\geq \ell_{\mathcal{ M}} + k$, the monomial $(\rho ^{(J)})^{n_1^{(J)}}(z ^{(J)})^{n_2^{(J)}}$ divides the system $\eta_{\ell,J}$. Moreover, put $ \eta_{\ell, J}':=(\rho^{n_1^{(J)}}z^{n_2^{(J)}})^{-1} \eta_{\ell, J}$. Then, if 
		$x$ is a coordinate with $\{  x=0 \}\subset E\cap C_J$, then 
		$$j^x_{k}(\eta_{\ell, J}')=j^x_{k}(\hat \eta_{ J}').$$
		\item 
		Let $k\in\mathbb{N}$ and $J\in  \{  I_{-\infty} ^{\mathcal M},  I_{\infty} ^{\mathcal M} \}$. For every $\ell\geq \ell_{\mathcal{ M}} + k$, the function $\xi_{\ell}^{(J)}(z)$ is divided by $z^{n^{(J)}}$, and $$ j^z_{k}\left({z^{-n^{(J)}}} \xi_\ell^{(J)}(z^{(J)})\right) = j^z_{k}\left ({z^{-n^{(J)}}} \hat \xi^{(J)} (z^{(J)})\right ).   $$

	\end{enumerate}
\end{cor}
\begin{proof}
	Both statements are direct consequence of the jet equality stated in Proposition~\ref{cor:relationjetsabove}. Since $k+\ell_{\mathcal{ M}}\geq  n_i^{(J)} + l +1$ for $i=1,2$ and for very $J\in \mathcal{J}\setminus \{  I_{-\infty} ^{\mathcal M},  I_{\infty} ^{\mathcal M} \}$ and $k+\ell_{\mathcal{ M}}\geq n^{(J)}+ l+1$ when $J\in \{  I_{-\infty} ^{\mathcal M},  I_{\infty} ^{\mathcal M} \}$, we have that the monomials of type $(\rho ^{(J)})^{n_1^{(J)}}(z ^{(J)})^{n_2^{(J)}}$ divide the system $\eta_{\ell, J}$ when $J\in \mathcal{J}\setminus \{  I_{-\infty} ^{\mathcal M},  I_{\infty} ^{\mathcal M} \}$, or $(z ^{(J)})^{n^{(J)}}$ divides $\xi_{\ell}^{(J)}(z^{(J)})$ when $J\in  \{  I_{-\infty} ^{\mathcal M},  I_{\infty} ^{\mathcal M} \}$.
\end{proof}

Recall from equation~\eqref{eq:conjugationpsiell} that there is a formal automorphism $\psi_\ell$ at $0\in \mathbb{R}^3$ that conjugates $\hat{\xi}$ and $\xi_{\ell}$, that is, $\hat \xi=\psi^\ast_\ell (\xi_{\ell})$, and that it is tangent to the identity up to order $\ell$, i.e. $j_\ell(\psi_\ell-Id)=0$. In the following result, we prove that we can lift this conjugation to any chart of a sequence of admissible blowing-ups $\mathcal{M}=(M,\pi,\mathcal A, \mathcal D)$. 

\begin{prop}\label{prop:lifting}
	For any $\ell\geq l+1$ and any $C_{J}\in \mathcal A$, the automorphism $\psi_\ell$ can be lifted to an automorphism that conjugates $\xi_{\ell}^{(J)}=(\pi\lvert_{C_J})^\ast \xi_{\ell}$  and $\hat{\xi}^{(J)}=(\pi\lvert_{C_J})^\ast \hat \xi$. More precisely, suppose that $\ell=k+l+1$ with $k\geq 1$. We have:
	\begin{enumerate}
		\item Suppose $J=I^\mathcal{ M}_{\pm \infty}$ and denote $(x,y,z)=(x^{(J)},y^{(J)},z^{(J)})$ then, there is a tuple $\psi_\ell^{(J)}$ in $\R [x,y] [[z]]^3$ that satisfies $j^z_k(\psi_\ell^{(J)}- (x,y,z))=0$ and such that:
		\begin{equation}
			\pi\lvert_{C_J} \circ \psi_\ell^{(J)}=\psi_\ell  \circ \pi\lvert_{C_J}.
		\end{equation}
		\item Suppose $J\neq I^\mathcal{ M}_{\pm \infty}$ and denote $(\theta, z,\rho)=(\theta, z^{(J)}, \rho^{(J)})$ where we assume that $\{ x=0 \}\subset E\cap C_J$ for $x=z$ or $x=\rho$. Then, there is a tuple $\psi_\ell^{(J)}=\psi_\ell^{(J)}(\theta,z,\rho)=(\theta + F_\theta, z+F_z,\rho+ F_\rho)$, satisfying $j_{k}^x(\psi_{\ell}^{(J)}-(\theta,z,\rho))=0$ and such that:
		\begin{equation}
			\pi\lvert_{C_J} \circ \psi_\ell^{(J)}=\psi_\ell  \circ \pi\lvert_{C_J}, 
		\end{equation}
		where each $F_i\in \R [\cos \theta, \sin \theta, z][[\rho]]$ if $x=\rho$ for $F_i\in \R[\cos\theta,\sin \theta,\rho][[z]]$ or $x=z$.
	\end{enumerate}
\end{prop}
\begin{proof}
	We start observing that $\xi_\ell$ and $\hat\xi$ are formally conjugated at $0\in \Rtres$, as shown in section~\ref{sec:formalnormal}. 
	Moreover, since $j_\ell(\varphi_\ell)=j_\ell(\hat \varphi)$, we have that $j_{\ell}(\psi_\ell)=Id$. We need to lift this formal conjugation to a neighborhood of $E\cap C_J$ between $\xi_{\ell}^{(J)}$ and $\hat \xi^{(J)}$. We can write $\pi=\sigma_0\circ \sigma_{1} \circ \cdots \circ \sigma_{r'} \circ \cdots \circ \sigma_{r} $ with $\sigma_i=\sigma_{\gamma_{I_i}}$ and $0\leq r'\leq r$, where $\gamma_{I_i}$ are points of the form $\gamma_{(\varepsilon\infty,\stackrel{i}{\ldots},\varepsilon\infty)}$ for $0\leq i \leq r'$ and $\varepsilon=\pm1$, and $\gamma_{I_{j}}$ are characteristic cycles for $1+r'\leq j\leq r$. We can also suppose that each $\gamma_{I_{i+1}}\subset \sigma_i ^{-1}(\gamma_{I_i})$, since otherwise $\sigma_{\gamma_{I_{i+1}}} ^\ast$ acts as the trivial automorphism.
	
	Suppose that $r'\geq 1$. First, the conjugation $\psi_\ell$ can be lifted to  $\psi_\ell^{(\infty)}$ at the point $\gamma_\infty$ (or correspondingly, to $\psi_\ell^{(-\infty)}$ at the point $\gamma_{-\infty}$): using the chart $(C_\infty, (x^{(\infty)}, y^{(\infty)} , z^{(\infty)}  ))$ and the quadratic expression of $\sigma_0\lvert_{C_\infty}$, the formal automorphism defined by 
	\begin{equation}\label{eq:proofliftofpsiellinfty}
		\psi_\ell^{(\infty)} (x^{(\infty)}, y^{(\infty)} , z^{(\infty)}  )=  \left(  \frac{x\circ \psi_\ell}{z\circ \psi_\ell},  \frac{y\circ \psi_\ell}{z\circ \psi_\ell} ,z\circ \psi_\ell    \right) \circ \sigma_0\lvert_{C_\infty} (x^{(\infty)}, y^{(\infty)} , z^{(\infty)}  ) 
	\end{equation}
	satisfies that $\sigma_0 \lvert_{C_\infty} \circ \psi_\ell^{(\infty)}= \psi_\ell \circ \sigma_0 \lvert_{C_\infty}  $. Moreover, using that $j_{\ell}(\psi_\ell)=Id$ and the explicit expression of $\psi_\ell^{(\infty)}$, we get that $j_{\ell-1}(\psi_\ell^{(\infty)})=Id$. Recursively, we obtain that there is a formal automorphism $\psi_\ell^{(\infty,\stackrel{r'}{\ldots}, \infty)}$ at $\gamma_{I_{r'}}$ such that $(\sigma_0\circ \ldots \circ \sigma_{\infty,\stackrel{r'-1}{\ldots}, \infty)})\lvert_ {C_{(\infty,\stackrel{r'}{\ldots}, \infty)}}  \psi_\ell^{(\infty,\stackrel{r'}{\ldots}, \infty)}= \psi_\ell \circ (\sigma_0\circ \ldots \circ \sigma_{\infty,\stackrel{r'-1}{\ldots}, \infty)})\lvert_ {C_{(\infty,\stackrel{r'}{\ldots}, \infty)}} $ and $j_{\ell-r'}(\psi_\ell^{(\infty,\stackrel{r'}{\ldots}, \infty)})=Id$. 
	
	To finish, we need to check that $\psi_\ell^{(\infty,\stackrel{r'}{\ldots}, \infty)}$ can be lifted to a neighborhood of any $\gamma_{I_i}$ for $r'+1\leq i\leq r$. Renaming the point $\gamma_{I_{r'}}$ as the origin when $r'\geq1$, we can assume that $r'=0$, so that $\psi_\ell^{(\infty,\stackrel{r'}{\ldots}, \infty)}=\psi_\ell$ and $\sigma_{r'}$ is the first blowing-up $\sigma_0$ given in the covering of the chart $C_0$ as $\sigma_{0}(\theta, z,\rho)=(\rho\cos \theta, \rho\sin \theta, \rho z)$ if $\hat \xi$ has isolated singularity.

	Write 
	$\psi_\ell(x,y,z)=(x+G_1,y+G_2,z+G_3)$ where each $G_i\in \R[[x,y,z]]$ has order at least $\ell+1$. 
	Write first that, omitting super-indices in $(z^{(0)}, \rho^{(0)})$, 
	\begin{equation}
		\begin{array}{lcl}
			\psi_\ell \circ \sigma_0\lvert_{C_0}(\theta, z,\rho)&=&(\rho \cos\theta + G_1(\rho\cos \theta, \rho\sin \theta,\rho z),\rho \cos\theta + G_2(\rho\cos \theta, \rho\sin \theta,\rho z),\\&&\rho z + G_3(\rho\cos \theta, \rho\sin \theta,\rho z))
		\end{array}
	\end{equation}
	then $G_i\circ \sigma_0\lvert_{C_0} \in \R[\cos\theta, \sin \theta,z ][[  \rho ]] $ and $\rho^{\ell+1}$ divides each $G_i\circ \sigma_0\lvert_{C_0}$. We introduce $\tilde{G}_i=\frac{1}{\rho}G_i\circ \sigma_0\lvert_{C_0} $.	
	We look for a formal automorphism of the form $$\psi_\ell^{(0)}(\theta, z,\rho)=(\theta+ \rho F_1(\theta, z,\rho), z+\rho F_2(\theta, z,\rho),\rho+ \rho F_3(\theta, z,\rho))$$ such that $\sigma_0\lvert_{C_0}\circ \varphi^{(0)}=\varphi \circ \sigma_0\lvert_{C_0}$. The series $F_1,F_2,F_3$ must fulfill:
	\begin{equation*}
		\begin{split}
			\rho \cos \theta + \rho \tilde G_1&=(\rho +  \rho F_3 ) \cos (\theta+ \rho F_1), \\
			\rho \sin \theta +\rho \tilde  G_2&=(\rho + \rho F_3) \sin (\theta+\rho  F_1), \\
			\rho z + \rho \tilde G_3&=(z + \rho  F_2 ) (\rho +\rho  F_3 ).
		\end{split}
	\end{equation*}

	Put $\tilde F_1:=\rho F_1$ and $\tilde F_2=\rho F_2$. Using classical formulas for trigonometric functions, and dividing the above expressions by $\rho$, we obtain the following system 
	with coefficients in the ring $\Rcossinzrho$ and in the unknowns $\tilde F=(\tilde F_1, \tilde F_2, F_3)$:
	\begin{equation*}
		\begin{split}
			\tilde G_1&= 	 \cos\theta   F_3- \sin\theta \tilde F_1- \sin \theta F_3 \tilde F_1 +  O(\tilde F_1^ 2) ,	\\
			\tilde G_2&= 	  \sin\theta F_3+ \cos \theta \tilde F_1+ \cos \theta F_3 \tilde F_1 +  O(\tilde F_1^ 2), \\
			\tilde 	G_3 &=	 z F_3 +\tilde  F_2 +  F_3\tilde  F_2.	
		\end{split}
	\end{equation*}
	The differential of the system with respect to $\tilde F$ is invertible, as a matrix with entries in $\Rcossinzrho$. We apply the  implicit function theorem to find a solution $\tilde F\in \Rcossinz [[\rho]]^3$. 
	Since $j_{\ell-1}^\rho(\tilde G_i)=0$ for $i=1,2,3$, 
	we find $\rho^{\ell-1}$ divides $\rho F_1, \rho F_2,F_3$ and thus, $j_{\ell-1}^\rho(F_i)=0$ and $j_\ell(\varphi^{(0)})=Id$. We get the desired lifting of $\psi_\ell$ to $\psi_\ell^{(0)}$ in the 
	chart $C_0$.
	
	We proceed studying the rest of the blowing-ups by recurrence, considering that these blowing-ups are quadratic, similarly to the first case.
	 We make the first step of the recurrence. Up to a translation in $\{\rho=0\}$, the expression of the blowing-up (omitting super-indices) is either $\sigma_{{I_i}}\lvert_{C_J}(\theta, \rho,z)=(\theta, \rho,\rho z)$ (in the non-corner chart) or $\sigma_{{I_i}}\lvert_{C_J}(\theta, \rho,z)=(\theta, z \rho, \pm z)$ (in a corner chart). We obtain the desired expression of $\psi^{(J)}_\ell$ proceeding as in~\eqref{eq:proofliftofpsiellinfty}. Suppose, for example, that we study the situation in the corner chart:
	$$\psi^{(J)}_\ell (\theta,z,\rho)=(\theta\circ \psi_{\ell}^{(0)},z\circ\psi_{\ell}^{(0)},\frac{\rho \circ \psi_{\ell}^{(0)}}{z\circ \psi_{\ell}^{(0)}} )\circ \sigma_{{I_i}}\lvert_{C_J}. $$
	 By the expression of $\sigma_{{I_i}}$, we have that the coefficients of $\psi_\ell^{(J)}$ belong to $\R [\cos \theta, \sin \theta, \rho][[z]]$. Since coefficients of $\psi_\ell^{(0)}$ belong to $\R [\cos \theta, \sin \theta, z][[\rho]]$ and the above explicit expression, we get that the coefficients also belong to $\R [\cos \theta, \sin \theta, z][[\rho]]$. By the fact that $j_{\ell-1}^\rho( \psi_{\ell}^{(0)} - (\theta,z,\rho)  )=0$ and the above expression, we deduce $j_{\ell-2}^x(\psi_{\ell}^{(J)}-(\theta,z,\rho))=0$ for $x=\rho,z$.
\end{proof}
It is possible to define an analytic conjugation $\psi_{\ell,\ell'}$ between $\xi_\ell$ and $\xi_{\ell'}$ for any $\ell'\geq \ell$. And, since $j_\ell(\varphi_{\ell'})=j_\ell(\psi_{\ell'})$, $\psi_{\ell,\ell'}$ is an analytic diffeomorphism tangent to the identity up to order $\ell$. In fact, this diffeomorphism can be lifted to any chart after a sequence of admissible blowing-ups when $\ell,\ell'\geq \ell_{\mathcal{ M}}+1$.
\begin{cor}\label{cor:psiellell'}
	With the above notations, there is an analytic lifting $\psi_{\ell,\ell'}^{(J)}$ of $\psi_{\ell,\ell'}$ to any chart $C_J$ that conjugates $\xi_{\ell}^{(J)}$ and $\xi_{\ell'} ^{(J)}$.
\end{cor}
The proof of this corollary is based on the idea that $j_{\ell'}(\varphi_{\ell'})=j_{\ell'}(\hat \varphi)$ and it is analogous to the proof of the lifting of the formal automorphism $\psi_{\ell}$.

\section{Characteristic cycles as limit sets} \label{sec:analyticcharcyc}\label{sec:noacum}
	In this section, we use the analytic approximations $\xi_\ell$ to the formal normal form $\hat \xi$ with an objective: proving that the characteristic elements of $\xi_\ell$ after a sequence of admissible blowing-ups $\mathcal M$ are the only limit sets of the family of local cycles of $\xi_\ell$ for $\ell$ large enough.
	
	Along this section, we fix a sequence of admissible blowing-ups $\mathcal{M}=(M,\pi,\mathcal A,\mathcal D)$ for $\hat \xi$,  with $\mathcal A=\{  C_J \}_{J\in \mathcal J}$ and $\mathcal D=\{  \gamma_I \}_{I\in \mathcal I}$. 
Denote by $E=\pi^{-1}(0)$ the total divisor of $\pi$. We define also the \emph{support of $\mathcal{D}$} as $\text{Supp}\mathcal D=\displaystyle \bigcup_{I\in \mathcal I} \gamma_I $. 
\begin{prop}\label{thm:accumulationlocus}
	Let $\ell\geq\ell_{\mathcal{ M}}+1$ and $W$ be a neighborhood of $\text{Supp}\mathcal D=\displaystyle \bigcup_{I\in \mathcal I} \gamma_I $. There is some neighborhood $U=U(W)$ of $0\in \Rtres$ such that $\pi^{-1}(\mathcal C_U(\xi_ \ell))\subseteq W$.
\end{prop}
To prove this result, we need to introduce new notation and a technical lemma. 
Consider the set  
	$$\dot{E}:=E\setminus \left( \displaystyle \bigcup_{I:\gamma_I\in \mathcal D \text{ corner}}  \{\gamma_I \} \cup \{\gamma_{I_\infty^{\mathcal M}}\} \cup \{\gamma_{I_{-\infty}^{\mathcal M}} \} \right).$$ It has a finite family of connected components denoted by $ \mathcal{E}_{\mathcal{M}}= \{  L_0, L_1, \ldots, L_{k_\mathcal{ M}}\}$. Each $L_i\in \mathcal{E}_{\mathcal{ M}}$ is open in $E$ and contained in a chart $C_{J_i}$ for $i=0,1,\ldots, k_{\mathcal M}$, therefore, we will call them \textit{open components (of $E$)}. In addition, in case $L_i$ is contained in two different charts, we choose $C_{J_i}$ in which $L_i\subseteq \{ \rho^{(J_i)}=0  \}$, which is always possible by the hypothesis (H3). Then, $L_i= \SSS^1\times (\lambda_i^-,\lambda^+_i)\times  \{0\}\subset C_{J_i}$ where $\lambda_i^-\in \R\cup \{-\infty\}$ and $\lambda_i^+\in \R\cup \{\infty\}$. An element $L_i\in \mathcal{E_M}$ is said to be dicritical (respectively, non-dicritical)if the component of $E$ that contains $L_i$ is dicritical (respectively, non-dicritical). 
	
	Fix $L=L_i\in \mathcal{E_M}$ and the corresponding chart $C_{J}$ with $J=J_i$. We explained in previous sections that the vector field $\hat \xi^{(J)}=(\pi\lvert_{C_{J}})^\ast\hat \xi$ in the chart $C_{J}$ has an associated two dimensional formal vector field $\hat \eta_J$ as in equation \eqref{eq:formal2dimgeneral}. For the purpose of this section, we write, removing super-indices in $(z,\rho)$:
		\begin{equation}\label{eq:formalsystemwithoutn2}
		\hat \eta_ J= \rho^{n_1^{(J)}} ( A^{(J)}_z(z,\rho)\z+ A^{(J)}_\rho(z,\rho)\rhoo)
	\end{equation}

	Notice that there is a small modification here with respect to equation \eqref{eq:formal2dimgeneral}: 
	we include $z^{n_2^{(J)}}$ in $A_j^{(J)}$ for $j=z,\rho$. However, we do not change the notation in this section. The reduced two dimensional vector field $\hat \eta_J'=\rho^{-n_1^{(J)}}\hat \eta_J$ has a finite number of adapted singularities along $\{ \rho=0 \}$, which determine the characteristic cycles contained in $L_i$. The
	$z-$coordinates of the characteristic cycles in $L$ are denoted by $\omega_1^{L}, \ldots, \omega^L_{m_{L}}$. We denote the associated characteristic cycles $\gamma^L_{1}, \ldots, \gamma^L_{m_L}$.
	
	 We define a collection of sets depending on two parameters $\varepsilon,\delta>0$: $$\mathcal V (L,\varepsilon,\delta):=  \{  V_0, V_1 ,\ldots, V_{m_{L_i}}, V_{m_{L_i} +1}  \},$$ which are given by:
	\begin{equation}\label{eq:boxesdefi}
		\begin{array}{lcl}
			V_{0} &=& \mathbb{S}^ 1 \times[\mu_-, \omega_1^{L}-\varepsilon] \times  (0,\delta],\\
			V_{j} &=&\mathbb{S}^ 1 \times [\omega_{j}^{L}+\varepsilon, \omega_{j+1}^{L}-\varepsilon] \times  (0,\delta] \quad j=1, \ldots, m_{L_i},\\
			V_{m_{L}+1} &=&\mathbb{S}^ 1 \times  [\omega_{m_L}^{L}+\varepsilon, \mu_+] \times (0,\delta],
		\end{array}
	\end{equation}
	where either $\mu_-=\lambda^-_i  + \varepsilon $ when $\lambda^-_i>-\infty$ or  $\mu_- =\omega_{1}^{L}-\frac{1}{\varepsilon}$ when $\lambda^-_i=-\infty$, and analogously for $\mu_+$. Define the surfaces $\partial_{min}V_j$ and $\partial_{max}V_j$ as follows:
	\begin{itemize}
		\item $\partial_{min}V_0= \SSS^1 \times \{   \mu_- \} \times [0,\delta]   $ and $\partial_{min} V_j= \SSS^1 \times \{   \omega_{j}^{L}+\varepsilon \} \times [0,\delta] $ for $j=1,2,\ldots,m_{L}$.
		\item $\partial_{max}V_j=  \SSS^1 \times \{   \omega_{j+1}^{L}-\varepsilon \} \times [0,\delta] $ for $j=0,1,\ldots,m_{L}-1$ and $\partial_{max}V_{m_{L}}= \SSS^1 \times \{   \mu_+ \} \times [0,\delta]   $
	\end{itemize}
	\begin{lem} \label{lem:boxesdef}
		
		Consider a connected component of the divisor $L\in \mathcal{E}_{\mathcal{ M}}$ and a chart $(C_{J},(\theta,z,\rho))$ such that $L\subset \{\rho=0\}$ as above. Assume $\ell\geq \ell_{\mathcal{ M}}+1$ and denote by $\xi_\ell^{(J)}=(\pi\lvert_{C_J})^\ast \xi_\ell$. There exists $\varepsilon_0>0$ such that for every small $\varepsilon$ with $\varepsilon_0>\varepsilon >0$, there exists $\delta=\delta(\varepsilon)>0$ such that $\mathcal V (L,\varepsilon,\delta):=  \{  V_0, V_1 ,\ldots, V_{m_{L}}, V_{m_{L} +1}  \}$ satisfies:

		\begin{enumerate}
			\item 
			In case $L$ is non-dicritical, the function $z$ is monotonic along the trajectories of $\xi_\ell^{(J)} $ in each $V_j$ for  $j \in \{0,1, \ldots ,m_{L}\}$. Otherwise, if $L$ is dicritical, the function $\rho$ is monotonic along the trajectories of $\xi_\ell^{(J)}=(\pi\lvert_{C_J}) ^\ast \xi_\ell $ in each $V_j$ for  $j \in \{0,1, \ldots ,m_{L}\}$.
	
	\item If $L$ is dicritical and $\rho^{n_1^{(J)}+1}$ does not divide $\xi_\ell^{(J)}(z)$, then $\xi_\ell^{(J)}(z)$ has constant sign along the surfaces $\partial _{min}V_j$ and $\partial _{max}V_j$.

	\item If $ L$ is dicritical and $\rho^{n_1^{(J)}+1}$ divides $\xi_\ell^{(J)}(z)$, denote $\mathcal V(L,\frac{\varepsilon}{2}, \delta)=\{  V_0',V_1',\ldots, V_{m_{L_i}+1}'  \}$. Then, each element $V'_j\in \mathcal V(L,\frac{\varepsilon}{2}, \delta)$ fulfills (1) and, moreover, any trajectory of $\pi^\ast \xi_\ell$ containing a point in $V_j$ remains inside $V_j'$ either for any positive time $t\geq 0$ or for any negative time $t\leq 0$.
		\end{enumerate}
	\end{lem}
Notice that $\mathcal V (L,\varepsilon,\delta')$ also fulfills (1-3) of the lemma for any $\delta'<\delta$.
\begin{proof}
		The vector field $\xi_\ell^{(J)}$ is described by a non-autonomous two dimensional system of ODEs (see equation \ref{eq:analytic2dimgeneral}), as follows:
	\begin{equation}
		\label{eq:systemODEsCI}
		\left\{ \begin{array}{lcll}
			\frac{d z}{d \theta} &=&  \rho^{n_1^{(J)}}  A^{\ell,(J)}_z (\theta,z, \rho),	\\
			\frac{d \rho}{d \theta} &=& \rho^{n_1^{(J)}}  A^{\ell,(J)}_\rho (\theta,z, \rho)
		\end{array} \right. 
	\end{equation}
	where $A^{\ell,(I)}_u(\theta,z, 0)=A^{(J)}_u(z,0)$ for $u=\rho,z.$ As for the formal system of ODEs \eqref{eq:formalsystemwithoutn2}, we include the factor $z^{n^{(J)}_2}$ in $A^{\ell,(I)}_u$.
	
	We choose $\varepsilon_0$ satisfying the following conditions:
	\begin{itemize}
		\item In any case, we require $\varepsilon_0< \tfrac{1}{2}\min_{i\neq j} \{  |\omega_i^L-\omega_j^L|      \}$.

		\item If L is dicritical we take $\varepsilon_0$ with an additional requirement. If ${A_z^{\ell,(I)}}(\theta,z,0)\neq 0$, we denote $t_1,\ldots, t_s$ its zeroes.
		Then, $\varepsilon_0$ must fulfill: $$\varepsilon_0< \frac{1}{2} \displaystyle \min \{ |\omega_j^L-t_{k}|\ :\  {1\leq j \leq m_L,\ 1\leq k \leq s, \ \omega^L_j\neq t_{k} }   \}.$$
		
	\end{itemize}

	In the non-dicritical case, since $\ell\geq \ell_{\mathcal{M}}+1$, the function $A_z^{\ell,(J)}(\theta,z,0)\equiv A_z^{(J)}(z,0) $, only depends on $z$ by Corollary~\ref{cor:elem}. Being its zeroes $\omega_1^{L}, \ldots, \omega^L_{m_{L}}$, the function $A_z^{\ell,(J)}(\theta,z,0)$ has constant sign when $z$ belongs to the interval of
	$\Omega_j(\varepsilon)=[\omega_j^{L}+\varepsilon, \omega_{j+1}^{L}-\varepsilon]$ for any $0<\varepsilon<\varepsilon_0$. By continuity and periodicity in $\theta$, 
	$A_z^{\ell,(J)}(\theta,z,\rho)$ has constant sign for $(\theta,z,\rho)$ in a set of the form $\mathbb{S}^1\times \Omega_j(\varepsilon)\times  (0,\delta_j]$ for some $\delta_j=\delta_j(\varepsilon)$. 
	Then, take $\delta$ fulfilling two properties: $\delta\leq \displaystyle \min_{i=0, \ldots, m_{L_i}+1} \{ \delta_i \}$ and $B_{\ell,\theta}^{(J)}$ has positive sign in $\mathbb{S}^1\times \Omega_j(\varepsilon)\times  (0,\delta]$ for every $j=0, \ldots, m_{L}+1$, which is possible since $B_{\ell,\theta}^{(J)}(\theta,0,0)=1$. With these considerations, we define $V_j=\mathbb{S}^1 \times \Omega_j(\varepsilon) \times  (0,\delta]$. 
	Taking into account that 
	$\xi_{\ell}^{(J)}(z)=\rho^{n_1^{(J)}} A^{\ell,(J)}_z(\theta,z,\rho) \cdot B_{\ell,\theta}^{(J)}( \theta,z, \rho)$, we obtain the property \textit{(1)}. 
	
	 In the dicritical case and when $\rho$ does not divide $A_z^{\ell,(J)}(\theta,z,0)\neq0$, we proceed in the same way, but we work with $A_\rho^{\ell,(J)}(\theta,z,0)$ instead of $A_z^{\ell,(J)}(\theta,z,0)$. Notice that, by the definition of $\ell_\mathcal{M}$, it is possible since $A_\rho^{\ell,(J)}(\theta,z,0)$ only depends on $z$ has the same zeros (namely, $\omega_1^L,\ldots,\omega_{m_L}^L$) as $A_\rho^{(J)}(z,0)$. We get that $\xi_\ell^{(J)}(\rho)$ has constant sign in each $V_j$. Then statement \textit{(1)} of the Lemma holds.
	
	Let us show \textit{(2)}. 
	By the choice of $\varepsilon_0$, we have that $
	 A_z^{\ell,(J)}
	 (\theta,z,0)$ does not vanish at any of the extreme values of $\Omega_j(\varepsilon)$. Since $\xi^{(J)}_\ell(z)=\rho^{n_1^{(J)}}  A^{(J)}_z(\theta,z,\rho) \cdot B_{\ell,\theta}^{(J)}( \theta,z, \rho)$, we obtain \textit{(2)}, up to taking a smaller $\delta$.
	
	Finally, we show \textit{(3)}. Suppose that $L$ is a dicritical open component and that
	 $A_z^{(J)}(\theta,z,0)\equiv 0$. Then, the associated to $\xi^{(J)}_\ell$ system~\eqref{eq:systemODEsCI} can be written as:
	\begin{equation}\label{eq:demlemacajas}
		\left\{      
		\begin{array}{lll}
						\frac{dz}{d\theta }&= & \rho^{n^{(J)}_1+1}\widetilde A_{z} ^{\ell,(J)}(\theta, z,\rho )\\
			\frac{d\rho}{d\theta }&= & \rho^{n^{(J)}_1}A_{\rho} ^{\ell,(J)}(\theta, z,\rho )
		\end{array}
		\right. ,
	\end{equation}
	where $A_{\rho} ^{\ell,(J)}(\theta, z,0 )$ does not depend on $\theta$ from Corollary~\ref{cor:elem}, vanishes for $z\in \{ \omega_1^L,\ldots,\omega_{m_L}^L  \}$, and $\widetilde{A}_{z} ^{\ell,(J)}( \theta,z,0 )\in \mathbb{R}[\cos\theta,\sin \theta,z]$. Proceeding as in the beginning of the proof, we take a constant $ \delta >0$ such that the collection $ \mathcal{V}(L,\frac{\varepsilon}{2},  \delta)=\{ V_0',V_1',\ldots, V_{m_L}'   \}$ fulfills \textit{(1)} so that $\rho$ is monotonic in every $V_j'$. 
	 The compactness of $\overline{V}_j$ implies that there are constants $a,K>0$  such that for any $V_j'\in\mathcal{V}(L,\frac{\varepsilon}{2},  \delta)$, we have 
	\begin{equation}
		\label{eq:inequalityprooflemmacajas}
	  \inf_{p\in V_j'}   \{   |A_{\rho} ^{\ell,(J)}(p) |    \} \geq a, \quad   \sup_{p\in V_j'}   \{   |\widetilde A_{z} ^{\ell,(J)}(p) |    \} \leq  K. 
   \end{equation}
	Fix $V_j'$ and suppose, for instance, that $A_{\rho} ^{\ell,(J)}\lvert_{V_j'}<0$. Then, if $\sigma:\mathbb{R}\longrightarrow M$ is a trajectory of $\xi_\ell^{(J)}$ parameterized as a solution $\sigma(\theta)=(\theta, z(\theta),\rho(\theta))$ of system \eqref{eq:demlemacajas}, as long as it remains in $V_j'\setminus L$, the function $\rho(\theta)$ is strictly decreasing. Hence, $\sigma$ can be parameterized by $\rho$ instead of $\theta$ and we obtain that
	$$ \left|  \frac{dz}{d\rho}   \right| \leq C\rho, \text{ where } C=\frac{K}{a}$$
	by \eqref{eq:demlemacajas} and \eqref{eq:inequalityprooflemmacajas}. Now, consider the collection $\mathcal{V}(L,\varepsilon,  \delta)=\{ V_0,V_1,\ldots, V_{m_L}  \}$ whose elements fulfill $V_j\subset V_j'$ for $j=0,1,\ldots, m_L$.  If the trajectory $\sigma$ starts at a point $p_0=(\theta_0,z_0,\rho_0)\in V_j\subset V_j'$ with $\rho_0>0$, it satisfies, for $\theta>\theta_0$:
	$$  |z(\theta)-z_0|\leq \frac{C}{2} |\rho(\theta)-\rho_0|^2\leq \frac{C}{2}\rho_0^2 \leq   \frac{C}{2} \delta^2 $$
	as long as $\text{Im}(\sigma\lvert_{[\theta_0,\theta]})\subset V_j'$. We obtain similar bounds for $|z(\theta)-z_0|$ when $A_{\rho} ^{\ell,(J)}\lvert_{V_j'}>0$. Imposing $\delta<\sqrt{\frac{\varepsilon}{C} }$, we can conclude that $|z_0-z(\theta)|<\frac{\varepsilon}{2}$ and guarantee, for any $j$ and for any $p_0\in V_j\in \mathcal V(L,\varepsilon,\delta)$, the trajectory $\sigma$ starting at $p_0$ satisfies $\text{Im}(\sigma\lvert_{[\theta_0,\infty)})\subset V_j'$ (or $\text{Im}(\sigma\lvert_{(-\infty,\theta_0]})\subset V_j'$ in case $A_{\rho} ^{\ell,(J)}\lvert_{V_j'}>0$).
\end{proof}

\begin{notation}
	With the above notations, we say that the box $V_j\in \mathcal V(L,\varepsilon, \delta)$ with $j=1,\ldots, m_L-1$ is \emph{adjacent} to the characteristic cycles $\{ z=z^L_{j}, \rho=0  \}$ and $\{ z=z^L_{j+1}, \rho=0   \}$ for $j=1,\ldots, m_{L}$. If, in the notations of section~\ref{sec:sucessive-blow-up}, these characteristic cycles are $\gamma_{I_{j}},\gamma_{I_{j+1}}$ with $I_j,I_{j+1}\in \mathcal{I}$, we denote $\partial_{I_j}V_j=\partial_{min}V_j$ and $\partial_{I_{j+1}}V_j=\partial_{max}V_j$. Similarly, we consider the corresponding corner characteristic cycles in $\bar{L}$, or points $\gamma_{I_{\infty}^{\mathcal M}}$ or $\gamma_{I_{-\infty}^{\mathcal M}}$, adjacent to $V_0$ and $V_{m_L}$.
\end{notation}

\begin{proof}[Proof of Proposition \ref{thm:accumulationlocus}.]
	Let $W$ be a neighborhood of $\text{Supp}\mathcal D$. For every $I\in \mathcal I$ we consider an open neighborhood $W_I\subset W$ of $\gamma_I$ such that $W_I\cap W_{I'}=\emptyset$ if $I\neq I'$. Consider the collection $\mathcal E_{\mathcal M}$, and apply Lemma~\ref{lem:boxesdef} to each $L_i\in \mathcal E_{\mathcal M}$, taking $\varepsilon$ and $\delta$  small enough so that each family $\mathcal V (L_i,\varepsilon,\delta)$ also satisfies: 
	\begin{itemize}
		\item For any $V\in \mathcal{V}(L_i,\varepsilon, \delta)$, we impose $V\cap W_I\neq \emptyset$ if and only if $\gamma_I$ is adjacent to $V$.
		\item The boundaries $\partial _{min}V$ and $\partial _{max} V$ are contained in the corresponding neighborhoods $W_{I_1}$ and $W_{I_2}$, where $\gamma_{I_1}$ and $\gamma_{I_2}$ are adjacent to $V$.
		\item The set
		\begin{equation*}
			\displaystyle \bigcup_{I\in \mathcal I} W_I  \cup \bigcup_{L\in \mathcal E_{\mathcal M}}    \bigcup_{V\in\mathcal V(L,\varepsilon ,\delta)} V
		\end{equation*}
		is a neighborhood of the divisor $E=\pi^{-1}(0)$ in $M$.
	\end{itemize}
	Now, we define a closed neighborhood $\wt{W}_I\subset W_I$ of $\gamma_I$ for each $I\in \mathcal{I}$ in such a way that (see Figure~\ref{fig:corsssectionprop-2}):
	\begin{enumerate}
		\item[(i)] The set
			\begin{equation*}
			\wt{U}=\text{int}\left(
			\displaystyle \bigcup_{I\in \mathcal I} \widetilde W_I  \cup \bigcup_{L\in \mathcal E_{\mathcal M}}    \bigcup_{V\in\mathcal V(L,\varepsilon ,\delta)} V\right)
		\end{equation*} 
	is a neighborhood of the divisor $E$ in $M$.
	\item[(ii)] For any $I\in \mathcal{I}$, $L\in \mathcal{E}_{\mathcal{ M}}$ and $V\in \mathcal{V}(L,\varepsilon, \delta)$, $\wt{W}_I\cap V$ is empty, in case $V$ is not adjacent to $\gamma_I$, or, otherwise, a closed non-empty subinterval of the corresponding boundary 
	$\partial_I V$, with one of the extremities at $\partial_I V\cap E$. This means that, if $\partial_IV=\SSS^1\times \{z_{\partial_I V}\}\times (0,\delta]$, then $\wt{W}_I\cap V=\SSS^1\times \{z_{\partial_I V}\}\times (0,\mu]$.
	\end{enumerate}

Now, the set $U:=\pi(\tilde U)$ is an open neighborhood of $0$ satisfying the requirements of the proposition. More precisely,
we claim that $\mathcal C_U\subset \displaystyle \bigcup_{I\in \mathcal I} \widetilde W_I   $.

\begin{figure}
	\centering
	\includegraphics[width=0.5\textwidth]{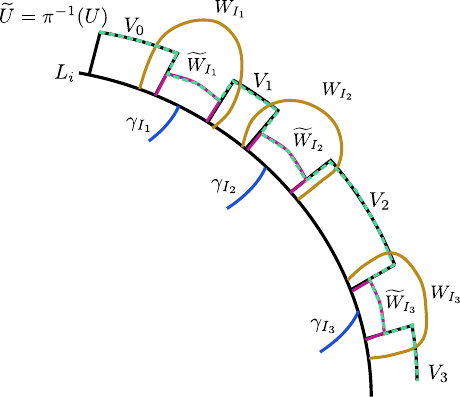}
	\label{fig:corsssectionprop}
	\caption{Cross-section of the neighborhoods $\widetilde W_{I_j}\subset W_{I_j}$ and of $\tilde{U}$.}
	\label{fig:corsssectionprop-2}
\end{figure}
	
	To see it, suppose that there is a cycle $Z$ of $\xi_\ell$ contained in $U$ and such that $\widetilde Z:= \pi^{-1}(Z)$ intersects some $V\in \mathcal V(L, \varepsilon, \delta)$ for some $L$. Consider a parametrization $\sigma:\R \longrightarrow \tilde U$ of $\widetilde Z$ as a trajectory of $\pi^\ast \xi_\ell$ such that $\sigma(0)\in V$. By the property \textit{(1)} of Lemma~\ref{lem:boxesdef}, one of the coordinates $z$ or $\rho$ is monotonic along $\sigma$  inside $V$, so it cannot be completely contained in $V$. 
	As a consequence, $\sigma$ leaves $V$ so that for some $t_0\geq 0$ we have $\sigma(t_0)\in \text{Fr}(V)\cap \wt{W}_I$, where $I\in \mathcal{I}$ and $\gamma_I$ is adjacent to $V$. By construction (cf. item (ii) above), $\sigma(t_0)$ belongs to the boundary $\partial _I V$. We have two cases to consider (notations as in Lemma~\ref{lem:boxesdef}). 
	\begin{itemize}
		\item $A_z^{\ell, (J)}(\theta,z,0)\neq 0$. By statement (2) of Lemma~\ref{lem:boxesdef}, the vector field $\pi^\ast \xi_\ell$ is transverse to $\partial_I V\setminus L$, so that, for instance, we have $\sigma((t_0-c,t_0))\subset \text{int}(V)$ and $\sigma((t_0,t_0+c))\subset \text{ext}(V)$ for some $c>0$. Since $\sigma$ is periodic, we must have that $\sigma$ crosses $\text{Fr}(V)$ at a first time $t_1>t_0$ necessarily along one of the boundaries $\partial_{min}V, \partial_{max}V $ where $\pi^\ast \xi_\ell$ points towards $\text{int}(V).$ 
		If we denote $\{  \partial_I V,\partial_{I'} V\}= \{ \partial_{min} V,\partial_{max} V\}$, we must have $\sigma(t_1)\in \partial_{I'} V$ and $\sigma((t_0,t_1))\subset \text{ext}(V)$, but this contradicts the fact that, by construction, $\tilde U\setminus V=\widetilde U_1 \cup \widetilde U_2$ where $U_1,\tilde U_2$ are non-empty open sets such that $\tilde U_1\cap \tilde U_2=\emptyset$ and the closure of each $\tilde U_i$ cuts $V$ along a unique interval among $\{\partial_{I'}V, \partial_{I} V\}$.
		\item $A_z^{\ell, (J)}(\theta,z,0)\equiv 0$. Using statement (3) of Lemma~\ref{lem:boxesdef}, we know that either $\sigma((t_0,\infty))$ or $\sigma((-\infty, t_0))$ is contained in the corresponding element $V'$ of the collection $\mathcal{V}(L, \frac{\varepsilon}{2},\delta)$ and $\rho\circ \sigma$ is monotonic along that interval. This is also a contradiction of $\sigma$ being periodic.
	\end{itemize}
	Consequently, we have proved that $\tilde Z\subset \displaystyle \bigcup_{I\in \mathcal I} \wt{W}_I$ (in fact, included in a single one $\wt{W}_I$ by connectedness). 
	Therefore, we have that:
	$$  \pi^{-1}(\cu)\subset  \pi^{-1}(U)\setminus \left( \displaystyle  \bigcup_{L\in \mathcal E_{\mathcal M}}    \bigcup_{V\in\mathcal V(L,\varepsilon ,\delta)} V \right) \subseteq  \bigcup_{I\in \mathcal I}\widetilde  W_I 
	\subseteq  \bigcup_{I\in \mathcal I} W_I\subseteq W,
	 $$
	 as we wanted to prove.

\end{proof}

\section{Analysis of final adapted simple singularities}
Along  this section, we consider some $\xi\in \mathcal{H}^3$, we fix a formal normal form $\hat \xi$ of $\xi$ and an adapted resolution of singularities $\mathcal{ M}=(M,\pi,\mathcal{A}, \mathcal{D})$ of $\hat{\xi}$ according to Proposition~\ref{prop:adaptedressing}. Denote by $E=\pi^{-1}(0)$ the exceptional divisor of $\mathcal M$.

\subsection{Infinitely near points of the rotational axis}

We see first that we can find a neighborhood of the two characteristic singular points that does not contain cycles of a jet approximation $\xi_{\ell}$ of $\hat \xi$. 
\begin{prop}\label{prop:zaxis} Suppose that $\xi$ has an isolated singularity at 0 and
consider $\xi_\ell$ with 
$\ell\geq \ell_{\mathcal{M}}+1$. There exist neighborhoods $W_\infty$ of $\gamma_{I_{-\infty}^{\mathcal M}}$ and $W_{-\infty}$ of $\gamma_{I_{\infty}^{\mathcal M}}$ in $M$ such that neither $W_\infty\setminus E$ nor $W_{-\infty}\setminus E$ contains cycles of $\pi^\ast \xi_{\ell}$.
\end{prop}
\begin{proof}
	According to the construction in Section~\ref{sec:sucessive-blow-up}, the point $\gamma_I=\gamma_{I_{\infty}^{\mathcal M}}$ is the origin of the chart $(C_{J},(x^{(J)} , y^{(J)}, z^{(J)}   ))$ with 
 $J=I_{\infty}^{\mathcal M}$ and $E\cap C_J=\{ z^{(J)}=0  \}$. 
 Being $\mathcal{M}$ and adapted resolution of singularities of $\hat \xi$ and by means of Proposition~\ref{cor:relationjetsabove}, we have in a neighborhood of $\gamma_{I_{\infty}^{\mathcal M}}$ that $\xi_\ell^{(J)}(z^{(J)})=(z^{(J)})^{n^{(J)}} \cdot F(x^{(J)} , y^{(J)}, z^{(J)})$ where $\xi_\ell^{(J)}=(\pi\lvert_{C_J})^\ast \xi_\ell$, $n ^{(J)}\in \mathbb{N}_{\geq 1}$ and $F(x^{(J)} , y^{(J)}, z^{(J)})\in \R[x^{(J)} , y^{(J)}][[ z^{(J)}]]$ converges and satisfies $F(0,0,0)\neq 0$. Take a neighborhood $W_\infty$ of $	\g_I$ in $M$ where $F$ has a constant sign, positive or negative. 
 We have that the trajectories of $\pi^\ast \xi_{\ell}$ in $W_\infty\setminus E$ can be parameterized by $z^{(J)}$, which avoids the existence of cycles of $\pi^\ast \xi_{\ell}$ in $W_\infty\setminus E$. The proof for $\gamma_{I_{-\infty}}^\mathcal{M}$ is analogous. 
 
\end{proof}
\subsection{Simple corner characteristic cycles}\label{sec:corner}
We prove that cycles of $\pi^\ast \xi_{\ell}$ cannot accumulate along corner characteristic cycles. Once again, the argument is to find a function, around such corner characteristic cycle, which is monotonic along the trajectories of $\pi^\ast \xi_{\ell}$, for $\ell$ sufficiently large. 
\begin{prop}\label{prop:corner}
Let $\ell\geq\ell_{\mathcal M}+1$. Consider a corner characteristic cycle $\gamma_I$ of $\hat \xi$ in $\mathcal M$. There exists a neighborhood $W_{I}$ of $\gamma_{I}$ in $M$ such that $\pi^\ast(\xi_{\ell})$ does not contain cycles in $W_{I}\setminus E$.
\end{prop}
\begin{proof}
	By construction, the corner characteristic cycle $\gamma_I$ is given by $\{  z ^{(J)}=\rho^{(J)} = 0 \}$ for $(C_J,(\theta, z ^{(J)}, \rho^{(J)})) \in \mathcal A$ for which $E\cap C_J=\{\rho^{(J)}z ^{(J)}=0  \}$. For simplicity, from now on, we remove the super-indices of the coordinates. By definition of $\pi$ being an adapted resolution of singularities, at least one of the two components of the divisor $\{  \rho=0  \}$ and $\{ z=0  \}$ is non-dicritical. More precisely, let $\hat \eta_J$ be the two dimensional vector field associated to $\hat\xi  ^{(J)}=(\pi\lvert_{C_J})^\ast(\xi_{\ell})$ and consider the strict transform $\hat \eta_J'=\frac{1}{\rho^{a}
		z^{b}} \hat \eta_J$. 
	Then, we have two cases:
	\begin{itemize}
		\item[a)] The origin is not a singular point of $\hat \eta_J'$ and one of the components, say $\{  z=0 \}$, is the solution of $\hat \eta_J'$. 
		\item[b)] The origin is a simple singularity of  $\hat \eta_J'$ and both components are invariant.
	\end{itemize}
	In the case a), write
	\begin{equation*}
		\hat \eta_J'= zF_z(z,\rho)\z + (\lambda_2 +  F_\rho(z,\rho))\rhoo
	\end{equation*}
	where $\lambda_2\neq 0$ and $F_\rho,F_z\in \R[\rho][[z]]\cap \R[z][[\rho]]$ with $F_\rho(0,0)= 0$.
	We have that
	\begin{equation}\label{eq:proofcorner1}
		\hat \xi ^{(J)}(\rho)=\rho^{a}
		z^{b}\cdot  (\lambda_2 + F_\rho(z,\rho))   \hat \xi ^{(J)}(\theta),
	\end{equation}
	and since $\ell\geq\ell_\mathcal{M}+1$, Corollary~\ref{cor:elem} implies that:
	 \begin{equation*}
		 \xi _\ell ^{(J)}(\rho)
	= \rho^{a}
		z^{b}\cdot  (\lambda_2 +F_\rho^\ell(\theta,z,\rho))     \xi_\ell ^{(J)}(\theta),
	\end{equation*}
	where $F^\ell_\rho$ is analytic and $F^\ell_\rho(\theta,0,0)=0$.
	Considering that the monomial $\rho^{a}
	z^{b}>0$ for $(z,\rho)\in \R_{>0}^ 2$, and taking into account that $\lambda_2\neq 0$ and $\xi_\ell ^{(J)}(\theta)>0$ along $\gamma_I$, there is a neighborhood $W_I$ of $\g_I$ such that $\xi_\ell ^{(J)}(z)$ has constant sign in $W_I\setminus E$. Hence, the trajectories of $\xi_\ell ^{(J)}$ can be parameterized by $\rho$ in $W_I\setminus E$ and thus $\xi_{\ell}^{(J)}$ cannot have cycles in $W_I\setminus E$.
	
	In the second case, being both components of the divisor invariant, we can write:
		\begin{equation*}
		\label{eq:formal2dimgeneralreducednondicritical}
		\hat \eta_J'= (	\lambda_1 z +z F_z(z,\rho))\z + (\lambda_2\rho +  \rho F_\rho(z,\rho))\rhoo,
\end{equation*}
	where $\lambda_1^ 2+ \lambda_2^2 \neq 0$ and $F_\rho,F_z\in \R[\rho][[z]]\cap \R[z][[\rho]]$ satisfy $F_\rho(0,0)=F_z(0,0)= 0$. 
Suppose without loss of generality that $\lambda_1\neq 0$. Then, we write:
\begin{equation*}
	\hat \xi ^{(J)}(z)
=\rho^{a}
	z^{b+1}\cdot  (\lambda_1 + F_z(z,\rho))   \hat \xi ^{(J)}(\theta),
\end{equation*}
and since $\ell\geq \ell_\mathcal{M}+1$, Corollary~\ref{cor:elem} implies that:
\begin{equation*}
	\xi_\ell ^{(J)}(z)= \rho^{a}
	z^{b+1}\cdot  (\lambda_1 +F_z^\ell(\theta,z,\rho))     \xi_\ell ^{(J)}(\theta),
\end{equation*}
where $F^\ell_z$ is analytic and $F^\ell_z(\theta,0,0)=0$.
As in the first case, we find that the trajectories of $\xi_\ell ^{(J)}$ can be parameterized by $z$ in $W_I\setminus E$ and $\xi_\ell ^{(J)}$ cannot have cycles in $W_I\setminus E$. 
\end{proof}

\subsection{Simple non-corner characteristic cycles} \label{sec:noncorner}
All along this subsection, we suppose that $\gamma_I$ is a non-corner characteristic cycle of $\mathcal{M}$ contained in a chart $C_J$ for which $\{ \rho^{(J)}=0 \}$ is the equation of $E\cap C_J$ and $\gamma_I=\{ \rho^{(J)}=0, z^{(J)}=w_{I}   \}\subset\{ \rho^{(J)}=0 \} $ for some $\omega_I\in \R$. 
Consider the transform $\hat \xi^{(J)}=(\pi\lvert_{C_J})^\ast \hat \xi$ in the translated coordinates $(z:=z^{(J)}-w_I,\rho:=\rho^{(J)})$. Its associated two-dimensional vector field is
$$\hat \eta_J:= \frac{\hat \xi^{(J)}(\rho)}{\hat \xi^{(J)}(\theta)}\rhoo +  \frac{\hat \xi^{(J)}(z)}{\hat \xi^{(J)}(\theta)}\z. $$
More precisely, we write $\hat \eta_J=\rho^a \hat \eta_J'$ where $a\geq 0$ and $\hat \eta'_J$ is a formal vector field in coordinates $(z,\rho)$ with a simple singularity at the origin, One of the separatrices of $\hat \eta'_J$ is the divisor $\{\rho=0\}$ and the other is smooth and transverse to the divisor, denoted by $\hat \Gamma_I$.

\subsubsection{Invariant formal surface along $\gamma_I$}
Being $\hat \Gamma_I$ a formal non-singular curve transverse to $\{\rho=0\}$, it can be expressed as a formal graph $z=\hat h_I(\rho)$, where $\hat h_I(\rho)\in \R[[\rho]]$. Since $L_{\thetaa} \hat{\xi}_I=0$, we have that $\hat S_I:= \hat \Gamma_I\times \mathbb{S}^1$ is a formal invariant non-singular surface of $\hat{\xi}_I$ supported along the cycle $\gamma_I$. Its vanishing ideal $\text{id}(\hat S_I)$ in the ring $\R[\cos\theta,\sin \theta][[\rho,z]]$ is generated by $H_I(z,\rho):=z-\hat h_I(\rho)$. 
 Using this surface, we can also construct an invariant surface for the transformed vector field $\xi_{\ell}^{(J)}=(\pi\lvert_{C_J})^\ast \xi_\ell$ along the characteristic cycle $\gamma_I$, when $\ell$ is sufficiently large. More precisely,

\begin{prop}\label{prop:non-corner}
	Suppose that $\ell\geq \ell_\mathcal{M}+1$. Then, there is a formal invariant surface $\hat S_{\ell,I}$ of $\xi_{\ell}^{(J)}$ along $\gamma_I$ expressed in coordinates $(\theta, z,\rho)$ as the ideal in $\R[\cos \theta, \sin \theta][[z,\rho]]$ generated by some series of the form $H_{\ell,I}(\theta, z,\rho):=z-h_{\ell,I}(\theta,\rho),$
	where $h_{\ell,I}\in \R[\cos \theta,\sin \theta][[\rho]]$ with $h_{\ell,I}(\theta,0)=0.$
\end{prop}

\begin{proof}

	Recall from Proposition~\ref{prop:lifting} that there is a $\R-$algebra formal automorphism $\psi_\ell:\R[[\rho,\theta,z]]^3\longrightarrow \R[[\rho,\theta,z]]^3$ defined by:
	\begin{equation*}
		(\theta,z,\rho)\circ \psi_{\ell}=(\psi_{\ell}^\theta,\psi_{\ell}^z,\psi_{\ell}^\rho)=(\theta + O(\rho^2),z+ O(\rho^2),\rho + O(\rho^2)),
	\end{equation*}
	conjugating $\hat{\xi}^{(J)}$ to $\xi_{\ell}^{(J)}$.
	We consider the formal surface $S_{\ell,I}$ whose defining ideal is $\text{id}(S_{\ell,I})=(\widetilde H_{\ell,I}(\theta, z,\rho))$ where $\widetilde H_{\ell,I}=\psi_\ell^\ast (\hat H_\ell)$:
	$$ \widetilde H_{\ell,I}(\theta, z,\rho):= H_I\circ \psi_{\ell}  =\psi_{\ell}^z-\hat h_{I}(\psi_{\ell}^\rho).  $$
	Using $\frac{\partial \widetilde H_{\ell,I}}{\partial z}\neq 0$ and 
	applying the implicit function theorem to $\widetilde H_{\ell,I}$, we find an expression of the form $H_{\ell,I}=z-h_{\ell,I}(\theta,\rho)$ for a generator of $id(S_{\ell,I})$.
\end{proof}

\subsubsection{Poincaré first-return map associated to $\gamma_I$}\label{sec:Poincare}
By Remark~\ref{rmk:characteristiccyclesarecycles}, $\gamma_I$ is a trajectory of the vector field $\xi_\ell^{(J)}=(\pi\lvert_{C_J})^\ast \xi_\ell$ for $\ell\geq \ell_{\mathcal{ M}}+1$. 
Let $P=P_{\ell,I}:\Delta \to \{ \theta=0  \}$ be the Poincaré first-return map of $\xi^{(J)}_\ell$ relatively to $\gamma_I$, where $\Delta$ is a sufficiently small neighborhood of $(z,\rho)=(0,0)$ in $\{  \theta=0 \}$ in which $P$ is analytic.

Notice that the Poincaré map does not depend on the parameterization of the trajectories of the vector field, and hence, we can define it using any equivalent vector field. In particular, we are going to consider the vector field $\tilde \xi_\ell^{(J)}$ equivalent to $\xi_\ell^{(J)}$ obtained by the multiplication by the inverse of $\xi_\ell^{(J)}(\theta)$. That is:
\begin{equation}
	\label{eq:xitilde}
	\tilde \xi_\ell^{(J)}=\thetaa + \chi,\quad 
	\chi=\frac{\xi_\ell^{(J)}(z)}{\xi_\ell^{(J)}(\theta)} \z+ \frac{\xi_\ell^{(J)}(\rho)}{\xi_\ell^{(J)}(\theta)}\rhoo. 
\end{equation}
Notice that the components of $\chi$ are the right members of the system of ODEs $\eta_{\ell, J}$
introduced in section~\ref{sec:approximationandjets}. They belong to the ring of formal power series $\R [\cos \theta, \sin \theta,z][[\rho]]$ and they  converge uniformly on $\SSS^1\times [-\delta,\delta]^2$ for some $\delta>0$. We consider then $\wt{\xi}_\ell^{(J)}$ as an analytic vector field on $\mathbb{R}\times (-\delta,\delta)^2$ that is $2\pi-$periodic in the variable $\theta$. Notice, moreover, from Corollary~\ref{cor:elem}, that $\rho$ divides $\chi$ and hence $\wt{\xi}^{(J)}_\ell\lvert_{E}=\thetaa$.

 Denote by $\Phi^t:=\Phi^t_{\tilde \xi_{\ell}^{(J)}}$ the flow map of $\tilde \xi_{\ell}^{(J)}$. It is defined and analytic for $(t,(\theta,z,\rho))\in (-\varepsilon,2\pi +\varepsilon)\times\left( (-\varepsilon,2\pi +\varepsilon) \times V\right)$ where $V$ is a neighborhood of $0\in \Rdos$. Using that $\wt{\xi}_\ell^{(J)}(\theta)=1$, we obtain
 \begin{equation}\label{eq:flow}
 	\Phi^t(\theta,z,\rho)= (\theta+t, \Psi^t_z(\theta,z,\rho),\Psi^t_\rho(\theta,z,\rho)),
 \end{equation}
that is, the angle $\theta$ is the natural time for $\tilde \xi_{\ell}^{(J)}$. By definition, the Poincaré map is given by
\begin{equation}\label{eq:Poincare}
	P(z,\rho)=(\Psi_z^{2\pi}(0,z,\rho),\Psi_\rho^{2\pi}(0,z,\rho)).
\end{equation}
We are going to express the flow via the exponential map. To be precise, given any $G\in \Rcossin[[z,\rho]]$, we define:
\begin{equation*}
	\text{Exp}(t\tilde \xi_{\ell}^{(J)})(G):=\displaystyle \sum_{i=0}^{\infty}   \frac{t^i}{i!} (\tilde \xi_{\ell}^{(J)})^{(i)} (G),
\end{equation*}
where, for any vector field $\zeta$, $\zeta^{(0)}(G)=G$ and $\zeta^{(i)}(G)=\zeta(\zeta^{(i-1)}(G))$, if $i\geq 1$. Taking into account the above properties of the coefficients of $\tilde \xi_{\ell}^{(J)}$ in the basis $\thetaa,\z,\rhoo$, it is immediate to check that $ \text{Exp}(t\tilde \xi_{\ell}^{(J)})(G)\in \Rcossin [[t,z,\rho]] $. In the following result, we get some useful properties of this exponential map and its relation with the flow map.
Notice first that, if $G\in\Rcossin [[z,\rho]] $, then the composition $G\circ \Phi^t$, due to the analyticity of $\Phi^t$, has a formal Taylor expansion at $t=0$, denoted by $T_0(G\circ \Phi^t)$, as a formal power series in variables $(t,z,\rho)$, with analytic functions (trigonometric polynomials) of $\theta\in (-\varepsilon, 2\pi+\varepsilon)$ as coefficients.
\begin{proposition}\label{prop:expandflow}
	Let $G\in\Rcossin [[z,\rho]] $. We have:
	\begin{enumerate}
		\item $T_0(G\circ \Phi^t)=\text{Exp}(t\tilde \xi_{\ell}^{(J)})(G) \in \Rcossin [[t,z,\rho]]$
		\item For any $t_0\in [0,2\pi]$, the expression $\text{Exp}(t_0\tilde \xi_{\ell}^{(J)})(G)= \displaystyle \sum_{i=0}^{\infty}   \frac{t^i_0}{i!} (\tilde \xi_{\ell}^{(J)})^{(i)} (G)$ has a sense as a series in $\Rcossin [[z,\rho]]$ and we have
		\begin{equation}\label{eq:GcircFlow}
			G\circ \Phi^{t_0}=\text{Exp}(t_0\tilde \xi_{\ell}^{(J)})(G)
		\end{equation}
	
	\end{enumerate}
\end{proposition}
\begin{proof}
	We prove \textit{(1)} with the same arguments as the case in Loray's text  for holomorphic vector fields \cite[pag. 15]{loray-pseudo}: expand $G\circ \Phi^t$ as a Taylor series in $t$ at $t=0$, so that we get
	$$ T_0(G\circ \Phi^t)=\displaystyle \sum_{i=0}^\infty \frac{t^i}{i!} \left.\frac{\partial ^i(G\circ \Phi^t)}{\partial t^i}\right\lvert_{t=0}, $$
	and check that, for any $i\geq 1$, $\frac{\partial ^i(G\circ \Phi^t)}{\partial t^i}= (\tilde \xi_{\ell}^{(J)})^{(i)} (G)\circ \Phi^t$.
	
	Let us prove item \textit{(2)}. First, we show that there exists $\alpha>0$ such that \textit{(2)} is true for any $t_0\in [0,\alpha]$. For that, consider the particular case where $G$ is either the coordinate $z$ or $\rho$ (with the notations of \eqref{eq:Poincare}, $z\circ \Phi^{t_0}=\Psi_z^{t_0}$ and $\rho\circ \Phi^{t_0}=\Psi_\rho^{t_0}$). By analyticity of these functions and by item \textit{(1)}, there exists some $\alpha$ such that $\text{Exp}(t\tilde \xi_{\ell}^{(J)})(z), \text{Exp}(t\tilde \xi_{\ell}^{(J)})(\rho)\in \Rcossin\{ t \}_\alpha [[ z,\rho  ]]$, where $\Rcossin\{ t \}_\alpha$ denotes the subring of $\Rcossin [[t]]$ consisting in series whose sequence of partial sums converges to an analytic function on the domain $(t,\theta)\in (-\beta,\beta)\times \mathbb{R}$ for some $\beta>\alpha$.
	Moreover, we have $\Psi_z^{t_0}=z\circ \Phi^{t_0}= \text{Exp}(t_0\tilde \xi_{\ell}^{(J)})(z)$ and  $\Psi_\rho^{t_0}=\rho\circ \Phi^{t_0}= \text{Exp}(t_0\tilde \xi_{\ell}^{(J)})(\rho)$ for any $t_0\in [0,\alpha]$.

	Let $G\in \Rcossin [[z,\rho]]$ be any formal series and write
	$$ G=\displaystyle \sum_{u,v} G_{uv}(\theta)z^u\rho^v , \text{ with } G_{uv}(\theta)\in \Rcossin.$$
	Consider the series
	$$\bar G= \displaystyle \sum_{u,v} G_{uv}(\theta+t)z^u\rho^v  $$
	which belongs to $\Rcossin\{ t \}_\alpha [[ z,\rho  ]]$ since each $G_{uv}(\theta)$ is a trigonometric polynomial. Taking into account the expression of the flow $\Phi^t$, we have that $G\circ \Phi^t$ is the result of substituting in the series $\bar{G}$ the variables $z,\rho$ by $\Psi^t_z,\Psi^t_\rho$, respectively. Since the series $\Psi_z^t,\Psi_\rho^t$ belong to $\Rcossin\{ t \}_\delta [[ z,\rho  ]]$ and have positive order with respect the variables $z,\rho$, substitution has perfect sense and provides an element in $\Rcossin\{ t \}_\alpha [[ z,\rho  ]]$.
	Since, by item \textit{(1)}, $G\circ \Phi^t$ also coincides with $\text{Exp}(t\widetilde{\xi}_\ell^{(J)})(G)$ as a series in $\Rcossin [[t,z,p]]$, we conclude item \textit{(2)} and~\eqref{eq:GcircFlow} for $t_0\in [0,\alpha]$. 
	Notice that $\alpha>0$ does not depend on $G$. Let us show that we can extend the property\eqref{eq:GcircFlow} to any $t_0\in [0,2\alpha]$ (and hence similar extensions will prove \textit{(2)}). Let $t_0\in [\alpha,2\alpha]$ and write $t_0=s_0+\alpha$, where $s_0\in [0,\delta]$. We have $G\circ \Phi^{t_0}=(G\circ \Phi^{s_0})\circ \Phi^\alpha$.  Applying~\eqref{eq:GcircFlow} for the values $s_0$ and $\alpha$, and for $G$ and $G\circ \Phi^{s_0}$, respectively, we get
	\begin{equation*}
		\begin{split}
			 	G\circ \Phi^{t_0}&=\displaystyle \sum_i \frac{\alpha^i}{i!} (\tilde \xi_{\ell}^{(J)})^{(i)} (G\circ \Phi^{s_0})= \sum_i \frac{\alpha^i}{i!} (\tilde \xi_{\ell}^{(J)})^{(i)} \left( \sum_j \frac{\alpha^j}{j!} (\tilde \xi_{\ell}^{(J)})^{(j)} (G) \right)  \\
			 	&= \sum_k \left( \sum_{i+j=k}  \frac{\alpha^i}{i!} \frac{s_0^j}{j!}(\tilde \xi_{\ell}^{(J)})^{(k)}(G)     \right)= \sum_k \frac{(\alpha+s_0)^k}{k!}(\tilde \xi_{\ell}^{(J)})^{(k)}(G) = \text{Exp}(t_0 \tilde \xi_{\ell}^{(J)})(G), 
		\end{split}
	\end{equation*}
as it was to be proved.
\end{proof}

We can now prove two important features of the Poincaré map $P$ of $\tilde \xi_{\ell}^{(J)}$ along $\gamma_I$. 
\begin{lem}\label{lem:PnoId}\label{lem:Pinvariantcurve}
	There exists $\ell_I$ such that for any $\ell\geq \ell_I$, the Poincaré map $P$  satisfies:
\begin{enumerate}
	\item[a)] $P$ is tangent to the identity but $P\neq Id$ as a germ of diffeomorphisms at $(0,0)\in \Delta$.
	\item[b)] The formal curve $\Gamma=\Gamma_{\ell,I}:=\hat{S}_{\ell,I}\cap \Delta$ is invariant for $P$.
\end{enumerate}
\end{lem}
\begin{proof}
	Recall that the two-dimensional formal vector field $\hat \eta_J$ associated to the formal vector field $\hat \xi^{(J)}$ has an adapted simple singularity corresponding to the characteristic cycle $\g_I$, and suppose that the the formal invariant curve $\Gamma_I$ is given by the ideal generated by $z-\hat h_I(\rho)$ where $\hat h_I(\rho)=\sum_{i\geq 1}a_i\rho^i$.
	Therefore, we write
	\begin{equation*}\label{eq:formalsystemgeneralsimple}
		\hat \eta_J= \rho^{n} (z-\hat h_{I})^r \left(   (\lambda_1 (z-a_1\rho)+ B_1(z,\rho)) \z + (\lambda_2\rho+ B_2(z,\rho))\rhoo  \right),
	\end{equation*}
	where $n=n_1^{(J)}$, $r\in \mathbb{N}_{\geq 0}$, $(\lambda_1,\lambda_2)\neq (0,0)$ and $B_i\in \R[z][[\rho]]$ has order in $\rho$ greater or equal to 2 for $i=1,2$.
	
	Define $\ell_I=\ell_{\mathcal{ M}}+ r+1$. Applying Corollary~\ref{cor:elem} to $k=r+1$, we get, for any $\ell\geq \ell_{\mathcal{ M}}+k=\ell_I$: $$j_{r+1}^\rho (\eta_{\ell, J}')=j_{r+1}^\rho (\hat \eta_{ J}')$$
	By equation~\eqref{eq:xitilde} and the equality between the system $\eta_{\ell, J}$ and $\chi$ stated in the subsequent paragraph, we obtain an expression of the vector field $\tilde \xi_\ell^{(J)}$
	\begin{equation}
		\label{eq:xitilde2}
		\tilde \xi_\ell^{(J)}=\thetaa + \rho^n j_{r+1}^\rho\left((z-\hat{h}_I(\rho))^r \left( (\lambda_1 (z-a_1\rho) + B_1(\theta,z,\rho))\z+   (\lambda_2 \rho+ B_2(\theta,z,\rho)) \rhoo  \right)\right) +O(\rho^{n+r+2})  , 
	\end{equation}
	where $O(\rho^{n+r+2})$ denotes a vector field whose coefficients can be divided by $\rho^{n+r+2}$.
	
	Assume for instance that $\lambda_1\neq 0$. We compute the exponential of the vector field  $\tilde \xi_{\ell}^{(J)}$ applied to the coordinate function $z$:
	\begin{equation*}
		\begin{split}
			\text{Exp}(t \tilde \xi_{\ell}^{(J)})(z)& = z+ t\rho^n\lambda_1 \left( \sum_{s=0}^{r} \binom{r}{s}  z^s (-a_1)^{r-s}\rho^{r-s}   \right)\cdot (z-a_1\rho)+O({n+r+2})+\\
			&+ \frac{t^2}{2} \xi_{\ell}^{(J)} \left(  \rho^n\lambda_1 \left( \sum_{s=0}^{r} \binom{r}{s}  z^s (-a_1)^{r-s}\rho^{r-s}   \right)\cdot (z-a_1\rho)+O({n+r+2})  \right) +\cdots,
		\end{split}
	\end{equation*}
which implies that $$j_{n+r+1}(\text{Exp}(t \tilde \xi_{\ell}^{(J)})(z))=z+ t\rho^n\lambda_1 \left( \sum_{s=0}^{r} \binom{r}{s}  z^s (-a_1)^{r-s}\rho^{r-s}   \right)\cdot (z-a_1\rho). $$ Using Proposition~\ref{prop:expandflow} and equation~\eqref{eq:Poincare}, we find
$$z\circ P(z,\rho)=(z\circ \Psi^{2\pi})(0,z,\rho) = z+ 2\pi\rho^n\lambda_1 \left( \sum_{s=0}^{r} \binom{r}{s}  z^s (-a_1)^{r-s}\rho^{r-s}   \right)\cdot (z-a_1\rho)+ \text{h.o.t}.$$
Since $\lambda_1\neq 0$, we conclude that $z\circ P\neq z$ and thus, $P\neq Id$. Analogous arguments lead to the same conclusion if $\lambda_2\neq 0$. Since, by the above construction, $j_1^{(z,\rho)}(z\circ \Psi^{2\pi}, \rho\circ \Psi^{2\pi})=(z,\rho)$, we have that $P$ is tangent to the identity and finish the proof of a).

Let us show b). Let $H=H_{\ell,I}\in \Rcossin [[z,\rho]]$ be a generator of the ideal of the invariant surface $S_{\ell,I}$ obtained in Proposition~\ref{prop:non-corner}. The series $g(z,\rho):=H(0,z,\rho)$ is a generator of the formal plane curve $\Gamma=\hat{S}_{\ell,I}\cap \{  \theta=0 \}$. We need to check that the composition $g\circ P$ is divisible by $g$.
Using Proposition~\ref{prop:expandflow}, we have
\begin{equation}\label{eq:invariantcurve}
	H\circ \Phi^{2\pi} = \displaystyle \sum_{i\geq 0} \frac{(2\pi)^i}{i!} (\tilde \xi_{\ell}^{(J)}) ^{(i)}(H).
\end{equation}
Since $\hat{S}_{\ell,I}$ is invariant for $\tilde \xi_{\ell}^{(J)},$ we have $\tilde \xi_{\ell}^{(J)}(H)\in \text{id}(\hat{S}_{\ell,I})$, that is, $H$ divides $\tilde \xi_{\ell}^{(J)}(H)$. By recurrence, $H$ divides $(\tilde \xi_{\ell}^{(J)}) ^{(i)}(H)$ for any $i\geq 0$. Thus, from equation~\eqref{eq:invariantcurve}, we get
$$ H\circ \Phi^{2\pi}=H\cdot K , \ \text{where }K\in \Rcossin[[z,\rho]] .$$
We conclude, using that $P(z,\rho)=\Phi^{2\pi}(0,z,\rho)$,
$$ g\circ P=(H\circ \Phi ^{2\pi})\lvert_{  \{  \theta=0 \} }= H(0,z,\rho)  K(0,z,\rho) = g \cdot \tilde K,\ \tilde K\in \R[[z,\rho]], $$
as we wanted to prove.
\end{proof}

\subsubsection{Periodic orbits of the Poincaré map around the invariant curve}
 We are interested in the periodic orbits of $P$ near $(0,0)$ since, as we know, they correspond to cycles of $\xi_\ell^{(J)}$ near $\gamma_I$. As it was the case for defining cycles, periodic orbits depend on the domain of (a representative of) $P$. To be more  precise, if $P$ is defined in some neighborhood $W$ of $(0,0)$, a \emph{periodic orbit in }$W$ is a finite set $\{p_i\}_{i=0}^{n-1}$  contained in $ W$ such that $p_i=P(p_{i-1})$ for $i=1,\ldots, n-1$ and $p_0=P(p_{n-1})$.

Denote by $\fix (P)$ the (germ of the) locus of fixed points of $P$.
 We can distinguish two different situations:
\begin{enumerate}
	\item[a)] The invariant curve $\Gamma$ is not contained in $\fix (P)$.
	\item[b)] The invariant curve $\Gamma$ is contained in $\fix (P)$. (In particular, in this case, $\Gamma$ is convergent since $\fix(P)$ the locus of fixed points is an analytic set of dimension at most one.)
\end{enumerate}
In both cases, we investigate periodic orbits of $P$ in some  conic neighborhood of $\Gamma$. To be precise, consider a parametrization of $\Gamma$ of the form $ z=h(\rho)$, with $h(\rho)\in \R[[\rho]]$ and $h(0)=0$, For any $N\in \mathbb{N}_{\geq 1}$ and any $\delta>0$ sufficiently small, we define the set:
\begin{equation*}\label{eq:defconos}
	\Sigma_{N,\delta} (\Gamma)= \Sigma_{ N,\delta}^{(z,\rho)}(\Gamma):=\{  (z,\rho): |z-j_N(h(\rho))|<\rho^N, 0<\rho<\delta  \}.
\end{equation*}

\begin{remark}\label{rmk:coneschangeofcoord}
	By its definition, these conic neighborhoods depend on the chosen coordinates. However, after a simple change of variables consisting in a transformation of type $\bar z=z +\alpha(\rho) $ with $\alpha(\rho)\in \R[[\rho]]$, for any $N,\delta$ the sets $\Sigma^{(z,\rho)}_{N,\delta}$ and $\Sigma^{(\bar z,\bar \rho)}_{N,\delta}$ are identical under this change.
\end{remark}
\paragraph{\textbf{Case a)}}
In this case, we prove that there are no periodic orbits in such a conic neighborhood. The arguments are inspired by the papers \cite{lopez-sanz2018,lop-rai-rib-sanz-2019}, devoted to treat this case for holomorphic diffeomorphisms. 
\begin{lem}\label{lem:curvenotfixed}
	Suppose that $\Gamma$ is not contained in $\fix(P)$. Then, there is some $N\in \mathbb{N}_{\geq 1}$ and some $\delta>0$ such that (a representative of) $P$ does not have periodic points in $\Sigma_{ N,\delta}(\Gamma)$.
\end{lem}
\begin{proof}
	First, since the divisor is contained in $\fix(P),$ we have that $\rho\circ P -\rho$ can be divided by $\rho$. On the other hand, being $\Gamma$ invariant for $P$, there is a formal diffeomorphism $\Theta(\rho)=\rho+ O(\rho^2)\in \R[[\rho]]$ satisfying:
	\begin{equation*}
		P(h(\rho), \rho)=(h(\Theta(\rho)), \Theta(\rho))
	\end{equation*}
	The formal diffeomorphism $\Theta$ is called the \emph{restriction of $P$ to $\Gamma$}, denoted by $P\lvert_{\Gamma}:=\Theta$, see \cite{lopez-sanz2018}. The \emph{order} of $P\lvert_\Gamma$, defined as $\ord_\rho(\Theta(\rho)-\rho)-1$, does not depend on the coordinates nor the parametrization $(h(\rho),\rho)$ of $\Gamma$. In this case, it is a natural number $m<\infty$, because otherwise, $\Gamma$ would be contained in $\fix(P)$. We deduce that there is a maximal $k\geq 1$ such that $\rho^k$ divides $\rho\circ P-\rho$, so that we can write
	\begin{equation}\label{eq:rhocircP}
		\rho\circ P-\rho=\rho^k(A(\rho)+ zB(z,\rho)),
	\end{equation}
	where $A\in \mathbb R \{ \rho\}$, $B\in \mathbb R \{z, \rho\}$ and $\rho^k(A(\rho)+ zB(z,\rho))\neq 0$. Up to taking new coordinates $(\bar z, \rho)$ with $\bar z =z-j_{m+1}(h(\rho))$ (which we rename $(z,\rho)$ for simplicity) and using Remark~\ref{rmk:coneschangeofcoord}, we may assume from the begining that $\ord _\rho (h(\rho))\geq m+2$. From equation \eqref{eq:rhocircP}, the definition of $\Theta=P\lvert_\Gamma$ and its order $m$, we have
	\begin{equation*}
		\alpha \rho^{m+1} + \cdots= \Theta(\rho)-\rho=\rho^k (A(\rho)+  h(\rho)B(h(\rho),\rho) ), \ \alpha\neq 0,
	\end{equation*}
	which implies that $A(\rho)\neq 0$ writes as $A(\rho)=\rho^s \wt{A}(\rho)$ with $k+s=m+1$ and $\widetilde A(0)=\alpha$ (in particular $s$ is finite). Put $N=m+1$ and let us prove the required property for a cone $\Sigma_{ N,\delta}=\Sigma_{ N,\delta}^{(z,\rho)}(\Gamma)$ with some $\delta>0$. 
	 Notice first that for the chosen coordinates $(z,\rho)$, we have $j_N(h(\rho))=0$, so $\Sigma_{ N,\delta}$ is given simply by equations $|z|<\rho^N$ and $0<\rho<\delta$. On the other hand, $N=k+s>s$, since $k>0$. Assume for instance that $\alpha<0$ (analogous arguments apply if $\alpha>0$). Take a preliminary $\delta_1>0$ such that $\wt{A}(\rho)<\frac{\alpha}{2}$ if $\rho<\delta_1$ and let $K>0$ be a bound for $|B|$ in a neighborhood of $(0,0)$ that contains $\Sigma_{N,\delta_1}$. We have in $\Sigma_{N,\delta_1}$
	\begin{equation}\label{eq:boundsinSigmanofix}
		\rho^s \wt{A}(\rho)+ zB<\rho^s\wt{A}(\rho)+\rho^N K <\rho^s(\frac{\alpha}{2}+\rho^{N-s}K).
	\end{equation}
	Now take $\delta<\delta_1$ such that $\delta^{N-s}<\frac{|\alpha|}{4K}$, so that, from \eqref{eq:boundsinSigmanofix} we obtain
	\begin{equation}\label{eq:endlemmanofix}
		\rho\circ P-\rho=\rho^k (\rho^s \wt{A}(\rho)+ zB)<\rho^k (\rho^s\frac{\alpha}{4})<0 \ \text{in }\Sigma_{ N,\delta}.
	\end{equation}
	We conclude that, if $\{ p_0, p_1=P(p_0), p_2=P(p_1), \ldots   \}$ is an orbit (finite or not) of $P$ completely contained in $\Sigma_{ N,\delta}$, then, the sequence of $\rho-$coordinates decreases strictly:
	$$ \rho(p_0)>\rho (p_1)>\cdots    $$
	This proves the result.
\end{proof}
\paragraph{\textbf{Case b)}}
Suppose that $\Gamma$ is convergent and contained in $\fix(P)$. Convergence means that $\Gamma$ is an analytic curve given by a graph $\Gamma= \{ (h(\rho),\rho): \rho\in [0,\varepsilon)  \}$, where $h(\rho)\in \R \{ \rho\}$. Up to taking new analytic coordinates $(z-h(\rho),\rho)$, we will assume that $h(\rho)\equiv0$. Thus, the sets $\Sigma_{ N,\delta}=\Sigma_{ N,\delta}^{(z,\rho)}$ will be simply defined by the equations $|z|<\rho^N$ and $0<\rho<\delta$. With these assumptions, we prove the following result.
\begin{lemma}\label{lem:curveoffixed}
	Suppose that $\Gamma= \{ z=0   \}\subset \fix(P)$. Then, there is $N\in \mathbb N_{\geq 1}$ and $\delta>0$ such that the only periodic orbits of $P$ in $\Sigma_{ N,\delta}$ are the fixed points of the set $\Gamma\cap \Sigma_{ N,\delta}$.
\end{lemma}
\begin{proof}
	The two coordinate axis $\{ \rho=0\}$ and $\{  z= 0\}$ are contained in $\fix (P)$. Thus, both components $(z\circ P- z,\rho\circ P -\rho)$ of the map $P-id$ are divisible by a positive power of $z$ and by a positive power of $\rho$. In particular, we can write $z\circ P= z(1+ \psi (z,\rho))$ where $\rho$ divides $\psi$. From this, we prove the following observation.\\
	\noindent \textbf{Claim:} There is a neighborhood $V$ of $(0,0)$ such that if $\{ p_0,p_1,\ldots  \}$ is an orbit of $P$ contained in $V$, then the sign of the $z$ coordiante of its elements is constant.\\
	\textit{Proof of the claim:}
		Since $\psi(0,0)=0$, we can consider a neighborhood $V$ where we have 
		$(1+\psi(z,\rho))>\frac{1}{2}$.
		Hence, $  \text{Sign}(z\circ P)=\text{Sign}(z(1+\psi(z,\rho)))=\text{Sign}(z)  $ and the claim follows.	
	
	On the other hand, since $P\neq id$, the two components $z\circ P-z, \rho \circ P-\rho$ cannot be identically zero simultaneously.
	Suppose that $\rho\circ P-\rho\neq 0$. Then we can write
	\begin{equation}
		\label{eq:rhonofixed}
		\rho \circ P -\rho= \rho^{k_1}z^{k_2}(A(\rho)+ zB(z,\rho)),
	\end{equation}
	where $k_1,k_2\in \mathbb{N}_{\geq 1}$ and $A(\rho)$ is a convergent non-zero series. We write $A(\rho)=\rho^s (\alpha+ \cdots),$ where $s\geq 0$ and $\alpha\neq 0$. Analogously as in the proof of Lemma~\ref{lem:curvenotfixed}, if $N$ is any given natural number with $N>s$, then there exists $\delta>0$ such that the function $A(\rho)+ zB(z,\rho)$ has constant non-zero sign on $\Sigma_{ N,\delta}$. Taking into account the claim above, if $\delta$ is sufficiently small so that $\Sigma_{ N,\delta}\subset V$, we conclude form equation~\eqref{eq:rhonofixed} that if $\{  p_0,p_1,\ldots \}$ is an orbit of $P$ contained in $\Sigma_{ N,\delta}\setminus \Gamma= \left(\Sigma_{ N,\delta}\cap \{ z>0\}   \right)  \cup \left(\Sigma_{ N,\delta}\cap \{ z<0\}   \right)  $, then, the sequence $\{ \rho(p_0), \rho(p_1), \ldots  \}$ is strictly increasing or strictly decreasing. This proves the lemma in case $\rho\circ P -\rho \neq 0$.
	On the contrary, if $\rho\circ P -\rho = 0$ but $z\circ P - z\neq 0$, we obtain analogously that the $z-$coordinate of elements of an orbit in $\Sigma_{ N,\delta}\setminus \Gamma$ is strictly increasing or strictly decreasing, which proves the lemma.
\end{proof}

\section{Proof of the main theorem}
In this section we give the proof of Theorem~\ref{th:main}. It is enough to prove the result for some jet approximation $\xi_{\ell}$ of a formal normal form $\hat \xi$ of $\xi$, since all those vector fields are locally analytically conjugated to $\xi$ at $0\in \Rtres$. Fix an adapted reduction of singularities $\mathcal{ M}=(M,\pi,\mathcal{A}, \mathcal{D})$ of $\hat \xi$ given by Proposition~\ref{prop:adaptedressing} and denote by $\mathcal{D}^{nc}$ the subset of $\mathcal{D}$ consisting on non-corner characteristic cycles. For any $\gamma_I\in \mathcal{D}^{nc}$, we consider a chart $C_J$ and coordinates $(\theta, \bar{z}^{(I)}:=z^{(J)}-\omega_I, \rho^{(I)}=\rho^{(J)})$ as in section~\ref{sec:noncorner}, such that $\gamma_I\subset C_J$ is given $\g_I=\{ z^{(I)}=0, \rho^{(I)}=0 \}$. Let $\hat \Gamma_I$ be the formal plane curve at the origin $(z^{(I)},\rho^{(I)})=(0,0)$ of $\{  \theta=0 \}$, invariant for the associated vector field $\hat{\eta}_J$ and transversal to the divisor $\{  \rho^{(I)}=0 \}$. Consider the parameterization of $\hat \Gamma_I$  given in these coordinates as $z=\hat h_I(\rho^{(I)})$, $\hat h_I\in \R[[\rho^{(I)}]]$, $\hat h_I(0)=0$. Consider also the formal surface $\hat S_I=\SSS^1\times \hat \Gamma_I$ invariant for $\hat \xi ^{(J)}$ and given also by the same equation $z=\hat h_I(\rho^{(I)})$, but considering it in $\Rcossinz[[\rho]]$. Let $\ell\geq \ell_{\mathcal{ M}}+1$ be the first natural number such that $\ell\geq \ell_I$ for any $\gamma_I \in \mathcal{D}^{nc}$, where $\ell_I$ is given in Lemma~\ref{lem:PnoId}. Denote by $S_{\ell,I}$ the formal invariant surface of $\xi_\ell^{(J)}$ given in Proposition~\ref{prop:non-corner} with defining ideal $id(S_{\ell,I})=(H_{\ell,I})$, $H_{\ell,I}=z^{(I)}-h_{\ell,I}(\theta, \rho^{(I)})$ where $h_{\ell,I}\in \Rcossin[[\rho^{(I)}]]$. Denote by $P_{\ell,I}$ the Poincaré map of $\xi_\ell^{(J)}$ along $\g_I$ defined in some neighborhood of the origin of $\{ \theta=0 \}$ and $\G_{\ell,I}=S_{\ell,I}\cap \{  \theta=0 \}$ its corresponding formal invariant curve (by Lemma~\ref{lem:Pinvariantcurve}). Denote furthermore $\mathcal{D}^{nc}=\mathcal{D}^{nc,non-fix}\cup \mathcal{D}^{nc,fix}$ as a disjoint union, where $\g_I\in \mathcal{D}^{nc,fix}$ if and only if $\Gamma_{\ell,I}\subset \fix(P_{\ell,I})$.

Apply Lemma~\ref{lem:curvenotfixed} or~\ref{lem:curveoffixed} according to whether $\g_I\in\mathcal{D}^{nc,non-fix} $ or $\g_I\in \mathcal{D}^{nc,fix}$, and get conic neighborhoods $\Sigma_{N_{\ell,I},\delta}=\Sigma_{N_{\ell,I},\delta}(\Gamma_{\ell,I})$ (in coordinates $(z^{(I)},\rho^{(I)})$ and where we have unified $\delta>0$) where $P_{\ell,I}$ has no periodic orbits except for the fixed points $\Gamma_{\ell,I}\cap \Sigma_{N_{\ell,I},\delta}$ when $\g_I\in\mathcal{D}^{nc,fix} $.

This information is of course relevant in order to describe the union of cycles of the transform $\xi_{\ell}^{(J)}$ near $\g_I$, but it is not enough a priori, since $\Sigma_{N_{\ell,I},\delta}$ is not a full neighborhood of the origin at $\{  \theta=0 \}$. By further blowing-ups along the characteristic cycles $\g_I$, one can eventually open these conic neighborhoods to full neighborhoods of $\g_I$ but, if the order of tangency $N_{\ell,I}$ is too large, the sequence of blowing-ups to be done may not follow the formal curve $\hat \Gamma_I$, and thus we could skip the context of sequences of admissible blowing-ups. We can take a bigger $\ell'$ so that the formal curve $\Gamma_{\ell,I}$ approximates $\hat{\Gamma}_I$ better than the curve $\Gamma_{\ell,I}$ does. But the order $N_{\ell',I}$ may increase with $\ell'$ a priori.

In the strategy that follows we overcome these difficulties. We have to consider first the same kind of conic three-dimensional neighborhoods of the surfaces $S_{\ell,I}$.
With more generality, consider coordinates $(\theta,z,\rho)$ in $\SSS^1\times \R\times\R_{\geq 0}$, $\g=\SSS^1\times \{ (0,0) \}$ and let $S$ be a formal non-singular surface along $\g$ given by an equation of the form 
$$z-h(\theta,\rho)=0 , \ \text{where } h(\theta,\rho)\in \Rcossin[[\rho]]. $$
For $N\in \mathbb{N}_{\geq 0}$ and constants $C,\delta>0$, we define
$$\wt{\Sigma}_{N,\delta,C} (S)= \wt{\Sigma}_{N,\delta,C}^{(\theta,z,\rho)} (S)= \{ (\theta,z,\rho): \rvert  z-j_N^\rho(h(\theta,\rho))  \lvert < C \rho^N  , 0<\rho<\delta \}. $$
In particular, notice that $\wt{\Sigma}_{N,\delta,1} (S_{\ell,I})=\Sigma_{N,\delta}(\Gamma_{\ell,I})\cap \{  \theta=0 \}$. We need two results, the first one is just a remark that follows from the construction of sequences of admissible blowing-ups in section~\ref{sec:adaptedres} and the definition of adapted simple singularities.

\begin{remark}\label{rmk:blupsimplesing}
	If $\g_I\in \mathcal{D}^{nc}$ and $\widetilde {\mathcal{M}}=(\widetilde M,\widetilde\pi,\widetilde{\mathcal A},\widetilde {\mathcal D})$ is the sequence of admissible blowing-ups obtained from $\mathcal {M}$
	 by blowing-up along $\g_I$ (that is, $\widetilde\pi=\pi\circ \sigma_{\gamma_{I}}$ in the notations of section~\ref{sec:adaptedres}), then $\widetilde {\mathcal{M}}$ is again an adapted reduction of singularities of $\hat \xi$ with $\widetilde {\mathcal{D}}=\mathcal{D}\setminus \{ \g_I\}\cup \{ \gamma_{I,-\infty}, \g_{I,\infty}, \g_{I,1}  \}$ and $\widetilde {\mathcal{D}}^{nc}=(\mathcal{D}^{nc}\setminus \{\g_I\}) \cup \{  \gamma_{I,1} \}
		$. We will say that $\g_{I,1}$ \emph{emerges from} $\g_I$. Moreover, the new non-corner characteristic cycle $\g_{I,1}$ is given by equations $\{   z^{(I,0)}=\rho^{(I,0)}=0 \}$ in coordinates for which $\sigma_{\gamma_I}$ is written as
		$$  z^{(I)} = \rho^{(I,0)}(z^{(I,0)}+a_{I,1}), \ \rho^{(I)}= \rho^{(I,0)} ,
		  $$
		  where $j_1^{\rho^{(I)}}(\hat h_I(\rho))=a_{I,1}\rho$. We deduce that if $N\in \mathbb{N}$ and $j_N^{\rho^{(I)}}(\hat h_I)=j_N^{\rho^{(I)}}(h_{\ell,I}),$ then the strict transform $\sigma_{\gamma_{I}}^\ast S_{\ell,I}$ by $\sigma_{\gamma_{I}}$ is a formal surface along $\gamma_{I,1}$ and
		  $$ \sigma_{\gamma_{I}}^{-1}(\wt{\Sigma}_{N,\delta,1}^{(\theta,z^{(I)},\rho^{(I)})} (S_{\ell,I}))  = \wt{\Sigma}_{N-1,\delta,1}^{(\theta,z^{(I,1)},\rho^{(I,1)})} (S_{\ell,I}) . $$
\end{remark}

\begin{lem}\label{lem:cones}
	Let $\varphi(\theta,z,\rho)=(\theta+ \varphi_\theta, z+ \varphi_z,\rho+ \varphi_\rho )$ be a diffeomorphism along $\g=\{  z=0,\rho=0  \}$, where $\varphi_\theta,\varphi_z,\varphi_\rho\in \Rcossinz[[\rho]]$ are of order at least two in $(\rho,z)$ and divisible by $\rho$.
	Let $S_i$, $i=1,2$ be formal surfaces with defining ideals $id(S_i)=(z-f_i(\theta,\rho))$ with $f_i(\theta,\rho)\in \Rcossin[[\rho]]$ and such that $\varphi^\ast (S_2)=\varphi(S_1)$. Then, for every cone $\wt{\Sigma}_{N,C_1,\delta_1}(S_1)$ there exists some $C_2,\delta_2>0$ such that:
	$$\varphi^{-1}(\widetilde{\Sigma}_{N,C_2,\delta_2}(S_2)) \subseteq \widetilde{\Sigma}_{N,C_1,\delta_1}(S_1).$$
\end{lem}
\begin{figure}
\centering
\includegraphics[width=0.5\linewidth]{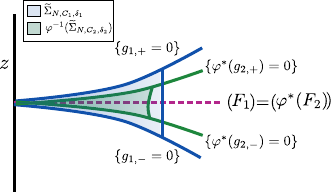}
\caption{$\wt{\Sigma}_{N,C_1,\delta_1}(S_1)$ and $\varphi^{-1}(\wt{\Sigma}_{N,C_2,\delta_2}(S_2))$}
\label{fig:cones}
\end{figure}
See figure~\ref{fig:cones} for an illustration of this lemma.
 \begin{proof}
 	Consider the functions $g_{i,\epsilon}=z-j_N^{\rho}(f_i(\theta,\rho))-\epsilon C_i\rho^N$ for $i=1,2$ and $\epsilon\in \{ -1,+1  \}$, so that the equations of the boundary surfaces of the cones $\widetilde{\Sigma}_{N,C_i,\delta_i}(S_i)$ are $g_{i,+}=0$ and $g_{i,-}=0$. It is enough to prove that there exists $C_2,\delta_2>0$ such that the function $\rho \circ \varphi >\delta_2$ in the points $D=\{ \rho=\delta_1 \}\cap \overline{ \wt{\Sigma}_{N,C_1,\delta_1}(S_1)  }$, the function
 	 $g_{2,+}\circ \varphi >0$ in the points that fulfill $\{ g_{1,+}=0, 0<\rho\circ \varphi  <\delta_2  \}$ and the function $g_{2,-}\circ \varphi <0$ in the points that fulfill $\{ g_{1,-}=0, 0<\rho \circ \varphi <\delta_2 \}$.

 	 A first upper bound to define $\delta_2$ follows from $\rho \circ \varphi (\theta,z,\delta_1)
 	 \geq K_1>0$ in the compact set $D$ where it is being evaluated. Hence, we must have $\delta_2<K_1$.
 	 
 	 Notice that we have $z-j_{N} ^\rho (f_i)= F_i + O(\rho^{N+1}), \ i=1,2,$ and $F_2\circ \varphi=F_1\cdot U $ where $U$ is a unit in $\Rcossinz[[\rho]] $ with $U=1+T,$ such that $\rho$ divides $T$.
	Then, we deduce
 $$g_{2,+}\circ \varphi =(F_2 + O(\rho^{N+1})-C_2\rho^N )\circ \varphi= F_1\cdot U + ( -C_2\rho^N+ O(\rho^{N+1}))\circ \varphi. $$
 	Using that $\varphi^\ast (\rho)= \rho + \varphi_\rho=\rho +\rho \tilde \varphi_\rho$, we have:
 	$$g_{2,+}\circ \varphi =F_1\cdot U -C_2\rho^N +O(\rho^{N+1}).$$
 	Now, we evaluate $g_{2,+}\circ \varphi$ in the points of $g_{1,+}=0$, in other words, when $z=j_{N}^\rho(f_1)+ C_1\rho^N$. Considering that $F_1(\theta,j_{N}^\rho(f_1)+ C_1\rho^N,\rho)=f_1+O(\rho^{N+1})+C_1\rho^N -f_1= C_1\rho^N+O(\rho^{N+1})$, we get:
 	$$  g_{2,+}\circ \varphi (\theta,j_{N}^\rho(f_1)+ C_1\rho^N,\rho)= (C_1\rho^N+O(\rho^{N+1})) \cdot (1+T(\theta,j_{N}^\rho(f_1)+ C_1\rho^N,\rho)) -C_2\rho^N + O(\rho^{N+1}),
 	$$
 	Since $\rho$ divides $T$, we obtain finally:
 	$$\varphi^\ast (g_{2,+})(\theta,j_{N}^\rho(f_1)+ C_1\rho^N,\rho)=(C_1-C_2)\rho^N +O(\rho^{N+1}). $$
 	Thus, taking $C_2<C_1$, we have that $ g_{2,+}\circ \varphi(\theta,j_{N}^\rho(f_1)+ C_1\rho^N,\rho)>0$ for $\rho$ small enough, as we required. Consider the upper bound of $\rho$ to be $\widetilde K_2$. Then, we obtain a second upper bound for $\delta_2$, that is, $K_2=\sup \{  \rho\circ \varphi (p): \ p\in \widetilde{\Sigma}_{N,\delta_1,C_1}, \ \rho(p)<\widetilde{K}_2  \}$.
 	Up to taking $C_2$ smaller, similar calculations follow to verify that $\varphi^\ast(g_{2,-})<0$ in the points that fulfill $g_{1,-}=0$, for $\rho$ smaller than $\widetilde K_3$ and again, we find an upper bound $K_3$ for $\delta_2$. Then, we take some $\delta_2<\min \{ K_1,K_2, K_3  \}$ and finish the proof.
 \end{proof}
\strut
Denote for simplicity $\mathcal{D}^{nc,fix}=\{ \g_1,\ldots, \g_r  \}$ and $\mathcal{D}^{nc,no-fix}=\{  \g_{r+1} ,\ldots, \g_s \}$ and let $I_j\in \mathcal{I}$ be defined by $\g_j=\g_{I_j}$. Let $J_j\in \mathcal J$ be the index of the chart $(C_{J_j}, \theta,z^{(I_j)},\rho^{(I_j)})$ as in the beginning of section~\ref{sec:noncorner}. Denote also $N_j=N_{\ell,I_j}$. 

Consider the sequence of adimissible blowing-ups $\mathcal{ M}'=(M',\pi',\mathcal{A}',\mathcal{D}')$ over $\mathcal{M}$ constructed as follows. For each $j=1,\dots, s$ let $\tau_j$ be the composition
$$ \tau_j=\sigma_{j,1}\circ \cdots \circ \sigma_{j,N_j},  $$
where $\sigma_{j,1}$ is the admissible blowing-up whose center is the characteristic cycle $\g_j=\g_{I_j}$, $\sigma_{j,2}$ is the admissible blowing-up whose center is the non-corner characteristic cycle $\g_{I_j,1}$ emerging form $\g_{I_j}$ (c.f. Remark~\ref{rmk:blupsimplesing}), and so on. Then, $\mathcal M'$ is the resulting sequence of admissible blowing-ups by setting $\pi'=\pi \circ \tau_1 \circ \cdots \circ \tau_s$.

Notice that $\mathcal M'$ is an adapted reduction of singularities of $\hat \xi$ with same number of non-corner characteristic cycles as $\mathcal M$. We put ${\mathcal{D}'}^{nc}=\{ \gamma_1',\cdots, \gamma_s'  \}$ where $\gamma_j'$ emerges from $\g_j$ by the composition of $\tau_j$ for any $j$. Now, we can take $\ell'\in \mathbb{N}$ satisfying $\ell'\geq \max \{ \ell,  \ell_{\mathcal{ M}'}+1 \}$.
\begin{proposition}
	The vector field $\xi_{\ell'}$ satisfies Theorem~\ref{th:main}.
\end{proposition}
\begin{proof}
	Choose $\delta>0$ sufficiently small so that the conic neighborhoods $V_j:=\wt{\Sigma}_{N_j,\delta,1}(S_{\ell,I_j})$ satisfy, for any $j=1,\ldots, s$
	\begin{itemize}
		\item The Poincaré map $P_{\ell,I_j}$ is defined in $V_j\cap \{ \theta=0  \}={\Sigma}_{N_j,\delta}(\Gamma_{\ell,I_j})$ and satisfies there the conclusions of Lemma~\ref{lem:curvenotfixed} or Lemma~\ref{lem:curveoffixed}, correspondingly.
		\item If $Z$ is a cycle of the transform $\xi_{\ell}^{(J_j)}$ contained in $V_j$, it intersects $ \{ \theta=0  \}$.
		\item If $j\in \{  1,\ldots,r  \}$ (that is, if $\gamma_{I_j}\in \mathcal{D}^{nc,fix}$) then $\Gamma_{\ell,I_j}\subset \fix(P_{\ell,I_j})$ admits a representative in $V_j\cap \{ \theta=0 \}$ as a connected analytic regular curve and such that, for any $a\in \Gamma_{\ell,I_j}\cap V_j$, the cycle of $\xi_{\ell}^{(J_j)}$ through $a$ is contained in $V_j$.
	\end{itemize}
The existence of such $\delta>0$ is guaranteed by standard arguments using the continuity of the flow of $\xi_{\ell}^{(J_j)}$ and the fact that each $\g_j$ is a cycle of $\xi_{\ell}^{(J_j)}$.

Define for $j=1,\ldots,r$ the set $\widetilde{S}_j$ given by the saturation of $\Gamma_{\ell,I_j}$ by the flow of $\xi_{\ell}^{(J_j)}$. By the above properties, $\wt{S}_j$ is an analytic submanifold of $V_j\subset M$
, intersecting the divisor $\pi^{-1}(0)$ along $\g_j$ and completely filled up with cycles of $\xi_{\ell}^{(J_j)}$. We have, furthermore:
\begin{equation}\label{eq:unionCu}
	\mathcal{C}_{ \bigcup_{j=1}^s V_j} (\pi^\ast \xi_{\ell}) =\wt{S}_1\cup \cdots \cup \wt{S}_r  .
\end{equation}
Since $\ell'\geq \ell\geq \ell_{\mathcal{ M}}+1$, we are in the conditions of Corollary~\ref{cor:psiellell'} and there is an analytic conjugation between $\xi_{\ell}^{(J_j)}$ and $\xi_{\ell'}^{(J_j)}$ in a neighborhood of $\g_j$. Moreover, $\psi_{\ell,\ell'}^{(I_j)}$ is in the conditions of Lemma~\ref{lem:cones} with respect to the coordinates $(\theta, z^{(I_j)},\rho^{(I_j)})$ and satisfies $\psi_{\ell,\ell'}^{(I_j)}(S_{\ell,I_j})=S_{\ell',I_j}$.

Thus, applying Lemma~\ref{lem:cones}, we obtain conic neighborhoods $W_j:=\wt{\Sigma}_{N_j,\delta_j, C_j}(S_{\ell',I_j})$ such that $(\psi_{\ell,\ell'}^{(I_j)})^{-1}(W_j)\subset V_j$. Up to shrinking $\delta$, we may assume that $\psi_{\ell,\ell'}^{(J_j)}$ is well defined in $V_j$ for any $j$. We obtain from~\eqref{eq:unionCu} and by conjugation
\begin{equation}
	\label{eq:unionCu2}
		\mathcal{C}_{ \bigcup_{j=1}^s W_j} (\pi^\ast \xi_{\ell'}) = \wt{S}_1 '\cup \cdots \cup \wt{S}_r'  ,
\end{equation}
where $\wt{S}_j:=\psi_{\ell,\ell'}^{(J_j)}(S_j)\cap W_j  $ for $j=1,\ldots, r$.
By Remark~\ref{rmk:blupsimplesing} and taking into account that $\ell\geq l + N_j +1 $ for any $j$, we have that $W_j':=\tau_j^{-1}(W_j)$, together with its closure along the divisor, $\tau_j^{-1}(\g_j)$, is a neighborhood of $\g_j'$ in ${M}'$. We may assume that $W_j'\cap W'_k=\emptyset$ if $j\neq k$. Complete the union $\displaystyle \bigcup_{j=1}^s W_j'$ to a neighborhood $W'$ of $\text{Supp}(\mathcal{D}')$ in $M'$ using Propositions~\ref{prop:zaxis}~\and~\ref{prop:corner} to find two by two disjoint neighborhoods of $\g\in \mathcal{D'}\setminus {\mathcal{D}'}^{nc}$.
We finally apply Proposition~\ref{thm:accumulationlocus} to $W'$ (recall $\ell'\geq \ell_{\mathcal{ M}'}+1$): there is a neighborhood $U$ of $0\in \mathbb{R}^3$ such that
\begin{equation*}
	(\pi')^{-1}(\cu(\xi_{\ell'})) \subset \displaystyle \bigcup_{j=1}^s W_j'.
\end{equation*}
This equation, together with~\eqref{eq:unionCu2} shows that, if we put $S_j:=\pi (\wt{S_j'})\cap U$ for $j=1,\ldots,r$, then
\begin{equation}\label{eq:unionsurfacesfinal}
	\cu(\xi_{\ell'})\subset S_1\cup \cdots \cup S_r.
\end{equation}
Notice that each $S_j$ is an analytic smooth surface in $U\setminus \{ 0\}$ and that $0\in \bar S_j$. Up to taking a smaller $U$ (for instance, such that $\pi^{-1}(U)\cap C_{J_j}\subset \{ \rho^{(I_j)}<\delta_j  \}$ for any $j$), we may assume that $\bar S_j\cap U=S_j\cup \{ 0 \}$. Moreover, $S_j$ is a subanalytic set, since $\pi$ is proper and $\wt{S}_j'$ is semi-analaytic. Finally, by construction, we have that $S_j\subset \pi (\psi_{\ell,\ell'}^{(J_j)}(\wt{S}_j))$ for $j=1,\ldots,r$, there two sets have the same germ at $0\in \mathbb{R}^3$ and the later is entirely composed of cycles for every $j=1,\ldots,r$. This implies, together with~\eqref{eq:unionsurfacesfinal} that if $V\subset U$ is any open neighborhood of $0\in \Rtres$ such that $\text{Fr}(V)\cap S_j$ coincides with one of such cycles for every $j=1,\ldots,r$, then we have
$$ \mathcal{C}_V(\xi_{\ell'})=(S_1\cup \cdots \cup S_r) \cap V .$$
This ends the proof.
\end{proof}

\end{document}